\documentclass[preprint,3p,12pt]{elsarticle}

\usepackage{mathrsfs}
\usepackage{amsmath}
\usepackage{stmaryrd}
\usepackage{bbding}
\usepackage{dcolumn}
\usepackage{graphicx}
\usepackage{amsfonts}
\usepackage{amssymb}
\usepackage{psfrag}
\usepackage{wrapfig}
\usepackage{makeidx}
\usepackage{bm}
\usepackage{epsf}
\usepackage{epsfig}
\usepackage{setspace}
\usepackage{epstopdf}
\usepackage{color}
\usepackage{comment}
\usepackage{algorithm}
\usepackage{algorithmic}
\usepackage{enumitem}
\usepackage{subcaption}
\usepackage{booktabs}
\usepackage{xcolor}
\newcommand \td {\mathrm{~d}}

\journal{}

\begin{document}

\begin{frontmatter}



\title{An asymptotic-preserving adjoint unified gas kinetic scheme for sensitivity analysis}

 \author[HKUST1]{Yue Zhang}
 \ead{yzhangnl@connect.ust.hk}	
  \author[HKUST1]{Junzhe Cao} 
 \ead{jcaobb@connect.ust.hk}
 \author[HKUST1]{Wenpei Long} 
 \ead{wlongab@connect.ust.hk}
 \author[HKUST1,HKUST2,HKUST3]{Kun Xu\corref{cor1}}
 \ead{makxu@ust.hk}

 \address[HKUST1]{Department of Mathematics, Hong Kong University of Science and Technology, Clear Water Bay, Kowloon, Hong Kong}
 \address[HKUST2]{Department of Mechanical and Aerospace Engineering, Hong Kong University of Science and Technology, Clear Water Bay, Kowloon, Hong Kong}
 \address[HKUST3]{Shenzhen Research Institute, Hong Kong University of Science and Technology, Shenzhen, China}
 \cortext[cor1]{Corresponding author}

\begin{abstract}
High-dimensional sensitivity analysis and uncertainty quantification for multiscale gas dynamics, spanning the continuum to rarefied regimes, require computationally efficient and mathematically consistent gradient evaluation. This paper develops a discrete adjoint method for the unified gas-kinetic scheme (UGKS) based on a dual-consistent formulation. The adjoint system is derived directly from the discrete microscopic velocity-distribution equation coupled with the macroscopic-moment compatibility conditions. To resolve the stiff cross-scale coupling, we propose an asymptotic-preserving (AP) adjoint formulation constructed via macroscopic-moment projection and microscopic lifting. Under this framework, the AP adjoint formulation eliminates the stiff collision coupling and removes the collision-time step restriction in the continuum regime. Numerically, a memory-efficient residual-evaluation algorithm that mirrors the forward UGKS cell-vertex data structure is implemented to bypass the memory bottleneck in velocity space. Furthermore, a macroscopic--microscopic predictor--corrector implicit marching scheme is designed to accelerate convergence without solving a globally coupled system. The accuracy, consistency, and robustness of the proposed AP-adjoint scheme are rigorously verified against an independent linearized UGKS solver across a wide range of Knudsen numbers, including lid-driven cavity heat conduction, microchannel thermal creep flow, and hypersonic flow past a circular cylinder.
\end{abstract}


\end{frontmatter}


\section{Introduction}

In the seminal NASA ``CFD Vision 2030'' study~\cite{slotnick2014cfd}, the development of robust, scalable, and mathematically rigorous techniques for high-dimensional sensitivity analysis (SA) and uncertainty quantification (UQ) was highlighted as a critical milestone for next-generation computational frameworks. Traditional computational fluid dynamics (CFD) has primarily focused on deterministic, single-point predictions. Complex engineering systems, especially those under extreme aerothermodynamic conditions, are nevertheless governed by numerous uncertain parameters such as manufacturing tolerances, material properties, and operational boundaries. Propagating these high-dimensional uncertainties through advanced solvers and performing grand-challenge design optimization, therefore demand a paradigm shift from single-point flow evaluations to gradient-aware computational methodologies~\cite{martins2021engineering}.

This need is particularly pronounced in rarefied gas dynamics and non-equilibrium flows, where the physics spans multiple scales and continuum assumptions break down. Existing approaches fall mainly into two categories. Stochastic methods, most notably the direct simulation Monte Carlo (DSMC) method~\cite{bird1998molecular}, track the microscopic distribution with representative particles. DSMC is efficient and robust for high-speed rarefied flows because free streaming and collisions conserve mass, momentum, and energy. Still, the inherent statistical noise makes low-speed simulations expensive, requiring either a large number of particles or extensive averaging. Deterministic discrete velocity methods (DVMs)~\cite{chu1965kinetic,yang1995rarefied}, by contrast, solve the Boltzmann equation or related kinetic models on a discrete velocity grid and yield noise-free solutions; the global velocity-space discretization, however, incurs substantial memory and cost, especially in three-dimensional high-speed applications.

A further limitation shared by classical stochastic and deterministic solvers is operator splitting between free transport and collisions, which introduces time-step-dependent numerical dissipation and constrains the mesh size and time step by the mean free path and mean collision time. Unified gas-kinetic methods overcome this difficulty by coupling collisions and transport within the numerical flux. Representative examples include the DVM-based unified gas-kinetic scheme (UGKS)~\cite{xu2010unified,juan-chen_huang_unified_2012}, the discrete UGKS (DUGKS)~\cite{guo2013dugks}, and the particle-based unified gas-kinetic wave-particle (UGKWP) method~\cite{liu2020ugkwp,zhu2019ugkwp}. Using the integral or finite-difference solution of the kinetic equation along characteristics, these schemes treat the time step as a physical observation scale in the flux evaluation, so that the physical laws are modeled directly in the discrete space~\cite{xu2015direct}. Their unified-preserving (UP) properties have been rigorously established~\cite{guo2023unified}: in the rarefied regime they behave as accurate Boltzmann solvers, while in the continuum regime they recover the Navier--Stokes equations without microscopic scale constraints.

Despite this fidelity, forward multiscale simulations still provide only single-point flow states. For high-dimensional SA, UQ, or design optimization, finite-difference or forward-sensitivity gradients scale linearly with the number of parameters and quickly become intractable. The adjoint method removes this bottleneck: the exact gradient of an objective functional with respect to an arbitrary number of design or uncertain variables can be obtained at a cost comparable to a single forward solve, independent of the parameter dimension~\cite{jameson1988aerodynamic,giles2000introduction,plessix2006review}. This property makes adjoint formulations especially powerful for high-dimensional engineering design~\cite{reuther1999constrained1,reuther1999constrained2}. According to the stage at which differentiation is performed, one distinguishes continuous and discrete adjoint approaches. The continuous framework derives analytic adjoint PDEs before discretization~\cite{jameson1988aerodynamic}, offering physical insight but often requiring difficult manual linearization for complex models. The discrete framework linearizes the discrete residual directly, guaranteeing consistency with the primal solver and facilitating implementation, including through automatic differentiation~\cite{griewank2008evaluating}. Within continuum kinetic CFD, Wu et al.~\cite{wu2026adjoint} recently developed a discrete adjoint gas-kinetic scheme (GKS) for turbulent aerodynamic shape optimization, demonstrating that discrete kinetic adjoint formulations can handle complex continuum regimes.

Parallel progress has been made for rarefied and multiscale flows, with particular emphasis on topology and shape optimization. Caflisch et al.~\cite{caflisch2021adjoint} established the pioneering discretize-then-optimize adjoint DSMC framework for Maxwell molecules, which was subsequently generalized by Yang et al.~\cite{yang2023adjoint} to arbitrary collision kernels by deriving a score-function-based adjoint formulation that bypasses the non-differentiability of rejection sampling. Sato et al.~\cite{sato2019topology} proposed a DVM-based adjoint method for topology optimization of rarefied gas flows. Guan et al.~\cite{guan2024topology} developed an implicit fully coupled discrete adjoint DVM with wall source terms that simplify the adjoint boundary treatment. They also formulated a density-based adjoint method for information-preserving DSMC (IP-DSMC)~\cite{guan2023topology}; because DSMC is inherently unsteady, adjoint optimization requires a full backward evolution and therefore prohibitive memory---approximately $75\,\mathrm{GB}$ even for a simple two-dimensional bend-pipe problem. Building on the discrete adjoint DVM framework, Yuan et al.~\cite{ruifeng2026multi} addressed drag-reduction optimization of very low Earth orbit (VLEO) satellites. For multiscale kinetic schemes, Yuan et al.\ constructed continuous adjoint methods for topology optimization~\cite{yuan2024design} and shape optimization~\cite{yuan2025adjoint} based on DUGKS. Combined with a macroscopic-predictor multiscale solver and the implicit marching of Zhu et al.~\cite{zhu2016implicit}, the approach markedly accelerated adjoint convergence and was subsequently applied to manifold optimization in continuum and rarefied flows~\cite{yuan2026wetted}. The convergence rate was further improved by incorporating the general synthetic iterative scheme (GSIS)~\cite{zhang2026fast}.

The primary objective of this work is to establish a rigorous and computationally efficient discrete adjoint framework for the unified gas-kinetic scheme (UGKS) that remains consistent across multiple scales. To resolve the stiff coupling between the microscopic velocity distribution function and macroscopic moments, a novel asymptotic-preserving (AP)\cite{larsen1983numerical,larsen1987asymptotic,jin1995runge} adjoint formulation is proposed. By constructing a macroscopic-moment projection and a corresponding microscopic lifting operator, it is shown that the discrete adjoint system naturally degenerates to the classical Euler adjoint equations in the continuum limit ($Kn \to 0$) and recovers the discrete-velocity-method (DVM) adjoint equations in the rarefied limit ($Kn \to \infty$).
To overcome the memory bottleneck inherent in high-dimensional phase space, a memory-efficient residual-evaluation algorithm is designed that avoids storing the global Jacobian matrix. Additionally, a macroscopic--microscopic predictor--corrector implicit marching scheme is implemented to accelerate convergence. The accuracy, convergence properties, and computational efficiency of the proposed method are systematically validated through several numerical benchmarks across different flow regimes, demonstrating its capability for gradient-based sensitivity analysis in multiscale gas dynamics.

The remainder of this paper is organized as follows.
Section~2 first formulates sensitivity evaluation via the adjoint method, and then constructs the asymptotic-preserving discrete adjoint UGKS.
Section~3 describes the numerical implementation, with an emphasis on a memory-efficient residual evaluation procedure and an implicit discretization.
Section~4 reports numerical test cases.
Section~5 concludes the paper.
\section{An asymptotic-preserving adjoint method for the unified gas-kinetic scheme}

\subsection{The unified gas-kinetic scheme}
In this work, the Bhatnagar-Gross-Krook(BGK) relaxation model is used, which is a simplified model of the Boltzmann equation. The BGK equation is given by

\begin{equation}
  \label{BGK}
  \frac{\partial f}{\partial t} + \mathbf{v}\cdot\nabla f =\frac{f^+-f}{\tau},
\end{equation}
where $f=f(\mathbf{x},t,\mathbf{v},\mathbf{\xi})$ is the distribution function for gas molecules at physical space $\mathbf{x}$ and time $t$. Here $\mathbf{v}$ is the particle velocity, $\mathbf{\xi}$ is the velocity of the internal degrees of freedom, $\tau$ is the particle collision time, and $f^+$ is the modified equilibrium distribution function. The modified equilibrium distribution function is given by the Shakhov model

\begin{equation*}
 f^+(\mathbf{W},\mathbf{q})=g\left[ 1+(1-\text{Pr})2\lambda\mathbf{c}\cdot\mathbf{q}\left(2\lambda c^2-5\right)/(5p)\right]=g+g^+,
\end{equation*}
where $g$ is the Maxwellian distribution, $\text{Pr}$ is the Prandtl number, $\mathbf{c}=\mathbf{v}-\mathbf{V}$ is the peculiar velocity, $\mathbf{V}$ is the macroscopic velocity, $\mathbf{q}$ is the heat flux, $\lambda=m/(2k_B T)$, $m$ is molecule mass, $k_B$ is Boltzmann constant, $T$ is the temperature. The Maxwellian distribution is

\begin{equation*}
 g=\rho {\left(\frac{\lambda}{\pi}\right)}^\frac{K+D}{2}e^{-\lambda ((\mathbf{v}-\mathbf{V})^2+\mathbf{\xi}^2)},
\end{equation*}
where $D$ is the spatial dimension and $K$ is the inner degree of freedom. The heat flux is obtained through the distribution function
\begin{equation*}
  \mathbf{q}=\frac{1}{2}\int (\mathbf{v}-\mathbf{V})(|\mathbf{v}-\mathbf{V}|^2+\mathbf{\xi}^2)f\td\mathbf{v} \td \mathbf{\xi},
\end{equation*}
where $\mathbf{\xi}^2=\xi_1^2+\xi_2^2+\cdots+\xi_K^2$.
The collision term satisfies the compatibility condition

\begin{equation*}
 \int (f^+-f)\mathbf{\Psi}\td \mathbf{v}\td\mathbf{\xi}=0,
\end{equation*}
or equivalently in the discrete velocity formulation,
\begin{equation}
  \label{eq:BGK_compat}
  \sum_{k} w_k (f_{i,k}^+-f_{i,k})\mathbf{\Psi}_k=0,
\end{equation}
where the velocity space is divided into discrete points $\delta \mathcal{V}_k$, the index $k$ denotes the discrete points in the particle velocity space, and $w_k$ is the weight of the discrete point.

The finite volume method is used in this paper. For each cell $\Omega_i$, the governing equation of the distribution function is
\begin{equation}
  \label{UGKSmicro}
 \frac{f_{i,k}^{n+1}-f_{i,k}^{n}}{\Delta t} + \frac{1}{|\Omega_i|}\underbrace{  \frac{1}{\Delta t}\sum_{j\in N(i)} S_{ij}\int_{t^n}^{t^{n+1}}\mathbf{v}_k\cdot \mathbf{n}_{ij}f_{ij,k}\td t}_{\mathcal{R}_{i,k}^n}=\frac{f_{i,k}^{+,n+1}-f_{i,k}^{n+1}}{\tau},
\end{equation}
where $|\Omega_i|$ denotes the volume of cell $i$, $N(i)$ is the set of all interface-adjacent neighboring cells of cell $i$ and $j$ is one of the neighboring cells of $i$. The interface between cells $i$ and $j$ is labeled $ij$ and has area $S_{ij}$. $\mathbf{v}_k$ is the velocity of the $k$-th discrete point. $\mathbf{n}_{ij}$ is the normal vector of interface $ij$.
To construct the numerical flux at cell interface $\mathbf{x}_0=(0,0,0)^T$, the integral solution along the characteristic line $\mathbf{x}^\prime =\mathbf{x}_0-\mathbf{v}(t-t^\prime)$ of the BGK equation (\ref{BGK}) gives
\begin{equation}
  \label{integralSolution}
 f(\mathbf{x},t,\mathbf{\xi})=\frac{1}{\tau}\int_{t_0}^t f^+(\mathbf{x}^\prime,t^\prime)e^{-(t-t^\prime)/\tau}\td t^\prime + e^{-t/\tau}f_{0}(\mathbf{x}-\mathbf{v} t),
\end{equation}
where $f_{0}(\mathbf{x})$ is the initial gas distribution function at the beginning of each step $t_n$, and $f^+(\mathbf{x},t)$ is the effective equilibrium state distributed in space and time. The integral solution provides a multiscale model of the evolution from an initial non-equilibrium distribution $f$ to an equilibrium distribution $f^+$ via collisions. To achieve second-order accuracy, the initial distribution function $f_{0}(\mathbf{x})$ is approximated as

\begin{equation*}
 f_{0}(\mathbf{x}) = \begin{cases}
 f^l+\mathbf{x}\cdot \nabla f^l,&v_n>0,\\
 f^r+\mathbf{x}\cdot \nabla  f^r,&v_n<0,\\
 \end{cases}
\end{equation*}
where $f^l$ and $f^r$ are the reconstructed initial distribution functions at the left and right sides of the interface. The equilibrium state is approximated as

\begin{equation*}
 f^+(\mathbf{x},t) \approx g_{0}+ g^+_{0} + \mathbf{x}\cdot\nabla g_{0}+\frac{\partial g_{0}}{\partial t}t,
\end{equation*}
where $g_{0}$ is the Maxwellian distribution function obtained from the conservative flow variables of colliding particles from both sides of the interface.

One thing should be noted that actually the reduced distribution function $h_{i,k}^{n}=\int f_{i,k}^{n}\td\mathbf{\xi}$ and $b_{i,k}^{n}=\int \mathbf{\xi}^2f_{i,k}^{n}\td\mathbf{\xi}$ are used in the UGKS. The relation
$f_{i,k}^{n}\mathbf{\Psi}_k=[h_{i,k}^{n},\mathbf{v}_k h_{i,k}^{n},1/2(\mathbf{v}_k^2 h_{i,k}^{n}+b_{i,k}^{n})]^T$ is used. However, for the simplicity of the derivation, the full distribution function $f_{i,k}^{n}$ is used in this paper.

By taking moments of the microscopic governing equation Eq.~(\ref{UGKSmicro}) and using the compatibility condition of the collision term, the macroscopic governing equation of the conservative variables is given by
\begin{equation}
  \label{macroscopicUpdata}
  \frac{\mathbf{W}_i^{n+1}-\mathbf{W}_i^n}{\Delta t} + \frac{1}{|\Omega_i|}\sum_k w_k \mathcal{R}_{i,k}^n\mathbf{\Psi}_k=0.
\end{equation}
More details are demonstrated in earlier work \cite{xu2010unified,xu2015direct}.

\subsection{Gradient Evaluation via the Adjoint Formulation}

With the primal UGKS discretization established, we first formulate the sensitivity problem through a Lagrangian, then identify the discrete residual constraints and the associated dual framework; the explicit expressions of the adjoint operators are derived in the following subsections.

The primary objective of the adjoint method is to efficiently evaluate sensitivities. Let $J$ denote a scalar objective function (e.g., drag, lift, or heat flux), and let $\mathbf{p}$ represent the vector of dependent parameters, which typically define boundary shapes, flow characteristics, or governing source terms. The goal is to determine the total derivative of $J$ with respect to $\mathbf{p}$ under the governing constraints $\mathcal{L}_{\text{micro}}=0$ and $\mathcal{L}_{\text{macro}}=0$. To this end, the Lagrange multiplier method is employed: the macroscopic and microscopic Lagrange multipliers $\mathbf{\Lambda}_i^{*,n}$ and $\lambda_{i,k}^{*,n}$ are introduced, and the Lagrangian is formed as $\mathcal{L}=J+\mathcal{I}$, where $\mathcal{I}$ is the full time-space inner product defined below.

The macroscopic governing equation~\eqref{macroscopicUpdata} can be recovered by taking moments of the microscopic equation~\eqref{UGKSmicro} together with the collision compatibility condition~\eqref{eq:BGK_compat}.
Accordingly, the discrete constraints employed in the adjoint formulation are the macroscopic compatibility residual and the microscopic UGKS residual:
\begin{equation}
  \label{eq:L_macro}
  \mathcal{L}_{\text{macro}}(\mathbf{W}_i^{n+1},f_{i,k}^{n+1})=\sum_k w_k(f_{i,k}^{+,n+1}- {f}_{i,k}^{n+1})\mathbf{\Psi}_k = 0,\quad f_{i,k}^{+,n+1}=f_{i,k}^{+,n+1}(\mathbf{W}_i^{n+1},\mathbf{q}_i^n,\mathbf{v}_k).
\end{equation}
\begin{equation}
  \label{eq:L_micro}
  \begin{aligned}
 \mathcal{L}_{\text{micro}}(\mathbf{W}_i^{n},f_{i,k}^{n})=\frac{f_{i,k}^{n+1}-f_{i,k}^{n}}{\Delta t}+\frac{1}{|\Omega_i|}\mathcal{R}_{i,k}^n - \frac{ f_{i,k}^{+,n+1}-f_{i,k}^{n+1}}{\tau_i^{n+1}}= 0.
  \end{aligned}
\end{equation}
An alternative adjoint derivation that starts directly from the primal macroscopic update~\eqref{macroscopicUpdata} and the microscopic equation~\eqref{UGKSmicro} is given in~\ref{app:ill_conditioned}; that route is less natural for the subsequent asymptotic-preserving reconstruction, but the two approaches lead to equivalent results.

The full time-space inner product is defined as
\begin{equation*}
  \mathcal{I}=\langle \mathcal{L}_{\text{macro}}(\mathbf{W}_i^{n},f_{i,k}^{n}), \mathbf{\Lambda}_{i}^{*,n} \rangle_{\text{macro}} + \langle \mathcal{L}_{\text{micro}}(\mathbf{W}_i^{n},f_{i,k}^{n}), {\lambda}_{i,k}^{*,n} \rangle_{\text{micro}},
\end{equation*}
and the Lagrangian is therefore
\begin{equation*}
  \begin{aligned}
    \mathcal{L}=J+\mathcal{I}=J+\langle \mathcal{L}_{\text{macro}}(\mathbf{W}_i^{n},f_{i,k}^{n}), \mathbf{\Lambda}_i^{*,n} \rangle_{\text{macro}} + \langle \mathcal{L}_{\text{micro}}(\mathbf{W}_i^{n},f_{i,k}^{n}), \lambda_{i,k}^{*,n} \rangle_{\text{micro}}.
  \end{aligned}
\end{equation*}
The full time-space inner products in the macroscopic and microscopic spaces are defined as
\begin{equation*}
  \langle \mathbf{\Phi},\mathbf{\Phi^*} \rangle_{\text{macro}}= \sum_{n=1}^{N} \sum_{i} (\mathbf{\Phi}_i^n)^T \mathbf{\Phi}_i^{*,n}\cdot |\Omega_i|,\quad \langle \phi,\phi^* \rangle_{\text{micro}}= \sum_{n=1}^{N-1} \sum_{i} \sum_{k} \phi_{i,k}^n \phi_{i,k}^{*,n}\cdot w_k\cdot |\Omega_i|,
\end{equation*}
where $N$ denotes the total number of time steps, while $\mathbf{\Phi},\mathbf{\Phi^*},\phi,\phi^*$ represent arbitrary macroscopic and microscopic variables, respectively.
The macroscopic and microscopic adjoint operators $\mathcal{L}_{\text{macro}}^*$ and $\mathcal{L}_{\text{micro}}^*$ are defined by transferring the primal residuals onto the dual variables through the duality relation
\begin{equation*}
  \begin{aligned}
  &\langle \mathcal{L}_{\text{macro}}(\mathbf{W}_i^{n},f_{i,k}^{n}), \mathbf{\Lambda}_i^{*,n} \rangle_{\text{macro}} + \langle \mathcal{L}_{\text{micro}}(\mathbf{W}_i^{n},f_{i,k}^{n}), \lambda_{i,k}^{*,n} \rangle_{\text{micro}} \\
  &\qquad=
  \langle \mathbf{W}_i^{n}, \mathcal{L}_{\text{macro}}^*(\mathbf{\Lambda}_i^{*,n},\lambda_{i,k}^{*,n}) \rangle_{\text{macro}} + \langle f_{i,k}^{n}, \mathcal{L}_{\text{micro}}^*(\mathbf{\Lambda}_i^{*,n},\lambda_{i,k}^{*,n}) \rangle_{\text{micro}},
  \end{aligned}
\end{equation*}
so that the full time-space inner product admits the equivalent representation
\begin{equation*}
  \mathcal{I}=\langle \mathbf{W}_i^{n}, \mathcal{L}_{\text{macro}}^*(\mathbf{\Lambda}_i^{*,n},\lambda_{i,k}^{*,n}) \rangle_{\text{macro}} + \langle f_{i,k}^{n}, \mathcal{L}_{\text{micro}}^*(\mathbf{\Lambda}_i^{*,n},\lambda_{i,k}^{*,n}) \rangle_{\text{micro}}.
\end{equation*}

With the Lagrangian and the duality relation in place, the total variation of $\mathcal{L}$ is computed by applying the chain rule to the objective function $J$ and the respective inner products. The variation operator is expanded with respect to the dependent parameters $\mathbf{p}$ and the state variables $\mathbf{W}_i^n, f_{i,k}^n$:
\begin{equation*}
  \begin{aligned}
  \delta \mathcal{L} =& \frac{\partial J}{\partial \mathbf{p}} \delta \mathbf{p}+\left\langle\delta \mathbf{W}_i^n,\frac{\partial J}{\partial \mathbf{W}_i^n}/|\Omega_i| \right\rangle_{\text{macro}}+\left\langle\delta f_{i,k}^n,\frac{\partial J}{\partial f_{i,k}^n}/(w_k|\Omega_i|) \right\rangle_{\text{micro}}\\
    &+\left(\delta\mathbf{p}\frac{\partial}{\partial \mathbf{p}}+\delta\mathbf{W}_i^n\frac{\partial}{\partial \mathbf{W}_i^n} +\delta f_{i,k}^n\frac{\partial}{\partial f_{i,k}^n}\right) (\langle \mathcal{L}_{\text{macro}}, \mathbf{\Lambda}_i^{*,n} \rangle_{\text{macro}} + \langle \mathcal{L}_{\text{micro}}, \lambda_{i,k}^{*,n} \rangle_{\text{micro}}).
  \end{aligned}
\end{equation*}
Next, the variations of the dependent parameters $\mathbf{p}$ are isolated from those of the state variables. For the variations associated with the states $\mathbf{W}_i^n$ and $f_{i,k}^n$, the duality relation is used to transfer the linear operations onto the adjoint variables:
\begin{equation*}
  \begin{aligned}
  \delta \mathcal{L} =& \frac{\partial J}{\partial \mathbf{p}} \delta \mathbf{p}+\left\langle\delta \mathbf{W}_i^n,\frac{\partial J}{\partial \mathbf{W}_i^n}/|\Omega_i| \right\rangle_{\text{macro}}+\left\langle\delta f_{i,k}^n,\frac{\partial J}{\partial f_{i,k}^n}/(w_k|\Omega_i|)\right\rangle_{\text{micro}}\\
    &+\delta\mathbf{p}\frac{\partial}{\partial \mathbf{p}} (\langle \mathcal{L}_{\text{macro}}, \mathbf{\Lambda}_i^{*,n} \rangle_{\text{macro}} + \langle \mathcal{L}_{\text{micro}}, \lambda_{i,k}^{*,n} \rangle_{\text{micro}}) \\
    &+\left(\delta\mathbf{W}_i^n\frac{\partial}{\partial \mathbf{W}_i^n} +\delta f_{i,k}^n\frac{\partial}{\partial f_{i,k}^n}\right) (\langle \mathbf{W}_i^{n}, \mathcal{L}_{\text{macro}}^*\rangle_{\text{macro}} + \langle f_{i,k}^{n}, \mathcal{L}_{\text{micro}}^*\rangle_{\text{micro}}).
  \end{aligned}
\end{equation*}
Evaluating the partial derivatives and explicitly moving the variations into the inner products simplifies the equation to:
\begin{equation*}
  \begin{aligned}
  \delta \mathcal{L} =& \frac{\partial J}{\partial \mathbf{p}} \delta \mathbf{p}+\left\langle\delta \mathbf{W}_i^n,\frac{\partial J}{\partial \mathbf{W}_i^n}/|\Omega_i| \right\rangle_{\text{macro}}+\left\langle\delta f_{i,k}^n,\frac{\partial J}{\partial f_{i,k}^n}/(w_k|\Omega_i|)\right\rangle_{\text{micro}}  \\
    &+\delta\mathbf{p} \langle\frac{\partial}{\partial \mathbf{p}} \mathcal{L}_{\text{macro}}, \mathbf{\Lambda}_i^{*,n} \rangle_{\text{macro}} + \delta\mathbf{p} \langle \frac{\partial}{\partial \mathbf{p}}\mathcal{L}_{\text{micro}}, \lambda_{i,k}^{*,n} \rangle_{\text{micro}}\\
    &+ \langle \delta \mathbf{W}_i^{n}, \mathcal{L}_{\text{macro}}^*\rangle_{\text{macro}} + \langle \delta f_{i,k}^{n}, \mathcal{L}_{\text{micro}}^*\rangle_{\text{micro}}.
  \end{aligned}
\end{equation*}
Finally, the terms are grouped by the independent variations $\delta\mathbf{p}$, $\delta\mathbf{W}_i^n$, and $\delta f_{i,k}^n$ to isolate the sensitivity gradient and define the stationary conditions for the adjoint system:
\begin{equation*}
  \begin{aligned}
  \delta \mathcal{L} = &\left( \frac{\partial J}{\partial \mathbf{p}} + \langle \frac{\partial \mathcal{L}_{\text{macro}}}{\partial \mathbf{p}}, \mathbf{\Lambda}_i^{*,n} \rangle_{\text{macro}} + \langle \frac{\partial \mathcal{L}_{\text{micro}}}{\partial \mathbf{p}}, \lambda_{i,k}^{*,n} \rangle_{\text{micro}} \right) \delta \mathbf{p}\\
    &+\left\langle\delta \mathbf{W}_i^n,\frac{\partial J}{\partial \mathbf{W}_i^n}/|\Omega_i| + \mathcal{L}_{\text{macro}}^*\right\rangle_{\text{macro}}+\left\langle\delta f_{i,k}^n,\frac{\partial J}{\partial f_{i,k}^n}/(w_k|\Omega_i|)+ \mathcal{L}_{\text{micro}}^*\right\rangle_{\text{micro}}.
  \end{aligned}
\end{equation*}
By eliminating the dependence on the state variations, the total derivative of the objective function $J$ with respect to the dependent parameter vector $\mathbf{p}$ reduces to:
\begin{equation}
  \label{eq:total_gradient}
  \frac{\text{d} J}{\text{d} \mathbf{p}}=\frac{\partial J}{\partial \mathbf{p}}+ \langle \frac{\partial \mathcal{L}_{\text{micro}}}{\partial \mathbf{p}}, \lambda_{i,k}^{*,n} \rangle_{\text{micro}},
\end{equation}
with ${\partial \mathcal{L}_{\text{macro}}}/{\partial \mathbf{p}}=0$.
The corresponding adjoint system:
\begin{equation}
  \label{eq:adjoint_system_coupled}
  \begin{cases}
    \mathcal{L}_{\text{macro}}^*(\mathbf{\Lambda}_i^{*,n},\lambda_{i,k}^{*,n}) + \frac{\partial J}{\partial \mathbf{W}_i^n}/|\Omega_i| = 0, \\
    \mathcal{L}_{\text{micro}}^*(\mathbf{\Lambda}_i^{*,n},\lambda_{i,k}^{*,n}) + \frac{\partial J}{\partial f_{i,k}^n}/(w_k|\Omega_i|) = 0.
  \end{cases}
\end{equation}

\subsection{The Adjoint System of the Unified Gas-Kinetic Scheme}
The sensitivity formula~\eqref{eq:total_gradient} is closed once the explicit forms of the adjoint operators $\mathcal{L}_{\text{macro}}^*$ and $\mathcal{L}_{\text{micro}}^*$ in~\eqref{eq:adjoint_system_coupled} are available. This subsection derives those operators from the discrete residuals~\eqref{eq:L_macro} and~\eqref{eq:L_micro}.

By taking the first-order variations of the full time-space inner product with respect to the macroscopic state $\mathbf{W}_i^{n}$ and the microscopic distribution function $f_{i,k}^{n}$, the corresponding adjoint system can be derived. Specifically, as this paper focuses on steady-state problems, all time indices associated with the original variables are omitted in the following derivations.

The microscopic adjoint equation $\mathcal{L}_{\text{micro}}^* = 0$ is obtained as:
\begin{equation}
  \label{scaled_micro}
  \begin{aligned}
  \mathcal{L}_{\text{micro}}^*(\mathbf{\Lambda}_i^{*,n},\lambda_{i,k}^{*,n})=\frac{\lambda_{i,k}^{*,n-1}-\lambda_{i,k}^{*,n}}{\Delta t}&+\frac{1}{|\Omega_i|}\sum_j\frac{\partial \mathcal{R}_{j,k}}{\partial f_{i,k}}\lambda_{j,k}^{*,n}+ \frac{1}{|\Omega_i|} \frac{1}{w_k}\frac{\partial \mathbf{q}_i}{\partial f_{i,k}^{n}}\sum_j \sum_{k^\prime}w_{k^\prime}\frac{\partial \mathcal{R}_{j,k^\prime}^n}{\partial \mathbf{q}_i}\lambda_{j,k^\prime}^{*,n} \\
  &+\frac{\lambda_{i,k}^{*,n-1}}{\tau_i}-\frac{1}{\tau_i}\frac{1}{w_k}\frac{\partial \mathbf{q}_i}{\partial f_{i,k}}\sum_{k^\prime} w_{k^\prime}  \frac{\partial f_{i,k^\prime}^+}{\partial \mathbf{q}_i}  \lambda_{i,k^\prime}^{*,n}  
  -\mathbf{\Psi}_k^T \mathbf{\Lambda}_i^{*,n} = 0.
\end{aligned}
\end{equation}

Similarly, the adjoint compatibility condition $\mathcal{L}_{\text{macro}}^* = 0$ is derived as:
\begin{equation}
  \label{adjointmacro_final}
  \begin{aligned}
   \mathcal{L}_{\text{macro}}^*(\mathbf{\Lambda}_i^{*,n},\lambda_{i,k}^{*,n})=&\mathbf{\Lambda}_i^{*,n} + \frac{1}{|\Omega_i|}\sum_j\sum_{k} w_k\frac{\partial \mathcal{R}_{j,k}}{\partial \mathbf{W}_i}\lambda_{j,k}^{*,n}+ \frac{1}{|\Omega_i|}\frac{\partial \mathbf{q}_i}{\partial \mathbf{W}_i}\sum_j\sum_{k} w_k\frac{\partial \mathcal{R}_{j,k}}{\partial \mathbf{q}_{i}}\lambda_{j,k}^{*,n}\\
   &-\frac{1}{\tau_i}\sum_{k} w_k  \frac{\partial f_{i,k}^+}{\partial \mathbf{W}_i}\lambda_{i,k}^{*,n-1}-\frac{1}{\tau_i}\frac{\partial \mathbf{q}_i}{\partial \mathbf{W}_i}\sum_{k} w_k  \frac{\partial f_{i,k}^+}{\partial \mathbf{q}_i} \lambda_{i,k}^{*,n}\\
   &+\frac{1}{\tau_i}\frac{\partial \tau_i}{\partial \mathbf{W}_i}\sum_{k} w_k\frac{f_{i,k}^+-f_{i,k}}{\tau_i}\lambda_{i,k}^{*,n}= 0.
  \end{aligned}
\end{equation}
In Eq.~\eqref{adjointmacro_final}, the derivative of the scaled macroscopic flux with respect to the macroscopic conservative values $\mathbf{W}_i$ is expanded via the chain rule as follows:
\begin{equation*}
  \frac{\partial \mathcal{R}_{j,k}}{\partial \mathbf{W}_i} =  \frac{\partial \mathcal{R}_{j,k}}{\partial g_{i,k}}\frac{\partial g_{i,k}}{\partial \mathbf{W}_i} + \frac{\partial \mathcal{R}_{j,k}}{\partial \tau_{ij}}\frac{\partial \tau_{ij}}{\partial \mathbf{W}_i} ,
\end{equation*}
where $\tau_{ij}$ is calculated by the equilibrium state of the interface between cell $i$ and $j$.
In the present work, the chain rule is adopted to determine these derivatives; the corresponding partial-derivative expressions are given in ~\ref{app:jacobi}.

\subsection{The Asymptotic-Preserving Adjoint Method}

The coupled adjoint system Eqs.~\eqref{scaled_micro} and~\eqref{adjointmacro_final} contains a stiff feedback through $\partial f_{i,k}^+/\partial \mathbf{W}_i$ and the mixed time-level term $\lambda_{i,k}^{*,n-1}$. All AP reconstructions and asymptotic analyses below are constructed algebraically from these two equations only.

To eliminate the stiff coupling while preserving the original adjoint operators, an intermediate macroscopic adjoint moment is introduced:
\begin{equation*}
  \widetilde{\mathbf{\Lambda}}_i^{*,n}=\sum_{k} w_k \frac{\partial f_{i,k}^+}{\partial \mathbf{W}_i}\lambda_{i,k}^{*,n},
\end{equation*}
which is the coefficient of $\lambda_{i,k}^{*,n-1}$ in the collision block of Eq.~\eqref{adjointmacro_final}.

\textbf{Route I: macroscopic-moment projection.}
Multiplying Eq.~\eqref{scaled_micro} by $w_k\partial f_{i,k}^+/\partial \mathbf{W}_i$ and summing over $k$ gives
\begin{equation} \label{asymptotic-preservingadjointmacro}
  \begin{aligned}
  \sum_{k} w_k \frac{\partial f_{i,k}^+}{\partial \mathbf{W}_i}\frac{\lambda_{i,k}^{*,n-1}-\lambda_{i,k}^{*,n}}{\Delta t} &+ \frac{1}{|\Omega_i|}\sum_{k} w_k\frac{\partial f_{i,k}^+}{\partial \mathbf{W}_i}\sum_j \frac{\partial \mathcal{R}_{j,k}}{\partial f_{i,k}}\lambda_{j,k}^{*,n} \\
  &+ \frac{1}{|\Omega_i|}\sum_{k} \frac{\partial f_{i,k}^+}{\partial \mathbf{W}_i}\frac{\partial \mathbf{q}_i}{\partial f_{i,k}}\sum_j \sum_{k^\prime}w_{k^\prime}\frac{\partial \mathcal{R}_{j,k^\prime}}{\partial \mathbf{q}_i}\lambda_{j,k^\prime}^{*,n} \\
  &+ \frac{1}{\tau_i}\sum_{k} w_k \lambda_{i,k}^{*,n-1} \frac{\partial f_{i,k}^+}{\partial \mathbf{W}_i} - \frac{1}{\tau_i}\sum_{k} \frac{\partial f_{i,k}^+}{\partial \mathbf{W}_i}\frac{\partial \mathbf{q}_i}{\partial f_{i,k}}\sum_{k^\prime} w_{k^\prime}  \frac{\partial f_{i,k^\prime}^+}{\partial \mathbf{q}_i} \lambda_{i,k^\prime}^{*,n} \\
  &- \mathbf{\Lambda}_i^{*,n}=0,
  \end{aligned}
\end{equation}
where the last term follows from
\begin{equation*}
  \sum_{k} w_k \mathbf{\Psi}_k \frac{\partial f_{i,k}^+}{\partial \mathbf{W}_i} = \frac{\partial}{\partial \mathbf{W}_i } \sum_{k} w_k \mathbf{\Psi}_k f_{i,k}^+= \mathbf{I}
\end{equation*}
and the $-\mathbf{\Psi}_k^T \mathbf{\Lambda}_i^{*,n}$ term in Eq.~\eqref{scaled_micro}.
Adding Eq.~\eqref{adjointmacro_final} to Eq.~\eqref{asymptotic-preservingadjointmacro}, the terms proportional to $\partial f_{i,k}^+/\partial \mathbf{W}_i$ and $\lambda_{i,k}^{*,n-1}$ cancel exactly.
The remaining $\mathbf{q}$-dependent blocks share the common coefficient
\begin{equation*}
  \sum_{k} \frac{\partial f_{i,k}^+}{\partial \mathbf{W}_i}\frac{\partial \mathbf{q}_i}{\partial f_{i,k}}+\frac{\partial \mathbf{q}_i}{\partial \mathbf{W}_i}
  =\frac{\mathrm{d}\mathbf{q}_i}{\mathrm{d}\mathbf{W}_i}\bigg|_{f=f^+},
\end{equation*}
which is the total derivative of the heat flux along the equilibrium manifold.
Because the Shakhov construction enforces $\mathbf{q}[f^+(\mathbf{W},\mathbf{q})]=\mathbf{q}$ independently of $\mathbf{W}$, this total derivative vanishes identically, and the corresponding blocks drop out.
The projected macroscopic adjoint equation therefore reduces to
\begin{equation}
  \label{asymptotic-preservingadjointsystem}
  \begin{aligned}
  \sum_{k} w_k  \frac{\partial f_{i,k}^+}{\partial \mathbf{W}_i}\frac{\lambda_{i,k}^{*,n-1}-\lambda_{i,k}^{*,n}}{\Delta t} &+ \frac{1}{|\Omega_i|}\sum_{k} w_k\frac{\partial f_{i,k}^+}{\partial \mathbf{W}_i}\sum_j \frac{\partial \mathcal{R}_{j,k}}{\partial f_{i,k}}\lambda_{j,k}^{*,n}+ \frac{1}{|\Omega_i|}\sum_j\sum_{k} w_k \frac{\partial \mathcal{R}_{j,k}}{\partial \mathbf{W}_i}\lambda_{j,k}^{*,n}\\
  &+\frac{1}{\tau_i}\frac{\partial \tau_i}{\partial \mathbf{W}_i}\sum_{k} w_k\frac{f_{i,k}^+-f_{i,k}}{\tau_i}\lambda_{i,k}^{*,n} =0.
  \end{aligned}
\end{equation}

Substituting $\widetilde{\mathbf{\Lambda}}_i^{*,n}$ and $\widetilde{\mathbf{\Lambda}}_i^{*,n-1}$ into Eq.~\eqref{asymptotic-preservingadjointsystem}, the projected macroscopic adjoint evolution becomes
\begin{equation}
  \label{asymptotic-preservingadjointsystem2}
      \begin{aligned}
  \frac{\widetilde{\mathbf{\Lambda}}_i^{*,n-1}-\widetilde{\mathbf{\Lambda}}_i^{*,n}}{\Delta t}&+ \frac{1}{|\Omega_i|}\sum_{k} w_k\frac{\partial f_{i,k}^+}{\partial \mathbf{W}_i}\sum_j \frac{\partial \mathcal{R}_{j,k}}{\partial f_{i,k}}\lambda_{j,k}^{*,n}+ \frac{1}{|\Omega_i|}\sum_j\sum_{k} w_k \frac{\partial \mathcal{R}_{j,k}}{\partial \mathbf{W}_i}\lambda_{j,k}^{*,n}\\
  &+\frac{1}{\tau_i}\frac{\partial \tau_i}{\partial \mathbf{W}_i}\sum_{k} w_k\frac{f_{i,k}^+-f_{i,k}}{\tau_i}\lambda_{i,k}^{*,n}=0.
  \end{aligned}
\end{equation}

\textbf{Route II: microscopic lifting.}
Multiplying Eq.~\eqref{adjointmacro_final} by $\mathbf{\Psi}_k^T$ and adding the result to Eq.~\eqref{scaled_micro} for each velocity point $k$ gives
\begin{equation}
  \label{eq:micro_coupled_W}
  \begin{aligned}
  \frac{\lambda_{i,k}^{*,n-1}-\lambda_{i,k}^{*,n}}{\Delta t}&+\frac{1}{|\Omega_i|}\sum_j\frac{\partial \mathcal{R}_{j,k}}{\partial f_{i,k}}\lambda_{j,k}^{*,n}+\mathbf{\Psi}_k^T\frac{1}{|\Omega_i|}\sum_j\sum_{k^\prime} w_{k^\prime}\frac{\partial \mathcal{R}_{j,k^\prime}}{\partial \mathbf{W}_i}\lambda_{j,k^\prime}^{*,n} \\
  &+\frac{1}{|\Omega_i|}\left(\frac{1}{w_k}\frac{\partial \mathbf{q}_i}{\partial f_{i,k}}+\mathbf{\Psi}_k^T\frac{\partial \mathbf{q}_i}{\partial \mathbf{W}_i}\right)\sum_j\sum_{k^\prime}w_{k^\prime}\frac{\partial \mathcal{R}_{j,k^\prime}}{\partial \mathbf{q}_i}\lambda_{j,k^\prime}^{*,n}\\
  &-\frac{\mathbf{\Psi}_k^T\widetilde{\mathbf{\Lambda}}_i^{*,n-1}-\lambda_{i,k}^{*,n-1}}{\tau_i}\\
  &-\frac{1}{\tau_i}\left(\frac{1}{w_k}\frac{\partial \mathbf{q}_i}{\partial f_{i,k}}+\mathbf{\Psi}_k^T\frac{\partial \mathbf{q}_i}{\partial \mathbf{W}_i}\right)\sum_{k^\prime} w_{k^\prime}  \frac{\partial f_{i,k^\prime}^+}{\partial \mathbf{q}_i} \lambda_{i,k^\prime}^{*,n}\\
  &+\mathbf{\Psi}_k^T\frac{1}{\tau_i}\frac{\partial \tau_i}{\partial \mathbf{W}_i}\sum_{k^\prime} w_{k^\prime}\frac{f_{i,k^\prime}^+-f_{i,k^\prime}}{\tau_i}\lambda_{i,k^\prime}^{*,n}
  =0.
  \end{aligned}
\end{equation}
The collision term in Eq.~\eqref{eq:micro_coupled_W} is obtained by combining $\lambda_{i,k}^{*,n-1}/\tau_i$ in Eq.~\eqref{scaled_micro} with $\mathbf{\Psi}_k^T$ times the collision block $-\tau_i^{-1}\sum_{k^\prime} w_{k^\prime} (\partial f_{i,k^\prime}^+/\partial \mathbf{W}_i)\lambda_{i,k^\prime}^{n-1}$ in Eq.~\eqref{adjointmacro_final}, and the explicit $\mathbf{\Lambda}_i^{*,n}$ terms cancel through the $-\mathbf{\Psi}_k^T \mathbf{\Lambda}_i^{*,n}$ contribution in Eq.~\eqref{scaled_micro}.

Therefore, the asymptotic-preserving adjoint operators are given by the left-hand sides of the projected macroscopic equation~\eqref{asymptotic-preservingadjointsystem2} and the lifted microscopic equation~\eqref{eq:micro_coupled_W}, denoted here by $\mathcal{L}_{\text{macro}}^{\mathrm{AP}*}$ and $\mathcal{L}_{\text{micro}}^{\mathrm{AP}*}$, respectively.
Applying the same projection and lifting to the inhomogeneous system~\eqref{eq:adjoint_system_coupled}, the objective source terms transform accordingly, and the AP adjoint system takes the form
\begin{equation}
  \label{eq:adjoint_system_AP}
  \begin{cases}
    \displaystyle
    \mathcal{L}_{\text{macro}}^{\mathrm{AP}*}\!\left(\widetilde{\mathbf{\Lambda}}_i^{*,n},\lambda_{i,k}^{*,n}\right)
    +\frac{1}{|\Omega_i|}\left(
      \frac{\partial J}{\partial \mathbf{W}_i^n}
      +\sum_{k}\frac{\partial f_{i,k}^+}{\partial \mathbf{W}_i}\frac{\partial J}{\partial f_{i,k}^n}
    \right)=0, \\[1.2em]
    \displaystyle
    \mathcal{L}_{\text{micro}}^{\mathrm{AP}*}\!\left(\widetilde{\mathbf{\Lambda}}_i^{*,n},\lambda_{i,k}^{*,n}\right)
    +\dfrac{\partial J}{\partial f_{i,k}^n}/(w_k|\Omega_i|)
    +\mathbf{\Psi}_k^T\dfrac{\partial J}{\partial \mathbf{W}_i^n}/|\Omega_i|=0.
  \end{cases}
\end{equation}
The macroscopic source receives the moment projection of the microscopic objective derivative, while the microscopic source receives the lifting of the macroscopic objective derivative by $\mathbf{\Psi}_k^T$. These two transformed sources are consistent: multiplying the microscopic equation in~\eqref{eq:adjoint_system_AP} by $w_k\partial f_{i,k}^+/\partial\mathbf{W}_i$ and summing over $k$ recovers the macroscopic source.

\subsection{Asymptotic limits from Eqs.~\eqref{asymptotic-preservingadjointsystem2} and~\eqref{eq:micro_coupled_W}}

The asymptotic-preserving property is verified by taking the hydrodynamic limits directly in the AP reconstructed equations Eqs.~\eqref{asymptotic-preservingadjointsystem2} and~\eqref{eq:micro_coupled_W}.

\subsubsection{Fluid limit}
In the fluid limit, $f_{i,k}^+ \to f_{i,k}$ and $\partial \mathcal{R}_{j,k}/\partial f_{i,k}\to 0$ due to the coefficient of $f_0$ term in integral solution (Eq. (~\ref{integralSolution})) tends to zero.

\textbf{Step 1: multiply Eq.~\eqref{eq:micro_coupled_W} by $\tau_i$.}
\begin{equation}
  \label{eq:micro_times_tau}
  \begin{aligned}
  &\tau_i\frac{\lambda_{i,k}^{*,n-1}-\lambda_{i,k}^{*,n}}{\Delta t}
  +\frac{\tau_i}{|\Omega_i|}\sum_j\frac{\partial \mathcal{R}_{j,k}}{\partial f_{i,k}}\lambda_{j,k}^{*,n}
  +\tau_i\mathbf{\Psi}_k^T\frac{1}{|\Omega_i|}\sum_j\sum_{k^\prime} w_{k^\prime}\frac{\partial \mathcal{R}_{j,k^\prime}}{\partial \mathbf{W}_i}\lambda_{j,k^\prime}^{*,n} \\
  &+\frac{\tau_i}{|\Omega_i|}\left(\frac{1}{w_k}\frac{\partial \mathbf{q}_i}{\partial f_{i,k}}+\mathbf{\Psi}_k^T\frac{\partial \mathbf{q}_i}{\partial \mathbf{W}_i}\right)\sum_j\sum_{k^\prime}w_{k^\prime}\frac{\partial \mathcal{R}_{j,k^\prime}}{\partial \mathbf{q}_i}\lambda_{j,k^\prime}^{*,n}\\
  &-\bigl(\mathbf{\Psi}_k^T\widetilde{\mathbf{\Lambda}}_i^{*,n-1}-\lambda_{i,k}^{*,n-1}\bigr)
  -\left(\frac{1}{w_k}\frac{\partial \mathbf{q}_i}{\partial f_{i,k}}+\mathbf{\Psi}_k^T\frac{\partial \mathbf{q}_i}{\partial \mathbf{W}_i}\right)\widetilde{\mathbf{q}}_i^{*,n}\\
  &+\mathbf{\Psi}_k^T\frac{\partial \tau_i}{\partial \mathbf{W}_i}\sum_{k^\prime} w_{k^\prime}\frac{f_{i,k^\prime}^+-f_{i,k^\prime}}{\tau_i}\lambda_{i,k^\prime}^{*,n}
  =0,
  \end{aligned}
\end{equation}
where
\begin{equation}
  \label{eq:qtilde_def}
  \widetilde{\mathbf{q}}_i^{*,n}=\sum_{k} w_k\frac{\partial f_{i,k}^+}{\partial \mathbf{q}_i}\lambda_{i,k}^{*,n}.
\end{equation}

\textbf{Step 2: asymptotic closure.}
Rearranging Eq.~\eqref{eq:micro_times_tau} and identifying the time levels under the steady-state backward marching, the temporal and flux contributions (all proportional to $\tau_i$) together with the vanishing $(f_{i,k}^+-f_{i,k})$ term may be collected into a remainder of order $o(\tau_i)$. One therefore obtains
\begin{equation}
  \label{eq:shakhov_closure_n}
  \lambda_{i,k}^{*,n}
  =
  \mathbf{\Psi}_k^T\widetilde{\mathbf{\Lambda}}_i^{*,n}
  +
  \mathbf{\mathcal{G}}_{i,k}^{*,n}
  +
  o(\tau_i),
\end{equation}
where
\begin{equation*}
  \mathbf{\mathcal{G}}_{i,k}^{*,n}
  =
  \left(\frac{1}{w_k}\frac{\partial \mathbf{q}_i}{\partial f_{i,k}}+\mathbf{\Psi}_k^T\frac{\partial \mathbf{q}_i}{\partial \mathbf{W}_i}\right)\widetilde{\mathbf{q}}_i^{*,n}.
\end{equation*}

\textbf{Step 3: continuum flux and macroscopic limiting adjoint.}
In the fluid limit, the asymptotic closure~\eqref{eq:shakhov_closure_n} is substituted into the collision-frequency block of Eq.~\eqref{asymptotic-preservingadjointsystem2}:
\begin{equation*}
  \sum_{k} w_k\frac{f_{i,k}^+-f_{i,k}}{\tau_i}\lambda_{i,k}^{*,n}
  =
  \frac{1}{\tau_i}\bigl(\widetilde{\mathbf{\Lambda}}_i^{*,n}\bigr)^T\sum_{k} w_k\bigl(f_{i,k}^+-f_{i,k}\bigr)\mathbf{\Psi}_k
  +
  \frac{1}{\tau_i}\sum_{k} w_k\bigl(f_{i,k}^+-f_{i,k}\bigr)\mathbf{\mathcal{G}}_{i,k}^{*,n}
  +o(1).
\end{equation*}
The first contribution vanishes by the compatibility condition~\eqref{eq:BGK_compat}.
For the heat-flux correction, the definition of \(\mathbf{\mathcal{G}}_{i,k}^{*,n}\) yields
\begin{equation*}
  \begin{aligned}
  \sum_{k} w_k\bigl(f_{i,k}^+-f_{i,k}\bigr)\mathbf{\mathcal{G}}_{i,k}^{*,n}
  &=
  \Biggl[
  \sum_{k}\bigl(f_{i,k}^+-f_{i,k}\bigr)\frac{\partial \mathbf{q}_i}{\partial f_{i,k}}
  +
  \Bigl(\frac{\partial \mathbf{q}_i}{\partial \mathbf{W}_i}\Bigr)^{T}
  \sum_{k} w_k\bigl(f_{i,k}^+-f_{i,k}\bigr)\mathbf{\Psi}_k
  \Biggr]
  \cdot\widetilde{\mathbf{q}}_i^{*,n}.
  \end{aligned}
\end{equation*}
The second bracket vanishes again by~\eqref{eq:BGK_compat}.
The first bracket likewise vanishes: at fixed \(\mathbf{W}_i\) the heat-flux map is linear in \(f\), and the Shakhov construction enforces \(\mathbf{q}_i[f_i^+]=\mathbf{q}_i[f_i]\), so
\begin{equation*}
  \sum_{k}\bigl(f_{i,k}^+-f_{i,k}\bigr)\frac{\partial \mathbf{q}_i}{\partial f_{i,k}}
  =
  \mathbf{q}_i[f_i^+]-\mathbf{q}_i[f_i]
  =0.
\end{equation*}
Consequently,
\begin{equation*}
  \sum_{k} w_k\frac{f_{i,k}^+-f_{i,k}}{\tau_i}\lambda_{i,k}^{*,n}
  =0.
\end{equation*}
Moreover, $\partial \mathcal{R}_{j,k}/\partial f_{i,k}\to 0$, so the corresponding blocks in Eq.~\eqref{asymptotic-preservingadjointsystem2} vanish and one is left with
\begin{equation*}
  \frac{\widetilde{\mathbf{\Lambda}}_i^{*,n-1}-\widetilde{\mathbf{\Lambda}}_i^{*,n}}{\Delta t}
  + \frac{1}{|\Omega_i|}\sum_{j}\sum_{k} w_k\frac{\partial \mathcal{R}_{j,k}}{\partial \mathbf{W}_i}\lambda_{j,k}^{*,n}=0.
\end{equation*}
To identify the macroscopic limit, the continuum expansion of the UGKS flux is required.
Starting from the integral solution and the discrete equilibrium reconstruction recalled above, in the continuum limit $\tau/\Delta t\to 0$ the free-transport contribution of $f_0$ is exponentially damped, and the leading-order time-averaged interface distribution reduces to the local Maxwellian contribution
\begin{equation*}
  f_{ij,k}
  =
  g_{0,k}
  +
  \Delta t\, g_{0,k,t}.
\end{equation*}
The corresponding flux residual then admits the asymptotic form
\begin{equation*}
  \mathcal{R}_{i,k}
  =
  \sum_{j\in N(i)}
  S_{ij}\,
  \mathbf{v}_k\cdot\mathbf{n}_{ij}\,
  f_{ij,k},
\end{equation*}
so that the moment flux reduces to the Euler flux
\begin{equation*}
  \sum_{k} w_k \mathcal{R}_{j,k}\mathbf{\Psi}_k^T
  =
  \mathbf{F}_{j}^{\mathrm{E}},
\end{equation*}
where $\mathbf{F}_{j}^{\mathrm{E}}$ is generated by the Maxwellian $g_0$ together with the $\Delta t\, g_{0,t}$ term.
Inserting the asymptotic closure~\eqref{eq:shakhov_closure_n} into the projected macroscopic equation~\eqref{asymptotic-preservingadjointsystem2}, the heat-flux correction $\mathbf{\mathcal{G}}^*$ likewise drops out of the flux projection: it is weighted by $\partial f^+/\partial \mathbf{W}$ and cancels by
\begin{equation*}
  \sum_{k} w_k\frac{\partial f_{j,k}^+}{\partial \mathbf{W}_j}\mathbf{\mathcal{G}}_{j,k}^{*,n}
  =\frac{\mathrm{d}\mathbf{q}_j}{\mathrm{d}\mathbf{W}_j}\bigg|_{f=f^+}\widetilde{\mathbf{q}}_j^{*,n}=\mathbf{0}.
\end{equation*}
Retaining the leading-order moment $\mathbf{\Psi}_k^T\widetilde{\mathbf{\Lambda}}_j^{*,n}$ and differentiating the Euler flux then yields
\begin{equation}
  \label{continuousadjointsystem}
  \frac{\widetilde{\mathbf{\Lambda}}_i^{*,n-1}-\widetilde{\mathbf{\Lambda}}_i^{*,n}}{\Delta t}
  + \frac{1}{|\Omega_i|}\sum_{j}
  \frac{\partial \mathbf{F}_{j}^{\mathrm{E}}}{\partial \mathbf{W}_i}
  \widetilde{\mathbf{W}}_{j}^{*,n}=0,
\end{equation}
which is precisely the discrete adjoint of the macroscopic Euler solver recovered by UGKS in the continuum limit.

\subsubsection{Kinetic limit}
In the kinetic limit, the UGKS numerical flux reduces to the standard upwind discrete-velocity flux: the interface state is reconstructed from the microscopic distribution alone, so that
\begin{equation*}
  \frac{\partial \mathcal{R}_{j,k}}{\partial \mathbf{W}_i}\to 0,
  \qquad
  \frac{\partial \mathcal{R}_{j,k}}{\partial \mathbf{q}_i}\to 0,
\end{equation*}
while the BGK/Shakhov collision source $(f^+-f)/\tau$ is retained exactly as in a DVM update.

\textbf{Step 1: reduce Eq.~\eqref{asymptotic-preservingadjointsystem2}.}
With $\partial \mathcal{R}_{j,k}/\partial \mathbf{W}_i=0$, Eq.~\eqref{asymptotic-preservingadjointsystem2} becomes
\begin{equation}
  \label{eq:kinetic_macro_adjoint}
  \begin{aligned}
  \frac{\widetilde{\mathbf{\Lambda}}_i^{*,n-1}-\widetilde{\mathbf{\Lambda}}_i^{*,n}}{\Delta t}
  &+\frac{1}{|\Omega_i|}\sum_{k} w_k\frac{\partial f_{i,k}^+}{\partial \mathbf{W}_i}\sum_j \frac{\partial \mathcal{R}_{j,k}}{\partial f_{i,k}}\lambda_{j,k}^{*,n}\\
  &+\frac{1}{\tau_i}\frac{\partial \tau_i}{\partial \mathbf{W}_i}\sum_{k} w_k\frac{f_{i,k}^+-f_{i,k}}{\tau_i}\lambda_{i,k}^{*,n}
  =0.
  \end{aligned}
\end{equation}

\textbf{Step 2: reduce Eq.~\eqref{eq:micro_coupled_W}.}
Under the same flux limit, Eq.~\eqref{eq:micro_coupled_W} reduces to
\begin{equation}
  \label{eq:kinetic_micro_adjoint}
  \begin{aligned}
  \frac{\lambda_{i,k}^{*,n-1}-\lambda_{i,k}^{*,n}}{\Delta t}
  &+\frac{1}{|\Omega_i|}\sum_j\frac{\partial \mathcal{R}_{j,k}}{\partial f_{i,k}}\lambda_{j,k}^{*,n}
  -\frac{\mathbf{\Psi}_k^T\widetilde{\mathbf{\Lambda}}_i^{*,n-1}-\lambda_{i,k}^{*,n-1}}{\tau_i}\\
  &-\frac{1}{\tau_i}\left(\frac{1}{w_k}\frac{\partial \mathbf{q}_i}{\partial f_{i,k}}+\mathbf{\Psi}_k^T\frac{\partial \mathbf{q}_i}{\partial \mathbf{W}_i}\right)\sum_{k^\prime} w_{k^\prime}  \frac{\partial f_{i,k^\prime}^+}{\partial \mathbf{q}_i} \lambda_{i,k^\prime}^{*,n}\\
  &+\mathbf{\Psi}_k^T\frac{1}{\tau_i}\frac{\partial \tau_i}{\partial \mathbf{W}_i}\sum_{k^\prime} w_{k^\prime}\frac{f_{i,k^\prime}^+-f_{i,k^\prime}}{\tau_i}\lambda_{i,k^\prime}^{*,n}
  =0.
  \end{aligned}
\end{equation}
The first two terms are the adjoint of the upwind DVM transport residual, while the remaining $O(\tau_i^{-1})$ blocks are the adjoint of the DVM collision operator (including the Shakhov heat-flux coupling and the $\tau$-derivative contribution).
Thus Eq.~\eqref{eq:kinetic_micro_adjoint} is precisely the discrete adjoint of the DVM scheme with BGK/Shakhov collisions.

\textbf{Step 3: consistency with the macroscopic projection.}
Multiplying Eq.~\eqref{eq:kinetic_micro_adjoint} by $w_k\partial f_{i,k}^+/\partial \mathbf{W}_i$ and summing over $k$ recovers Eq.~\eqref{eq:kinetic_macro_adjoint}.
Hence the macroscopic adjoint~\eqref{eq:kinetic_macro_adjoint} is the moment projection of the microscopic DVM adjoint, and the AP system collapses to the classical DVM adjoint equations in the rarefied regime.
This confirms the asymptotic-preserving property as $Kn\to\infty$.

\subsection{Adjoint Boundary Conditions and Sensitivity Analysis for Solid Walls}

At the solid wall boundary $\partial\Omega_i$ with area $S_w$, the Maxwell reflection model with accommodation coefficient $\alpha\in[0,1]$ is employed.
The wall distribution $\tilde{f}_{k,w}$ combines diffuse and specular reflections:
\begin{equation*}
\tilde{f}_{k,w}=\begin{cases}
\tilde{f}_{k,\mathrm{in}}(\mathbf{W}_{\mathrm{in}},f_{\mathrm{in},k}), & u_k<0, \\
\alpha\,\rho_w g_k(1, U_w, \lambda_w)+(1-\alpha)\,\tilde{f}_{k^s,\mathrm{in}}, & u_k>0,
\end{cases}
\qquad
\rho_w = -\dfrac{\displaystyle\sum_{u_k<0}w_ku_k\tilde{f}_{k,\mathrm{in}}}{\displaystyle\sum_{u_k>0}w_ku_k g_k(1)},
\end{equation*}
where $k^s$ denotes the specularly reflected discrete-velocity index of $k$, $g_k(1)$ is the normalized wall Maxwellian with unit density, and $\lambda_w$ is the wall parameter (e.g., the inverse of the wall temperature).
The limits $\alpha=1$ and $\alpha=0$ recover pure diffuse and pure specular reflection, respectively.
Because the specular contribution conserves the wall-normal mass flux, the wall density $\rho_w$ retains the same impermeability formula as in the diffuse case.
The local boundary functional incorporating the macro-micro coupled constraints yields the inner product for the $i$-th cell:
\begin{equation*}
  \begin{aligned}
  \mathcal{I}_{\mathrm{wall}, i} = & \sum_{n=1}^{N} |\Omega_i|\left( \mathbf{W}_{\mathrm{in}, i}^n - \sum_{k}w_k f_{\mathrm{in},k,i}^n \mathbf{\Psi}_k \right)^T \mathbf{\Lambda}_{\mathrm{in}, i}^{*,n}  \\
  &+ \sum_{n=1}^{N-1} \sum_{k} \left( \frac{f_{\mathrm{in},k,i}^{n+1} - f_{\mathrm{in},k,i}^n}{\Delta t} w_k |\Omega_i| + S_w w_k u_k \tilde{f}_{k,w,i}^n \right) \lambda_{k,\mathrm{in},i}^{*,n}.
\end{aligned}
\end{equation*}

By invoking the variational principle $\delta \mathcal{I}_{\mathrm{wall}, i} = 0$ at time step $n$ ($n=2, \dots, N-1$) and utilizing the integration-by-parts shift for the temporal derivative, the coupling terms associated with the incident flux $\tilde{f}_{k,\mathrm{in}}^n$ ($u_k<0$) are grouped.
Setting the coefficients of $\delta \mathbf{W}_{\mathrm{in},i}^n$ and $\delta f_{\mathrm{in},k,i}^n$ to zero yields the macroscopic wall adjoint condition and the microscopic adjoint wall-flux contribution.
For convenience, the diffuse adjoint weight is denoted by
\begin{equation*}
\mathcal{A}_{k}^{*,n}
=
\lambda_{k,\mathrm{in},i}^{*,n}
-
\frac{\displaystyle\sum_{u_{k^\prime}>0} w_{k^\prime} u_{k^\prime} g_{k^\prime}(1)\, \lambda_{k^\prime,\mathrm{in},i}^{*,n}}{\displaystyle\sum_{u_{k^\prime}>0} w_{k^\prime} u_{k^\prime} g_{k^\prime}(1)}.
\end{equation*}

\paragraph{Macroscopic Adjoint Equation:}
\begin{equation}
\mathbf{\Lambda}_{\mathrm{in}, i}^{*,n} 
=
-\frac{S_w}{\left|\Omega_i\right|}
\sum_{u_k<0} w_k u_k
\left( \frac{\partial \tilde{f}_{k,\mathrm{in}}^n}{\partial \mathbf{W}_{\mathrm{in}}^n} \right)^T
\left[
\alpha\,\mathcal{A}_{k}^{*,n}
+(1-\alpha)\bigl(\lambda_{k,\mathrm{in},i}^{*,n}-f_{k^s,\mathrm{in},i}^{*,n}\bigr)
\right].
\end{equation}

\paragraph{Microscopic Adjoint Wall Flux:}
The microscopic adjoint contribution of the Maxwell wall appears only through the boundary flux and is the sum of the macroscopic-$\mathbf{W}$ part and the microscopic-$f$ part,
\begin{equation}
\mathcal{F}_{\mathrm{bound}}^n
=
\mathcal{F}_{\mathbf{W}}^n
+
\mathcal{F}_{f}^n,
\end{equation}
where
\begin{equation*}
\mathcal{F}_{\mathbf{W}}^n
=
-(\mathbf{\Psi}_k)^T \mathbf{\Lambda}_{\mathrm{in}, i}^{*,n} ,
\end{equation*}
and, for $u_k<0$,
\begin{equation*}
\mathcal{F}_{f}^n
=
\frac{S_w u_k}{\left|\Omega_i\right|}
\frac{\partial \tilde{f}_{k,\mathrm{in}}^n}{\partial f_{\mathrm{in},k}^n}
\left[
\alpha\,\mathcal{A}_{k}^{*,n}
+(1-\alpha)\bigl(\lambda_{k,\mathrm{in},i}^{*,n}-f_{k^s,\mathrm{in},i}^{*,n}\bigr)
\right],
\end{equation*}
while $\mathcal{F}_{f}^n=0$ for $u_k>0$.
When $\alpha=1$, the specular coupling in $\mathcal{F}_{f}^n$ vanishes and the pure-diffuse adjoint wall flux is recovered.

\section{Implementation of the Adjoint Method}

This section discusses the numerical implementation of the adjoint system. For discrete-velocity methods such as the UGKS, memory consumption is a major bottleneck; therefore, the efficiency- and memory-saving UGKS programming paradigm~\cite{zhang2025efficiency} is adopted first for adjoint residual evaluation. To accelerate the adjoint solve, the implicit UGKS marching strategy~\cite{zhu2016implicit} is then employed.

\subsection{Calculation of adjoint-system residuals}
Consider the microscopic adjoint equation in the form of Eq.~\eqref{eq:micro_coupled_W}. The contributions associated with $\mathbf{W}_i$ and $\mathbf{q}_i$ involve summations over the full velocity space, whereas the contribution associated with $f_{i,k}$ depends only on the local velocity point $k$. Moreover, the $\mathbf{W}$- and $\mathbf{q}$-related terms can be decomposed into a velocity-space sum that is common to all points and coefficients that vary with $k$. To avoid redundant summations over the velocity space, the forward calculation strategy in~\cite{zhang2025efficiency} is followed: the common parts of the $\mathbf{W}$ and $\mathbf{q}$ contributions are precomputed, and the velocity-space loop is then used to assemble the remaining pointwise terms. In the adjoint residual evaluation, the same idea is applied in reverse order within the velocity loop: the common part is evaluated first, and the residual at each velocity point is subsequently completed by adding its local contribution.

\begin{figure}[ht]
  \centering
  \includegraphics[width=0.5\textwidth]{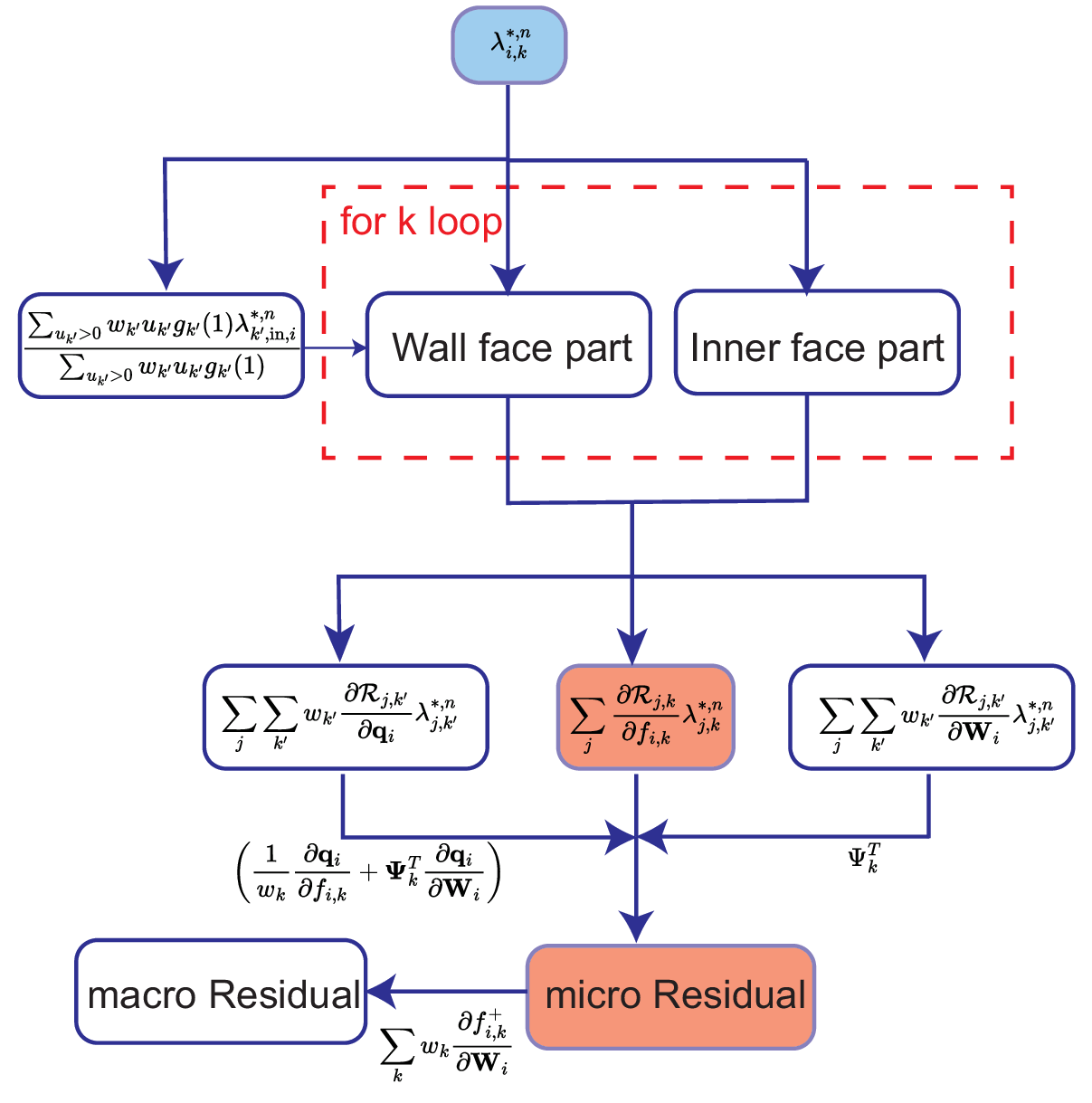}
  \caption{Organization of adjoint residual evaluation following the UGKS programming paradigm~\cite{zhang2025efficiency}.}
  \label{fig:resAlg}
\end{figure}

A similar treatment is required for the adjoint boundary conditions, where the Maxwell wall condition introduces coupled incident, specular and diffuse reflected flux contributions. In the forward UGKS, the incident flux is computed first, and the reflected flux is added afterward. In the adjoint implementation, the reflected contribution must be evaluated before the incident contribution. The overall residual-evaluation procedure is summarized in Fig.~\ref{fig:resAlg}. With the above strategy, only the adjoint variables and their residuals need to be stored throughout the microscopic adjoint solve. In addition, the microscopic distribution functions must be retained to evaluate the required partial derivatives. Although this introduces some redundant computation, the overall memory consumption remains under control. All remaining variables are independent of the velocity-space size, so the memory requirement of the adjoint implementation stays manageable.

\subsection{Implicit marching}

As shown in the asymptotic analysis above, the macroscopic adjoint equation recovers the adjoint system of the macroscopic Euler equations in the continuum limit, whereas the microscopic adjoint equation recovers the adjoint system of the discrete velocity method in the free-molecular limit. To solve the coupled adjoint system implicitly, the predictor--corrector marching algorithm in~\cite{zhu2016implicit} is adopted. At each time step, the macroscopic adjoint equation is solved first to obtain the predicted macroscopic adjoint variable $\widetilde{\mathbf{\Lambda}}_i^{*,n+1/2}$; the microscopic adjoint equation is then solved to obtain the predicted microscopic adjoint variable $f_{i,k}^{*,n+1}$. A final correction step yields the updated adjoint variables $\widetilde{\mathbf{\Lambda}}_i^{*,n+1}$. This procedure avoids a global matrix inversion over the full velocity space and therefore substantially reduces the computational cost.

In the predictor step, the macroscopic adjoint equation is updated by the lower--upper symmetric Gauss--Seidel (LU-SGS) scheme~\cite{yoon1988lower}. Linearizing the macroscopic adjoint residual at time level $n$ gives
\begin{equation*}
\frac{1}{\Delta t}\Delta \widetilde{\mathbf{\Lambda}}_i^{*,n}+\frac{1}{|\Omega_i|}\sum_{j\in\mathcal{N}(i)}S_{ij}\Delta \hat{\mathbf{F}}_{ij}^{*,n}=\widetilde{\mathbf{R}}_{\text{macro},i}^{*,n},
\end{equation*}
where $\widetilde{\mathbf{R}}_{\text{macro},i}^{*,n}$ is the macroscopic adjoint residual. The variation of the numerical flux is approximated by the Euler-equation-based flux variation. At the interface between cells $i$ and $j$, $\Delta \hat{\mathbf{F}}_{ij}^{*,n}$ is approximated as
\begin{equation*}
\Delta \hat{\mathbf{F}}_{ij}^{*,n} \approx \frac{1}{2}\left[\left(\frac{\partial \mathbf{F}_{\text{Euler},i}}{\partial \mathbf{W}_i}\right)^{\!T}\Delta \widetilde{\mathbf{\Lambda}}_i^{*,n}+\left(\frac{\partial \mathbf{F}_{\text{Euler},j}}{\partial \mathbf{W}_j}\right)^{\!T}\Delta \widetilde{\mathbf{\Lambda}}_j^{*,n}-\Gamma_{ij}\left(\Delta \widetilde{\mathbf{\Lambda}}_i^{*,n}-\Delta \widetilde{\mathbf{\Lambda}}_j^{*,n}\right)\right],
\end{equation*}
where $\Gamma_{ij}$ is the spectral radius of the Euler flux Jacobian at the interface,
\begin{equation*}
\Gamma_{ij}=\left|\mathbf{V}_{ij}\cdot\mathbf{n}_{ij}\right|+a_{ij}.
\end{equation*}
The resulting linearized system can be solved by the standard LU-SGS method. For the adjoint formulation, the neighbor coupling appears through the transpose of the primal flux Jacobian, the cell sweep order is reversed, and the sign of the $\Gamma_{ij}$ term is opposite to that in the forward solver. After the LU-SGS solve, $\widetilde{\mathbf{\Lambda}}_i^{*,n-1/2}$ is obtained and passed to the microscopic adjoint update.

In the predictor step, the microscopic adjoint equation, Eq.~\eqref{eq:micro_coupled_W}, is rewritten in migration--relaxation form as
\begin{equation*}
\frac{\lambda_{i,k}^{*,n-1}-\lambda_{i,k}^{*,n}}{\Delta t}+\hat{r}_{i,k,\text{mig}}^{*,n}+s_{i,k}^{*,n}=\frac{\mathbf{\Psi}_k^T\widetilde{\mathbf{\Lambda}}_i^{*,n-1/2}-\lambda_{i,k}^{*,n-1}}{\tau_i}.
\end{equation*}
Adding $(\lambda_{i,k}^{*,n-1}-\lambda_{i,k}^{*,n})/\tau_i$ to both sides gives
\begin{equation*}
\frac{\lambda_{i,k}^{*,n-1}-\lambda_{i,k}^{*,n}}{\Delta t}+\frac{\lambda_{i,k}^{*,n-1}-\lambda_{i,k}^{*,n}}{\tau_i}+\hat{r}_{i,k,\text{mig}}^{*,n}+s_{i,k}^{*,n}=\frac{\mathbf{\Psi}_k^T\widetilde{\mathbf{\Lambda}}_i^{*,n-1/2}-\lambda_{i,k}^{*,n}}{\tau_i},
\end{equation*}
so that the unknown $\lambda_{i,k}^{*,n}$ appears on the right-hand side and can be treated implicitly together with the time derivative. Here the migration part collects the spatial flux contributions,
\begin{equation*}
  \begin{aligned}    
    \hat{r}_{i,k,\text{mig}}^{*,n}={}& \frac{1}{|\Omega_i|}\sum_j\frac{\partial \mathcal{R}_{j,k}}{\partial f_{i,k}}\lambda_{j,k}^{*,n} + \mathbf{\Psi}_k^T\frac{1}{|\Omega_i|}\sum_j\sum_{k^\prime} w_{k^\prime}\frac{\partial \mathcal{R}_{j,k^\prime}}{\partial \mathbf{W}_i}\lambda_{j,k^\prime}^{*,n} \\
    &+ \frac{1}{|\Omega_i|}\left(\frac{1}{w_k}\frac{\partial \mathbf{q}_i}{\partial f_{i,k}}+\mathbf{\Psi}_k^T\frac{\partial \mathbf{q}_i}{\partial \mathbf{W}_i}\right)\sum_j\sum_{k^\prime}w_{k^\prime}\frac{\partial \mathcal{R}_{j,k^\prime}}{\partial \mathbf{q}_i}\lambda_{j,k^\prime}^{*,n},
  \end{aligned}
\end{equation*}
and the remaining source terms are treated explicitly,
\begin{equation*}
s_{i,k}^{*,n}=-\frac{1}{\tau_i}\left(\frac{1}{w_k}\frac{\partial \mathbf{q}_i}{\partial f_{i,k}}+\mathbf{\Psi}_k^T\frac{\partial \mathbf{q}_i}{\partial \mathbf{W}_i}\right)\sum_{k^\prime} w_{k^\prime}  \frac{\partial f_{i,k^\prime}^+}{\partial \mathbf{q}_i} \lambda_{i,k^\prime}^{*,n}
+\mathbf{\Psi}_k^T\frac{1}{\tau_i}\frac{\partial \tau_i}{\partial \mathbf{W}_i}\sum_{k^\prime} w_{k^\prime}\frac{f_{i,k^\prime}^+-f_{i,k^\prime}}{\tau_i}\lambda_{i,k^\prime}^{*,n}.
\end{equation*}
Only the time derivative, the shifted relaxation term, and the migration flux are treated implicitly. Their linearization gives
\begin{equation*}
\frac{1}{\Delta t}\Delta \lambda_{i,k}^{*,n}+\frac{1}{\tau_i}\Delta \lambda_{i,k}^{*,n}+\frac{1}{|\Omega_i|}\sum_{j\in\mathcal{N}(i)}S_{ij}\Delta \hat{\mathcal{F}}_{ij,k}^{*,n}=\frac{\mathbf{\Psi}_k^T\widetilde{\mathbf{\Lambda}}_i^{*,n-1}-\lambda_{i,k}^{*,n-1}}{\tau_i}-\hat{r}_{i,k,\text{mig}}^{*,n}-s_{i,k}^{*,n},
\end{equation*}
where the terms $(\lambda_{i,k}^{*,n-1}-\lambda_{i,k}^{*,n})/\Delta t$ and $(\lambda_{i,k}^{*,n-1}-\lambda_{i,k}^{*,n})/\tau_i$ on the left-hand side are linearized as $(1/\Delta t)\Delta \lambda_{i,k}^{*,n}$ and $(1/\tau_i)\Delta \lambda_{i,k}^{*,n}$, and the $\lambda_{i,k}^{*,n}/\tau_i$ term on the right-hand side is moved to the left-hand side during the implicit solve. The variation of the microscopic numerical flux is approximated by the upwind form. At the interface between cells $i$ and $j$,
\begin{equation*}
\Delta \hat{\mathcal{F}}_{ij,k}^{*,n} \approx \frac{1}{2}\left[(\hat{w}_{ij,k}-|\hat{w}_{ij,k}|)\Delta \lambda_{i,k}^{*,n}+(\hat{w}_{ij,k}+|\hat{w}_{ij,k}|)\Delta \lambda_{j,k}^{*,n}\right],
\end{equation*}
with
\begin{equation*}
\sigma_{ij,k}=\left|\hat{w}_{ij,k}\right|, \qquad \hat{w}_{ij,k}=\mathbf{v}_k\cdot\mathbf{n}_{ij}.
\end{equation*}
The resulting linearized system is solved by the same pointwise iteration strategy as in~\cite{zhu2016implicit}. For the adjoint formulation, the flux Jacobian appears in transposed form and the sign of the $\sigma_{ij,k}$ term is opposite to that in the forward solver. After convergence, $\lambda_{i,k}^{*,n-1}$ is obtained and used in the subsequent correction step.

In the correction step, the updated macroscopic adjoint moment is reconstructed from the microscopic adjoint variable as
\begin{equation*}
\widetilde{\mathbf{\Lambda}}_i^{*,n-1}=\widetilde{\mathbf{\Lambda}}_i^{*,n-1/2}+\sum_{k} w_k \frac{\partial f_{i,k}^+}{\partial \mathbf{W}_i}( \lambda_{i,k}^{*,n-1}-\mathbf{\Psi}_k^T \widetilde{\mathbf{\Lambda}}_i^{*,n-1/2}),
\end{equation*}
which enforces consistency between the macroscopic and microscopic adjoint variables at the new time level.

\section{Test cases}
This section verifies the correctness and effectiveness of the adjoint UGKS.
Comparisons are made against the linearized UGKS solver (L-UGKS) described in ~\ref{app:linear_UGKS}.
\subsection{Cavity flow}
The first test case is the heat conduction problem of lid-driven cavity flow. Two representative flow regimes are first considered: $Re=100$ and $Kn=0.075$; the method is then further examined at $Re=1000$ and $Kn=1$. The wall temperature is fixed at $T_w=273\text{K}$ for all cases. The velocity space is discretized with a $60\times60$ grid over the domain $(-5\sqrt{RT_w},5\sqrt{RT_w})^2$. In physical space, the $Kn=0.075$ and $Kn=1$ cases use a uniform grid of $60\times60$, while the $Re=100$ and $Re=1000$ cases employ a $70\times70$ grid with a near-wall refinement scale of $0.004$. The working gas is argon with a Prandtl number $Pr=2/3$.

The investigated problem is the sensitivity of the total heat flux on the right wall with respect to the temperature on the left wall. A linearized solver is used for comparison. For the linearized program, six perturbation points are selected at the following $y$ positions:
\begin{equation*}
x_n = \frac{1}{60}(25 -10n -\frac{1}{2} ), \quad n=0,1,2,3,4,5
\end{equation*}
This selection provides a regular distribution for analysis and facilitates extensibility. At each perturbation point, a small temperature disturbance $\delta T=1\times 10^{-3}\text{K}$ is imposed on the left wall, and numerical simulations are conducted to compute and analyze the sensitivity response of the total heat flux on the right wall using linearized approaches.

\begin{figure}[!htbp]
  \centering
  \begin{subfigure}[b]{0.35\textwidth}
    \centering
    \includegraphics[width=\linewidth]{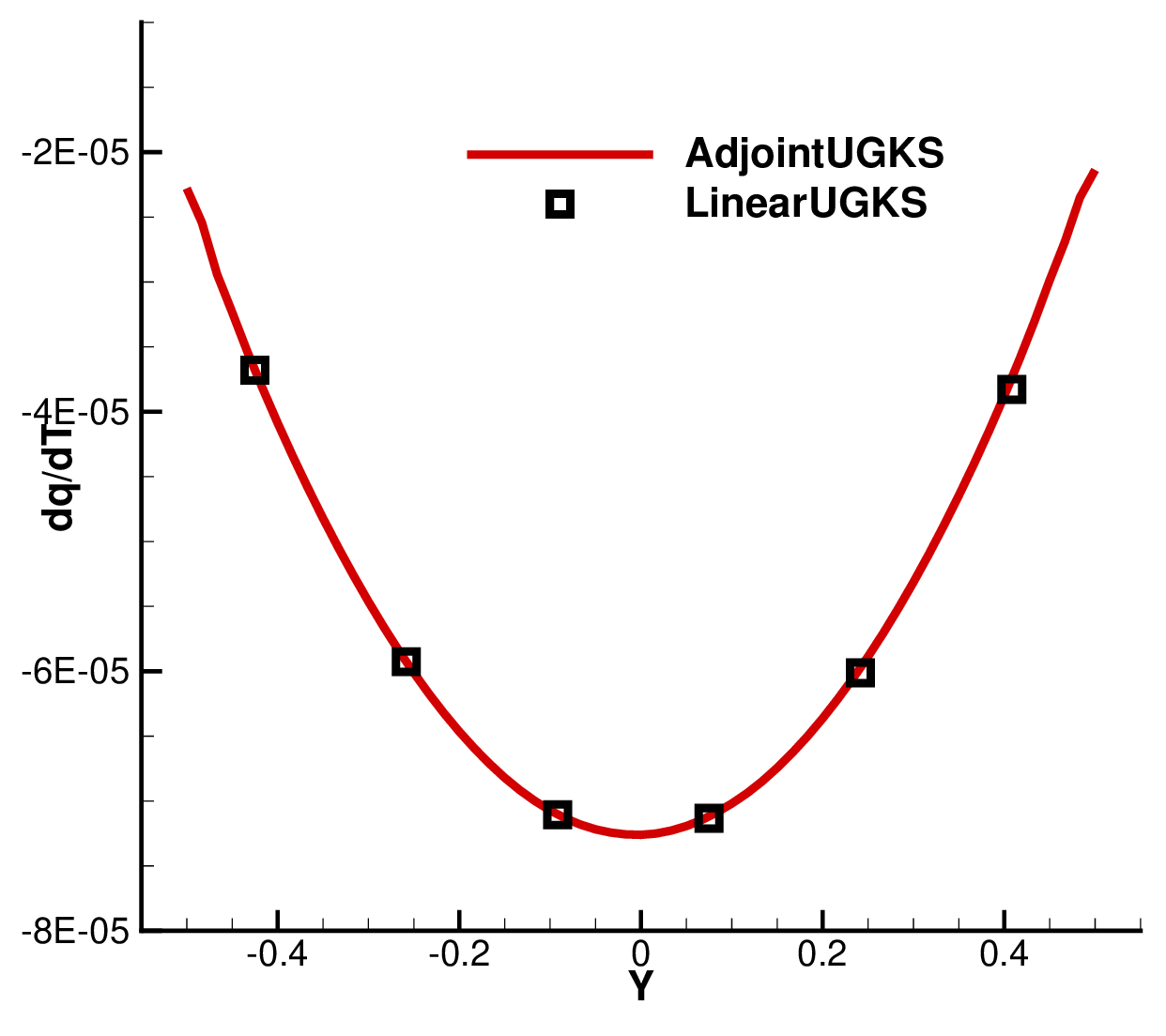}
    \caption{}
    \label{fig:cavityKn075_sensity}
  \end{subfigure}
  \quad
  \begin{subfigure}[b]{0.35\textwidth}
    \centering
    \includegraphics[width=\linewidth]{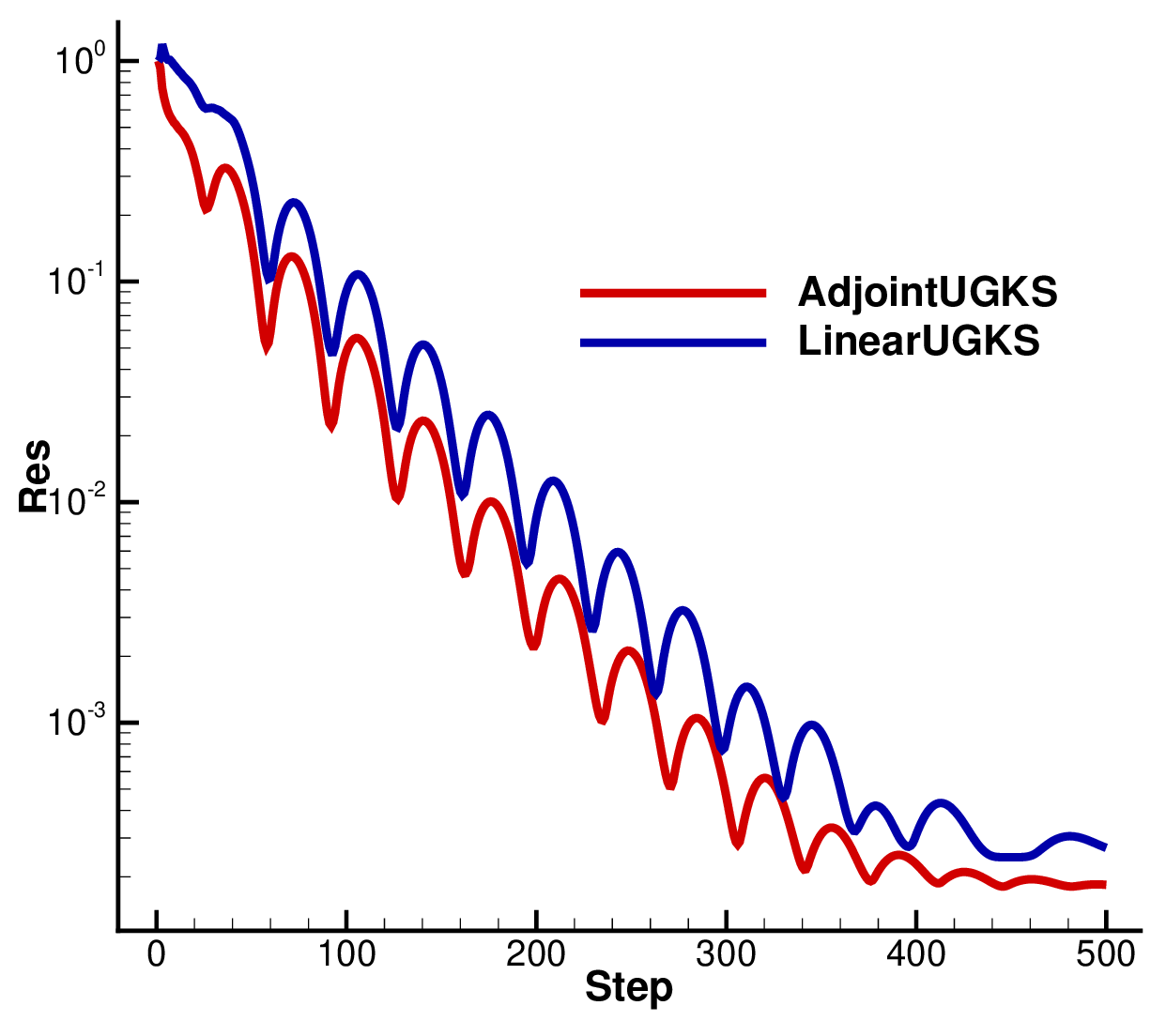}
    \caption{}
    \label{fig:cavityKn075_residual}
  \end{subfigure}
  \caption{Cavity flow at Knudsen number 0.075. (a) Sensitivity to the left wall temperature. (b) Residual of the adjoint system and the linear system.}
  \label{fig:cavityKn075}
\end{figure}

Figure~\ref{fig:cavityKn075} compares the sensitivities of the total heat flux on the right wall with respect to the left-wall temperature obtained by the linearized and adjoint methods for $Kn=0.075$. The convergence histories of the linearized and adjoint residuals at $n=0$ are also shown. The two sensitivity curves coincide almost exactly, and the residual histories are in excellent agreement. These results verify the correctness of the derived adjoint system and the consistency of its numerical implementation.

\begin{figure}[!htbp]
  \centering
  \begin{subfigure}[b]{0.32\textwidth}
    \centering
    \includegraphics[width=\linewidth]{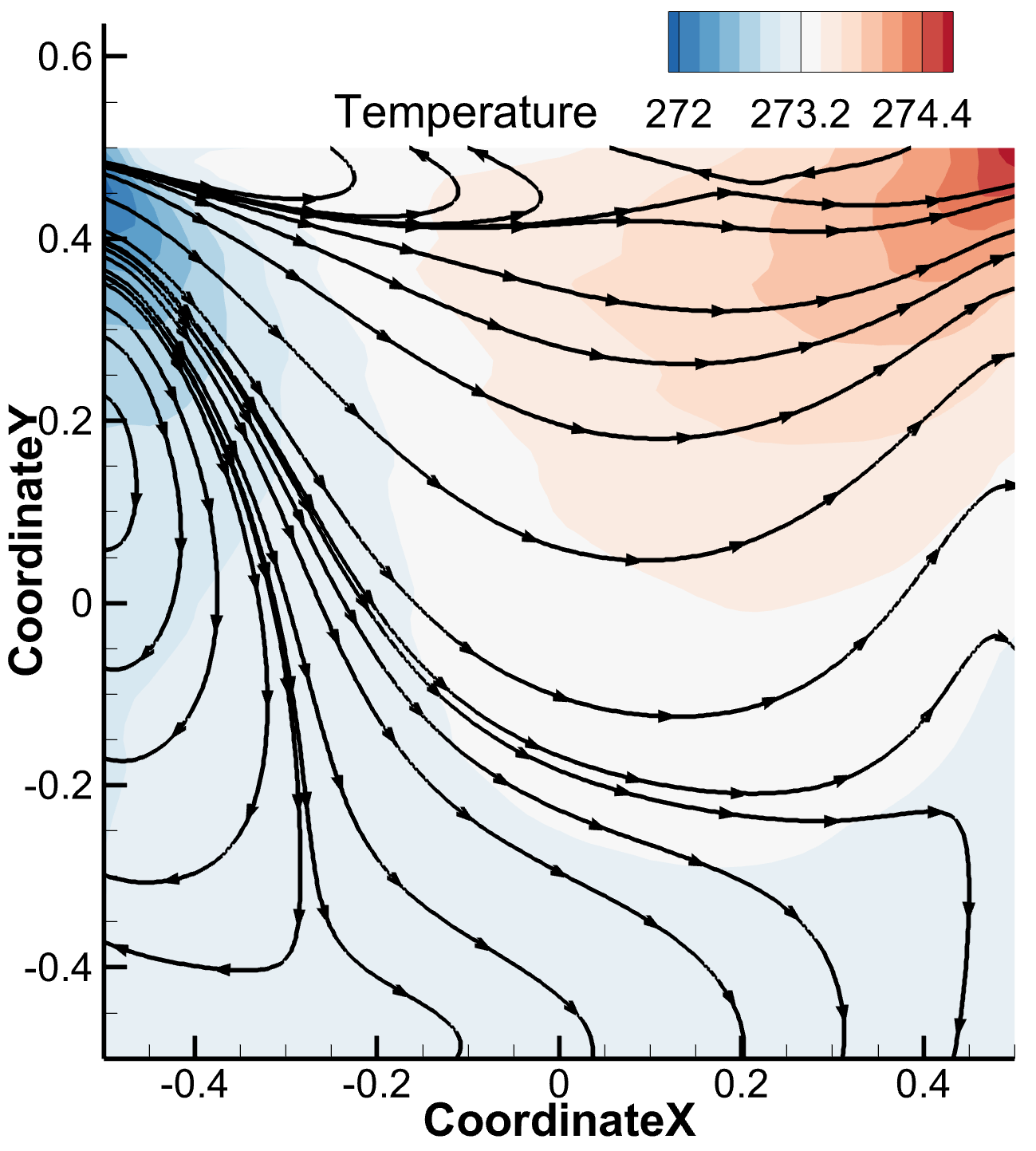}
    \caption{}
    \label{fig:cavityKn075_origin}
  \end{subfigure}
  \hfill
  \begin{subfigure}[b]{0.32\textwidth}
    \centering
    \includegraphics[width=\linewidth]{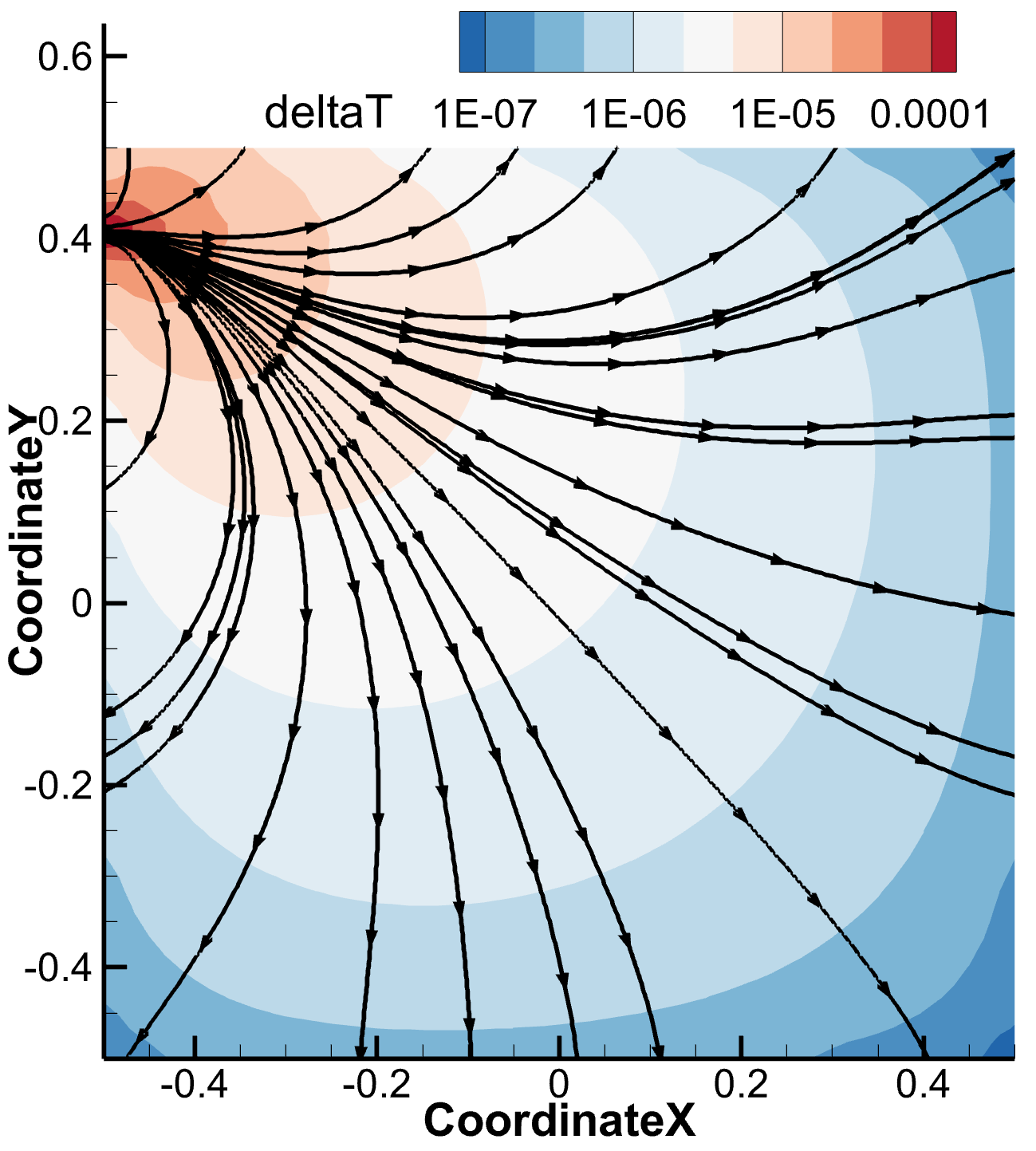}
    \caption{}
    \label{fig:cavityKn075_linear}
  \end{subfigure}
  \hfill
  \begin{subfigure}[b]{0.32\textwidth}
    \centering
    \includegraphics[width=\linewidth]{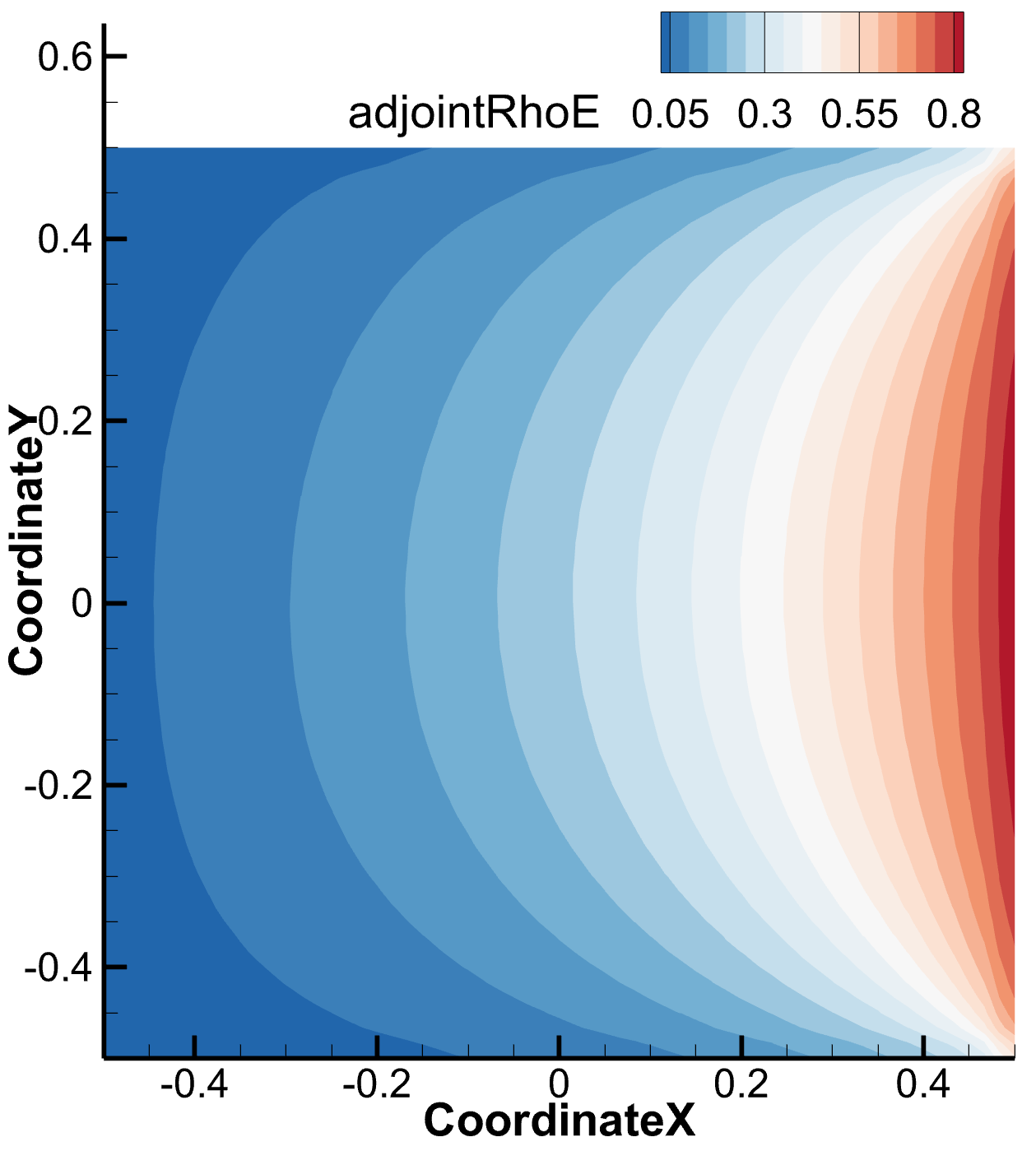}
    \caption{}
    \label{fig:cavityKn075_adjoint}
  \end{subfigure}
  \caption{Field distributions of cavity flow at Knudsen number 0.075. (a) Baseline temperature field and heat-flux streamlines. (b) Linearized field at perturbation location $y=0.40833$ (contours: perturbed temperature; streamlines: perturbed heat flux). (c) Adjoint energy field.}
  \label{fig:cavityKn075_contour}
\end{figure}

The field distributions for the $Kn=0.075$ case are shown in Fig.~\ref{fig:cavityKn075_contour}. The baseline temperature field and heat-flux streamlines in Fig.~\ref{fig:cavityKn075_origin} reveal a transitional-regime response in which the heat flux is not aligned with the temperature gradient, indicating a departure from Fourier's law. Figure~\ref{fig:cavityKn075_linear} presents the linearized field at $y=0.40833$, where the perturbed heat flux emanates directly from the disturbance location on the left wall. The corresponding adjoint energy field in Fig.~\ref{fig:cavityKn075_adjoint} exhibits the reverse propagation of the thermal response from the objective region toward the perturbed wall temperature.

\begin{figure}[!htbp]
  \centering
  \begin{subfigure}[b]{0.35\textwidth}
    \centering
    \includegraphics[width=\linewidth]{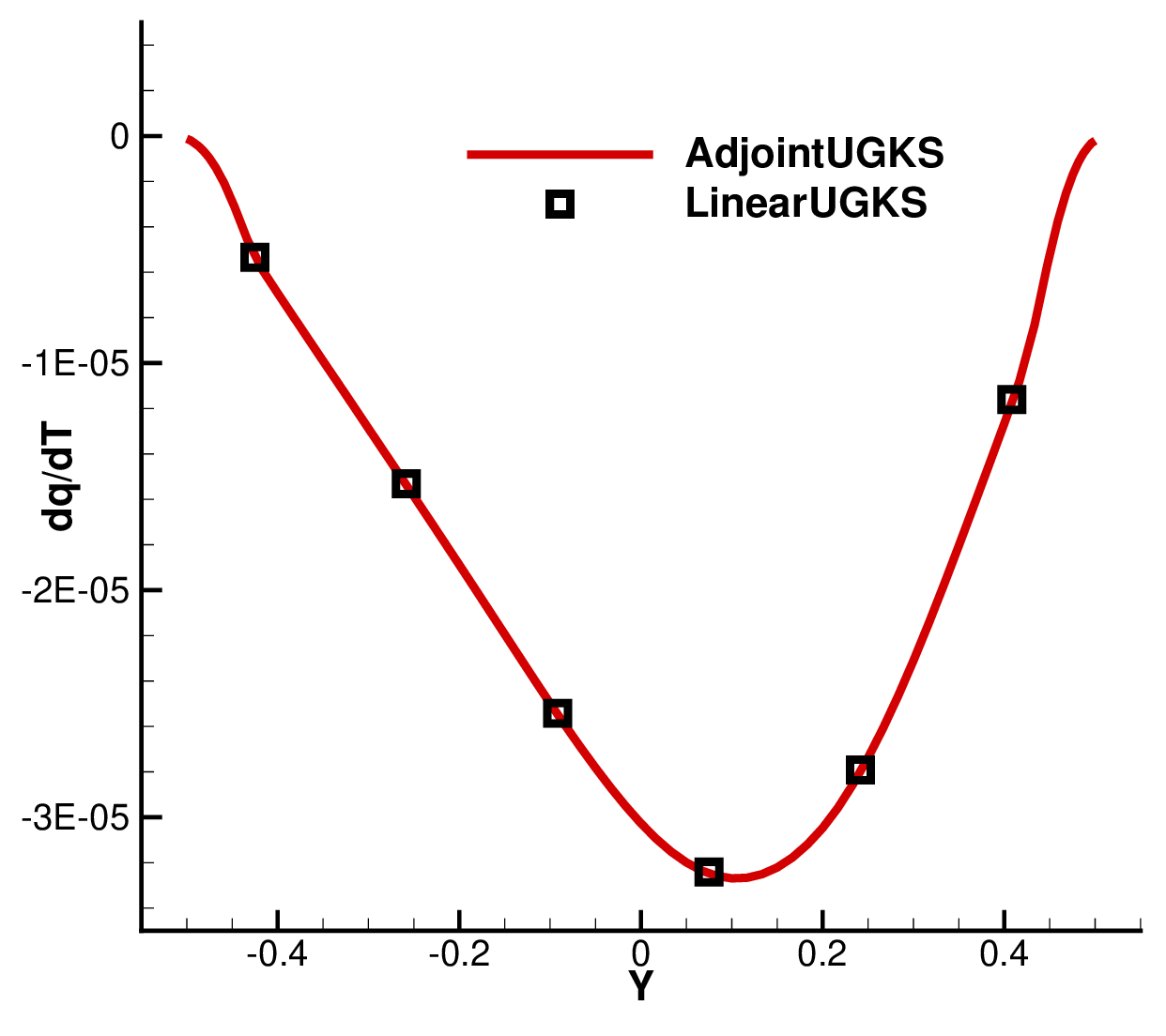}
    \caption{}
    \label{fig:cavityRe100_sensity}
  \end{subfigure}
  \quad
  \begin{subfigure}[b]{0.35\textwidth}
    \centering
    \includegraphics[width=\linewidth]{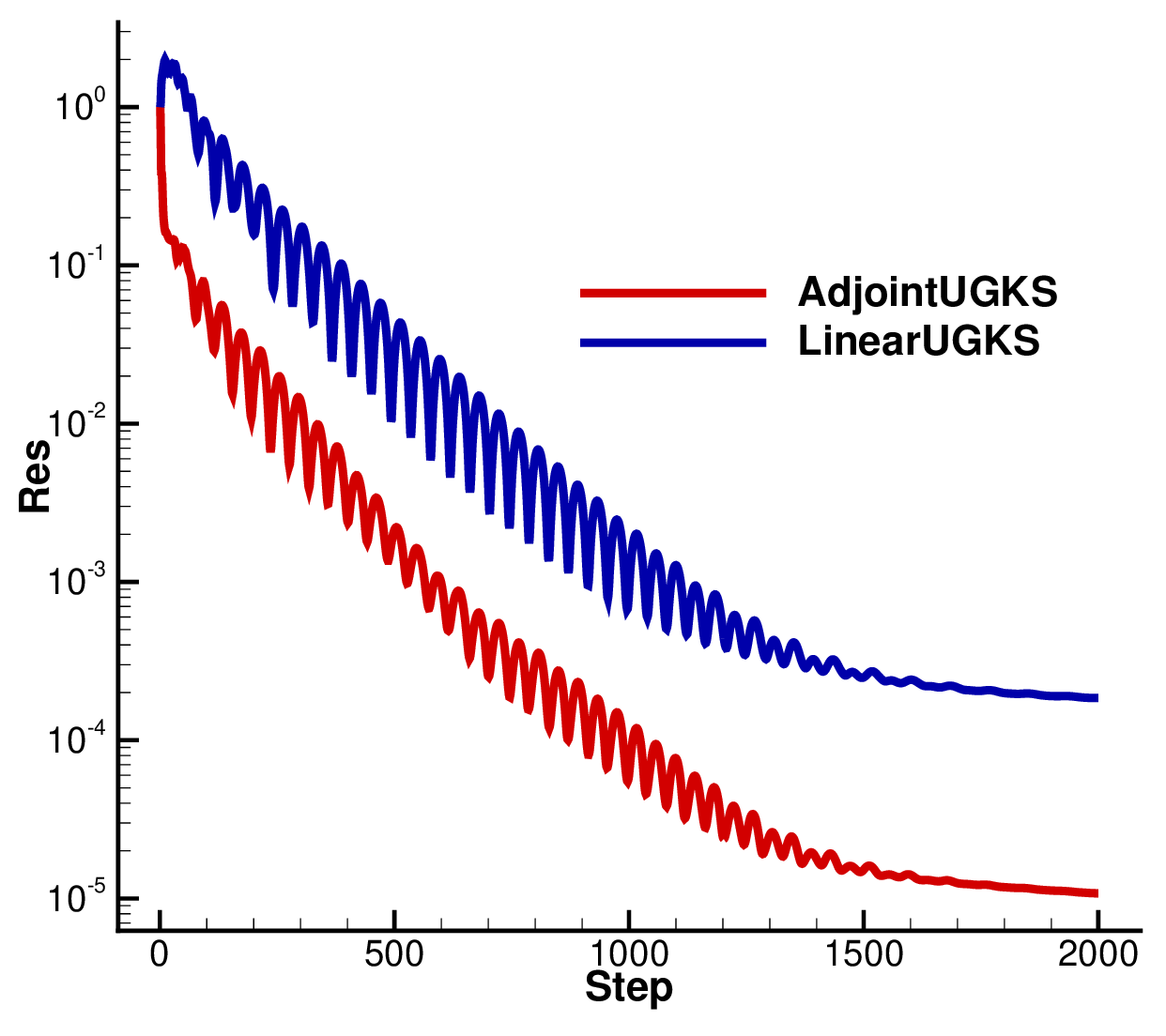}
    \caption{}
    \label{fig:cavityRe100_residual}
  \end{subfigure}
  \caption{Cavity flow at Reynolds number 100. (a) Sensitivity to the left wall temperature. (b) Residual of adjoint system and linear system.}
  \label{fig:cavityRe100}
\end{figure}

Figure~\ref{fig:cavityRe100} presents the corresponding results for the continuum case at $Re=100$. Again, the adjoint and linearized sensitivities agree closely, and the residual histories converge at essentially the same rate. This confirms that the adjoint formulation remains valid in the near-continuum regime.

\begin{figure}[!htbp]
  \centering
  \begin{subfigure}[b]{0.32\textwidth}
    \centering
    \includegraphics[width=\linewidth]{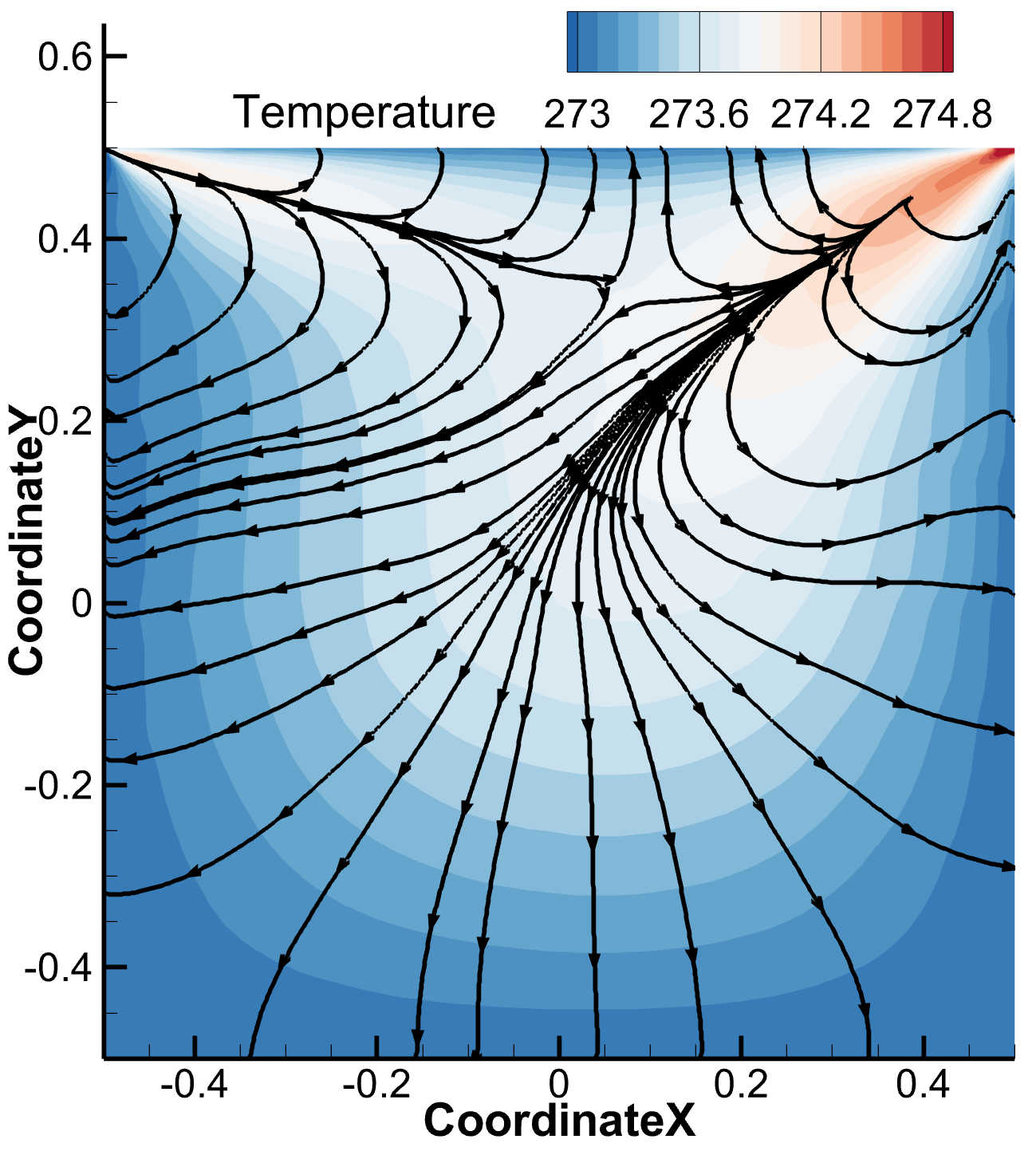}
    \caption{}
    \label{fig:cavityRe100_origin}
  \end{subfigure}
  \hfill
  \begin{subfigure}[b]{0.32\textwidth}
    \centering
    \includegraphics[width=\linewidth]{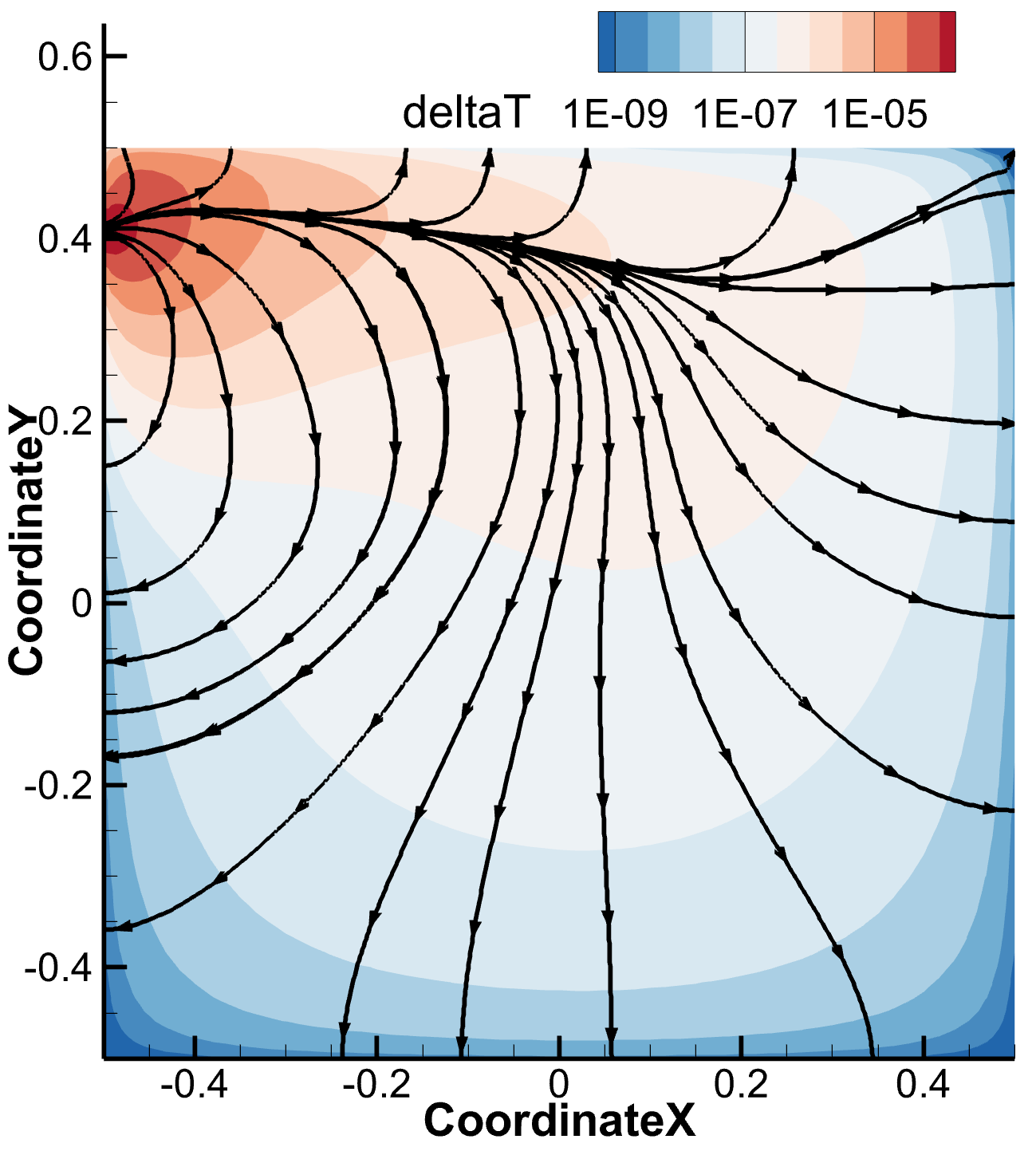}
    \caption{}
    \label{fig:cavityRe100_linear}
  \end{subfigure}
  \hfill
  \begin{subfigure}[b]{0.32\textwidth}
    \centering
    \includegraphics[width=\linewidth]{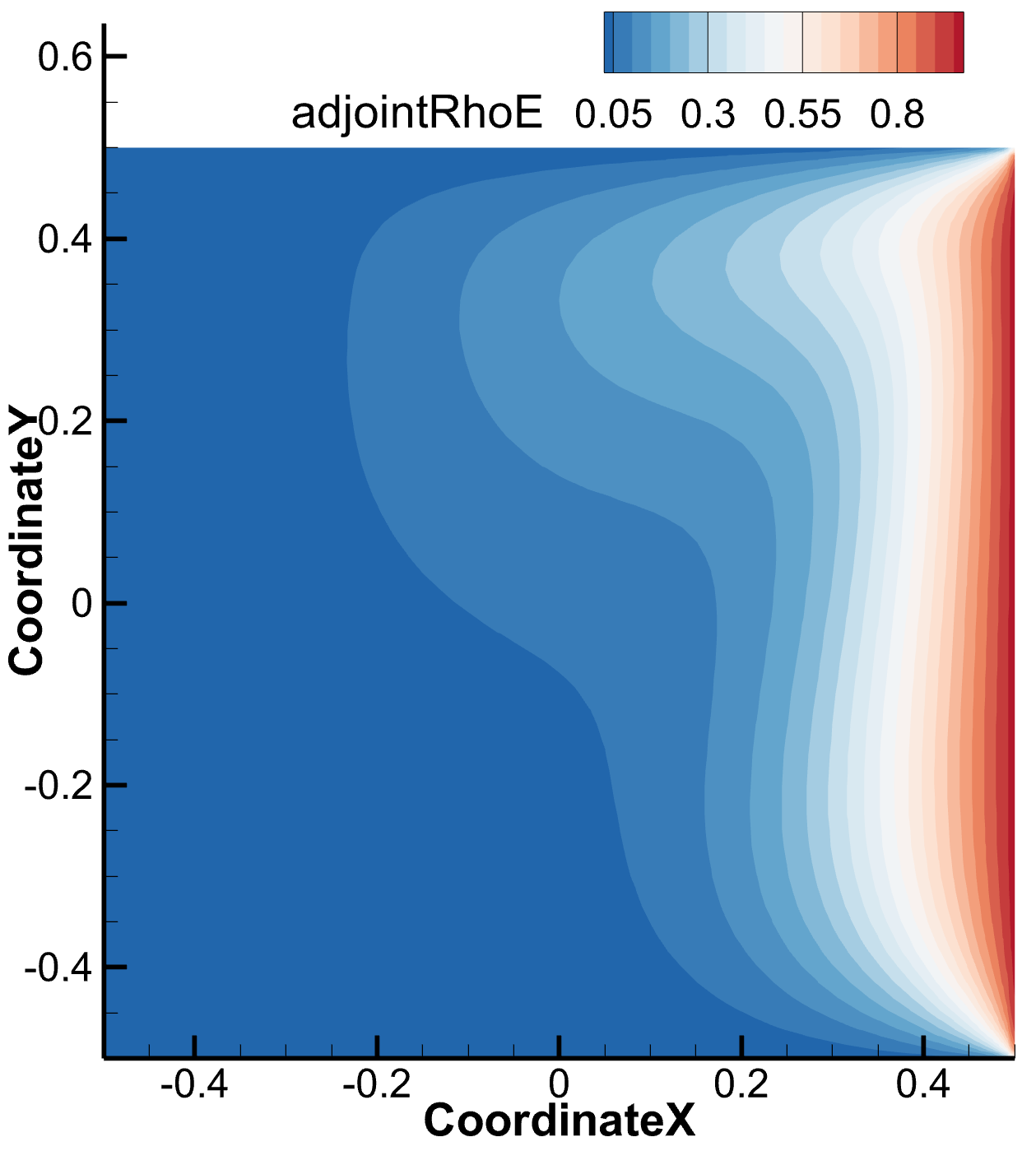}
    \caption{}
    \label{fig:cavityRe100_adjoint}
  \end{subfigure}
  \caption{Field distributions of cavity flow at Reynolds number 100. (a) Baseline temperature field and heat-flux streamlines. (b) Linearized field at perturbation location $y=0.40833$ (contours: perturbed temperature; streamlines: perturbed heat flux). (c) Adjoint energy field.}
  \label{fig:cavityRe100_contour}
\end{figure}

The field distributions for the $Re=100$ case are shown in Fig.~\ref{fig:cavityRe100_contour}. In contrast to the transitional case, the baseline field in Fig.~\ref{fig:cavityRe100_origin} displays heat-flux streamlines that are approximately aligned with the temperature gradient, consistent with Fourier's law in the continuum limit. Figure~\ref{fig:cavityRe100_linear} shows that, at $y=0.40833$, the perturbed thermal response is not emitted along straight paths; instead, it is carried by the cavity flow and organized into the streamline pattern of the baseline field. The adjoint energy field in Fig.~\ref{fig:cavityRe100_adjoint} also exhibits reverse propagation of the sensitivity information from the right wall, but its inward extension is strongly shaped by the flow structure. Compared with the $Kn=0.075$ case, the $Re=100$ fields are therefore influenced to a greater extent by advection associated with the macroscopic flow.

To further verify the present method, two additional regimes are simulated: $Re=1000$ and $Kn=1$.
The $Re=1000$ case employs the same physical mesh as $Re=100$, while the $Kn=1$ case uses the same mesh as $Kn=0.075$.

\begin{figure}[!htbp]
  \centering
  \begin{subfigure}[b]{0.35\textwidth}
    \centering
    \includegraphics[width=\linewidth]{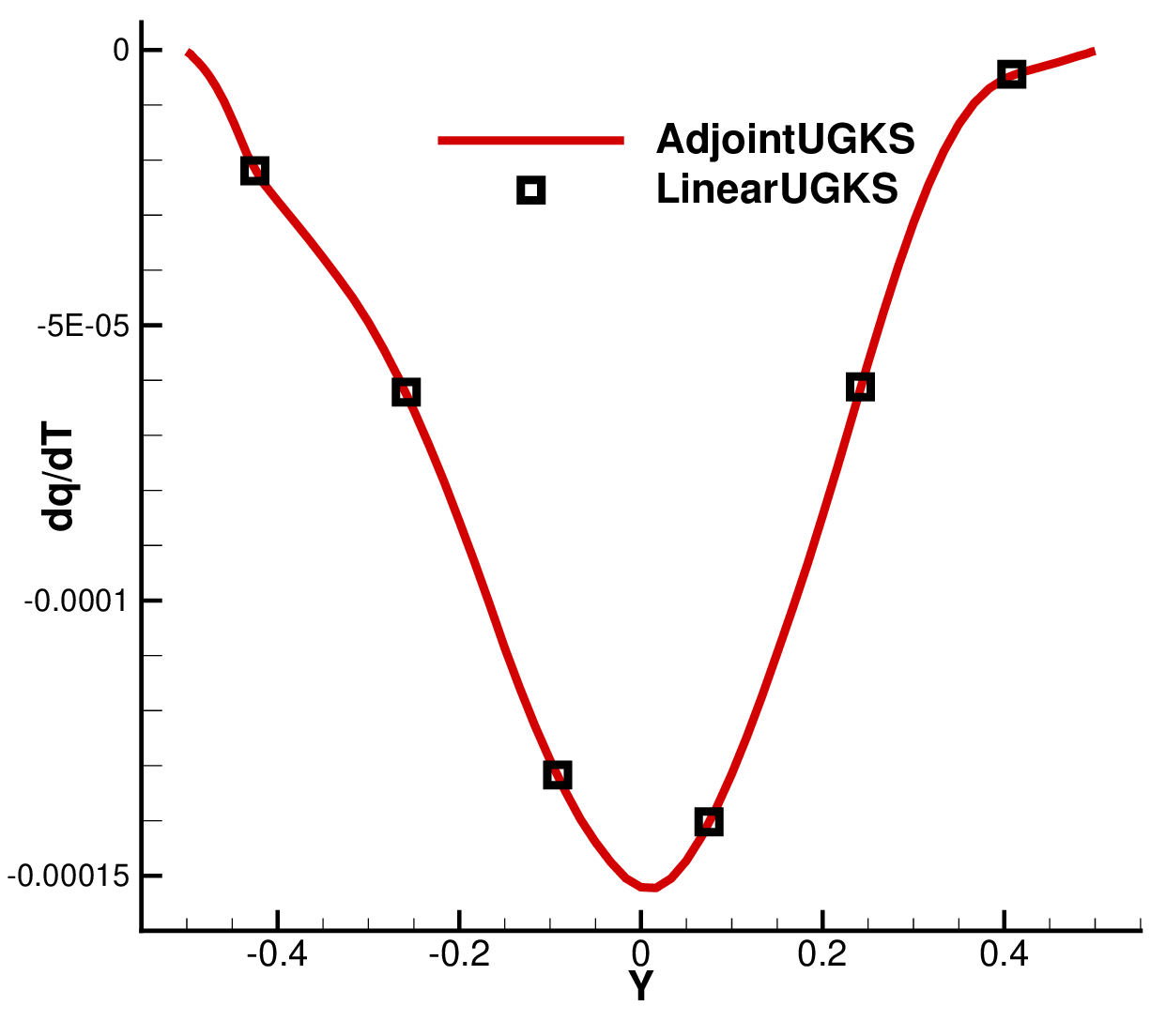}
    \caption{}
    \label{fig:cavityRe1000_sensity}
  \end{subfigure}
  \quad
  \begin{subfigure}[b]{0.35\textwidth}
    \centering
    \includegraphics[width=\linewidth]{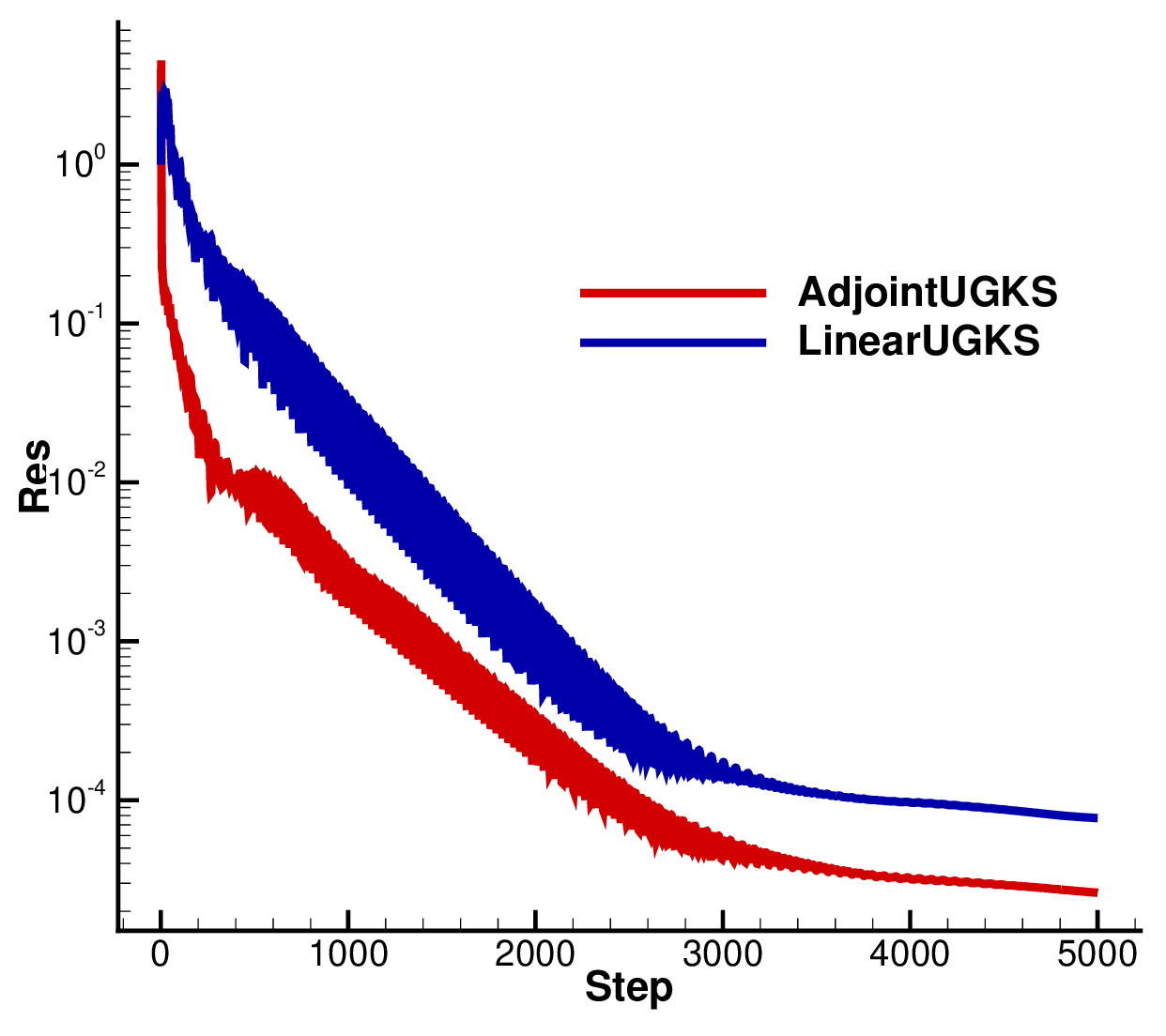}
    \caption{}
    \label{fig:cavityRe1000_residual}
  \end{subfigure}
  \caption{Cavity flow at Reynolds number 1000. (a) Sensitivity to the left wall temperature. (b) Residual of the adjoint system and the linear system.}
  \label{fig:cavityRe1000}
\end{figure}

\begin{figure}[!htbp]
  \centering
  \begin{subfigure}[b]{0.32\textwidth}
    \centering
    \includegraphics[width=\linewidth]{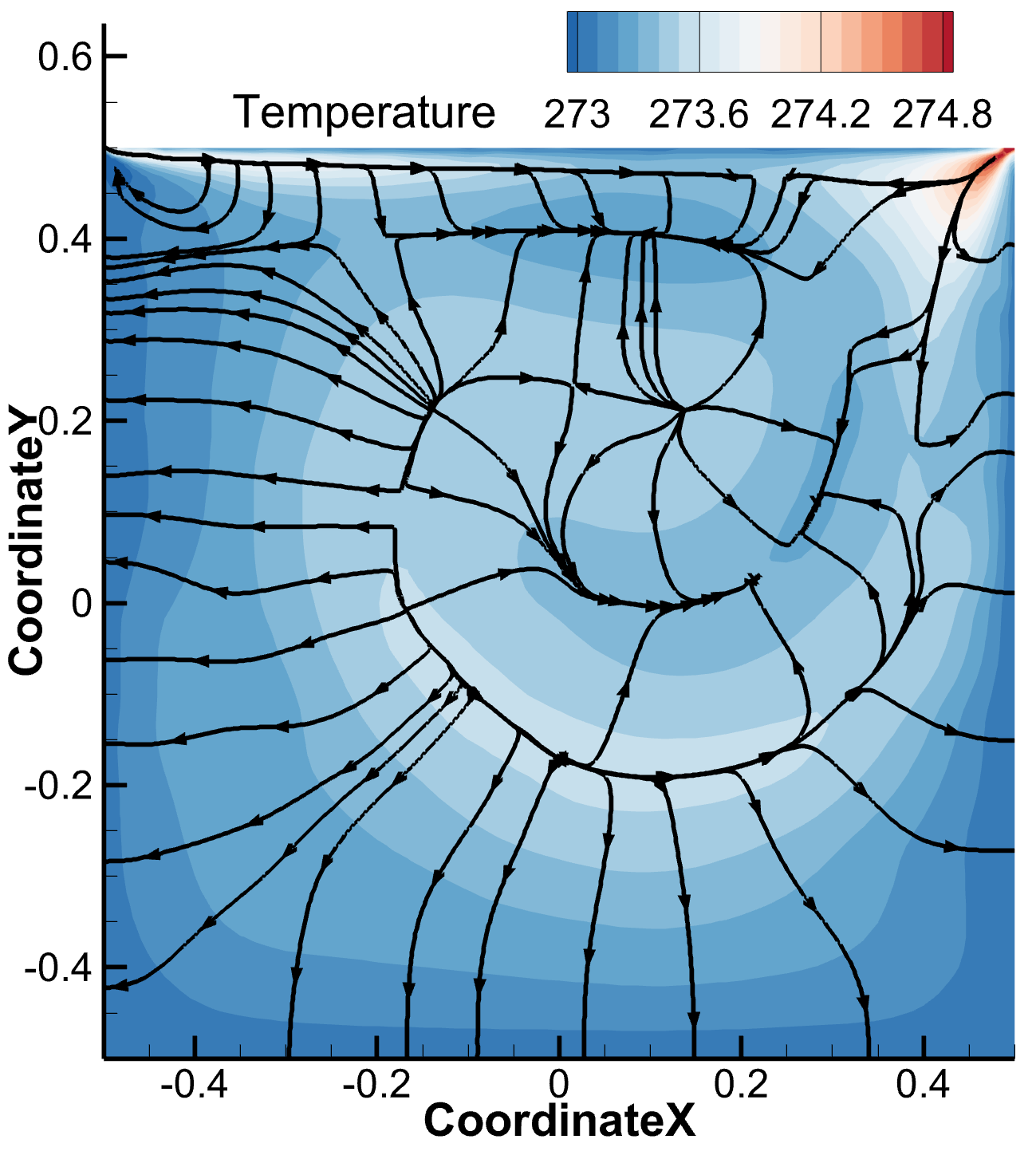}
    \caption{}
    \label{fig:cavityRe1000_origin}
  \end{subfigure}
  \hfill
  \begin{subfigure}[b]{0.32\textwidth}
    \centering
    \includegraphics[width=\linewidth]{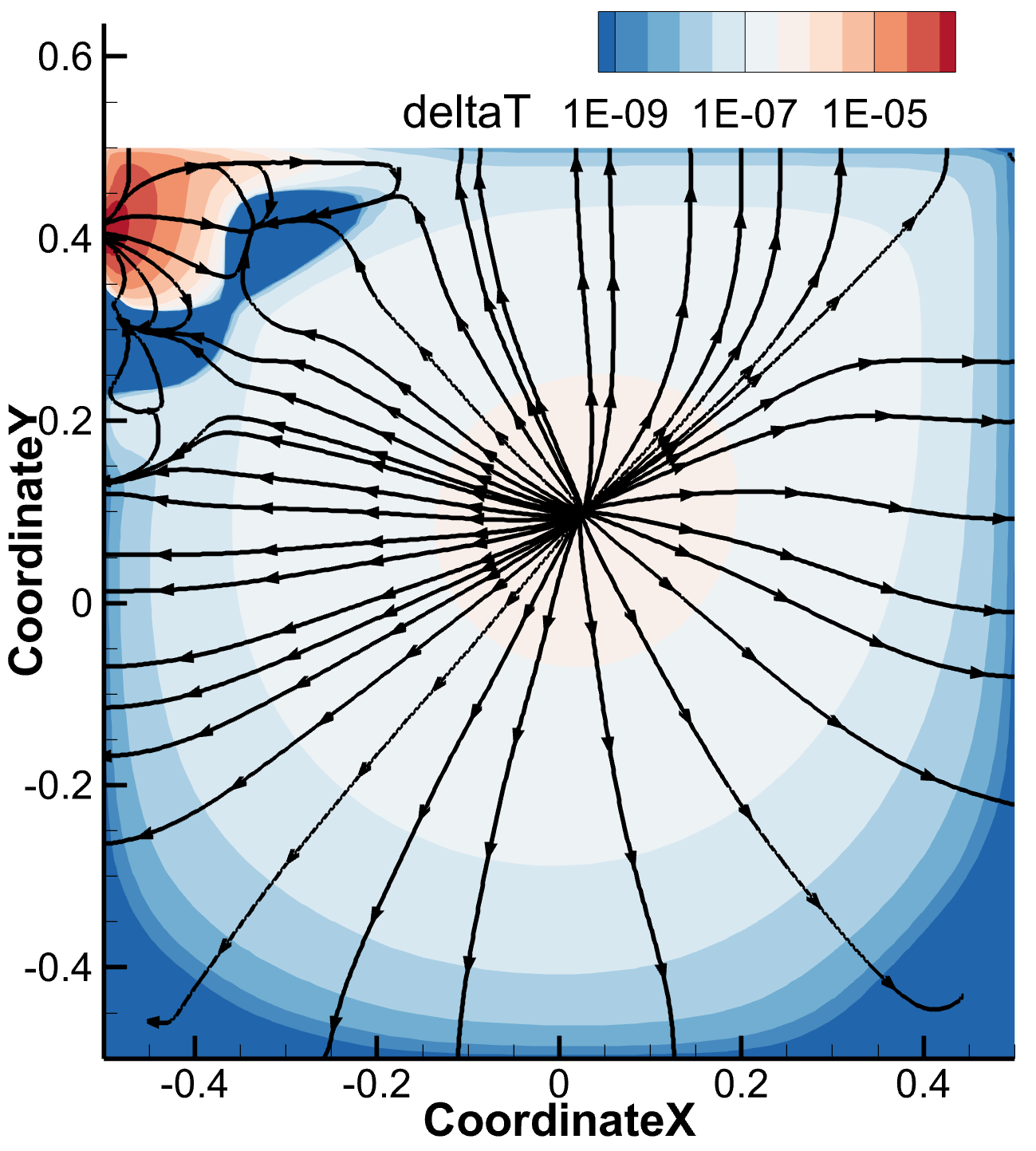}
    \caption{}
    \label{fig:cavityRe1000_linear}
  \end{subfigure}
  \hfill
  \begin{subfigure}[b]{0.32\textwidth}
    \centering
    \includegraphics[width=\linewidth]{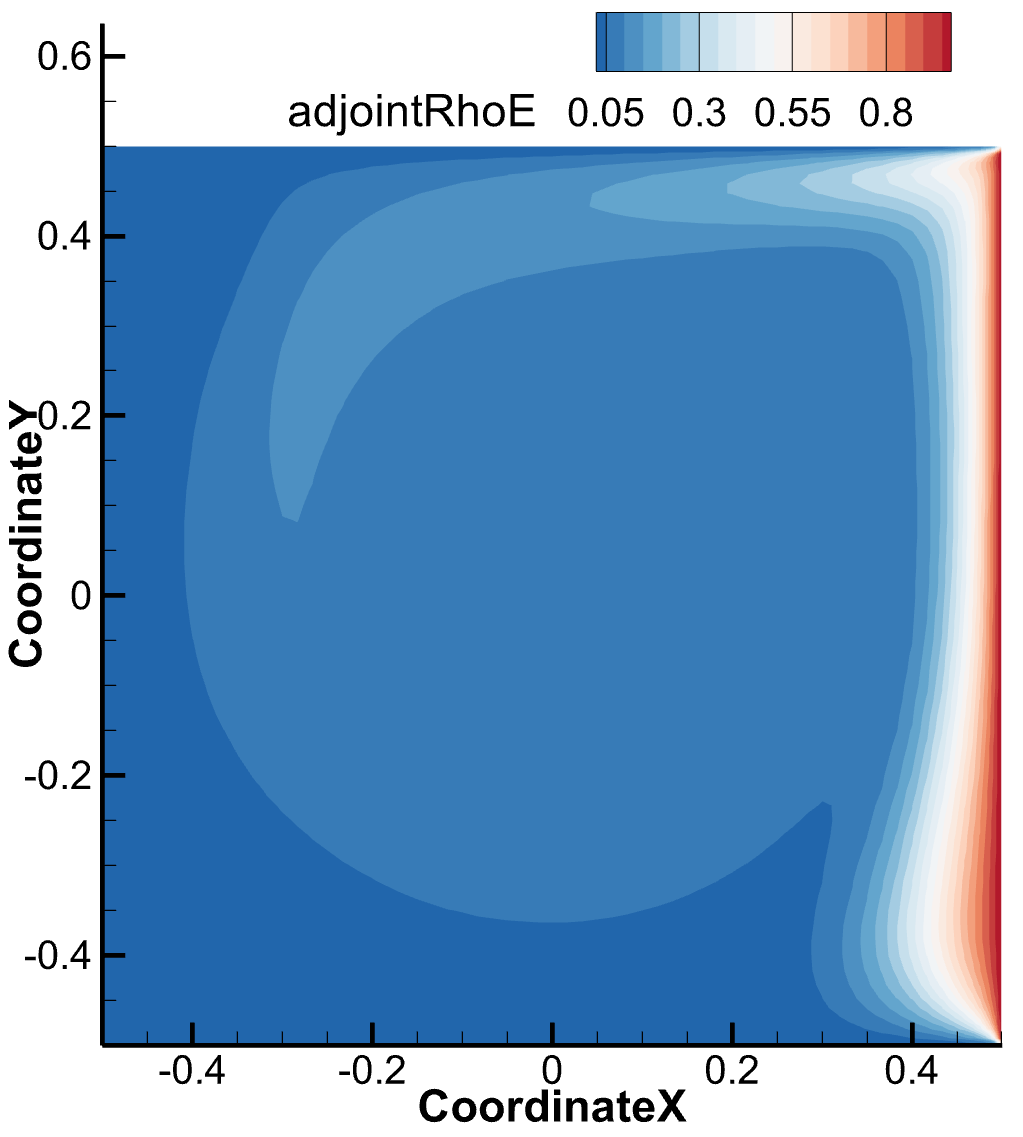}
    \caption{}
    \label{fig:cavityRe1000_adjoint}
  \end{subfigure}
  \caption{Field distributions of cavity flow at Reynolds number 1000. (a) Baseline temperature field and heat-flux streamlines. (b) Linearized field at perturbation location $y=0.40833$ (contours: perturbed temperature; streamlines: perturbed heat flux). (c) Adjoint energy field.}
  \label{fig:cavityRe1000_contour}
\end{figure}

Figure~\ref{fig:cavityRe1000} presents the sensitivity and residual histories for $Re=1000$.
The adjoint and linearized results agree closely in both the sensitivity profiles and the residual convergence curves.
The corresponding field distributions are shown in Fig.~\ref{fig:cavityRe1000_contour}.

\begin{figure}[!htbp]
  \centering
  \begin{subfigure}[b]{0.35\textwidth}
    \centering
    \includegraphics[width=\linewidth]{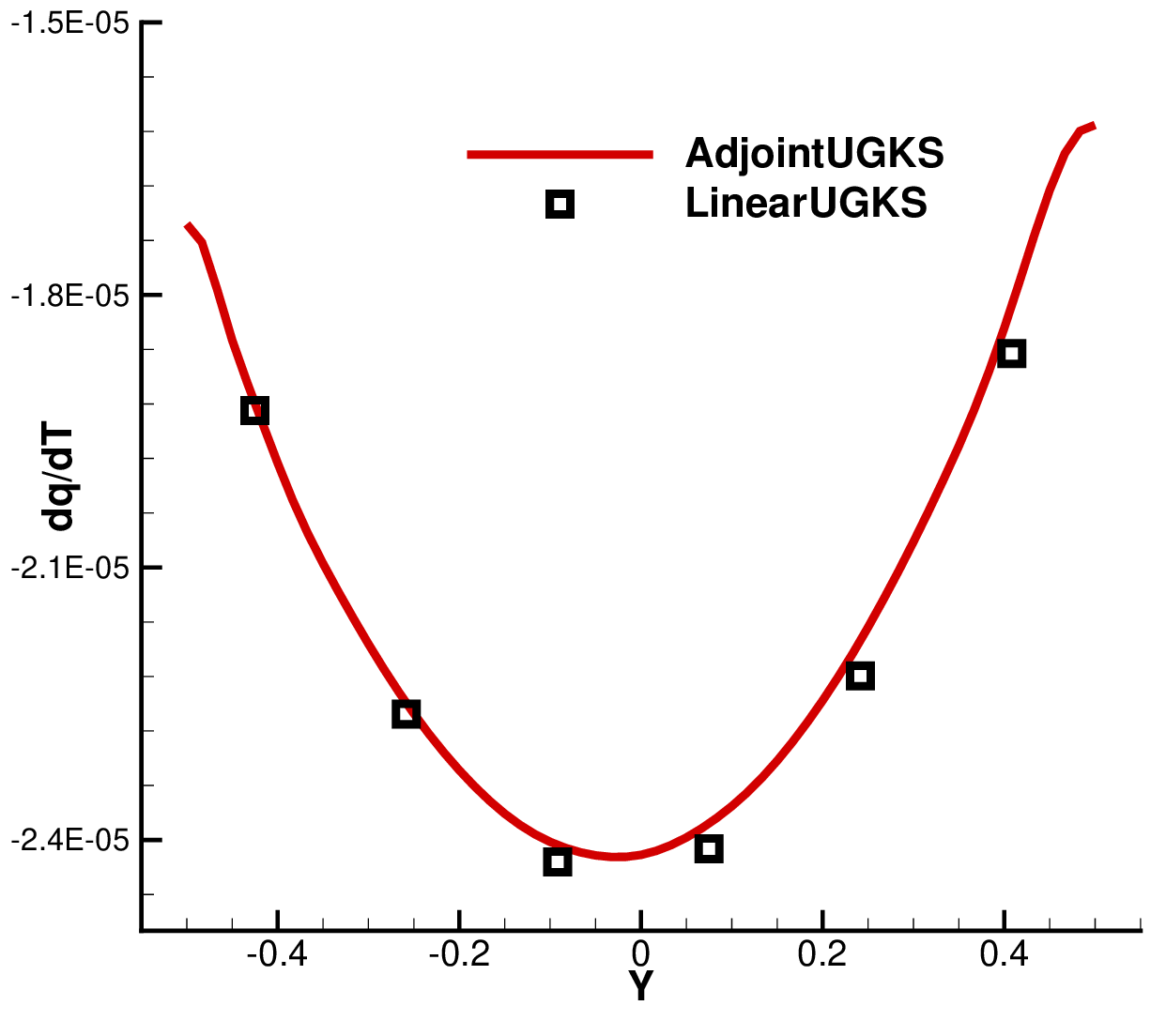}
    \caption{}
    \label{fig:cavityKn1_sensity}
  \end{subfigure}
  \quad
  \begin{subfigure}[b]{0.35\textwidth}
    \centering
    \includegraphics[width=\linewidth]{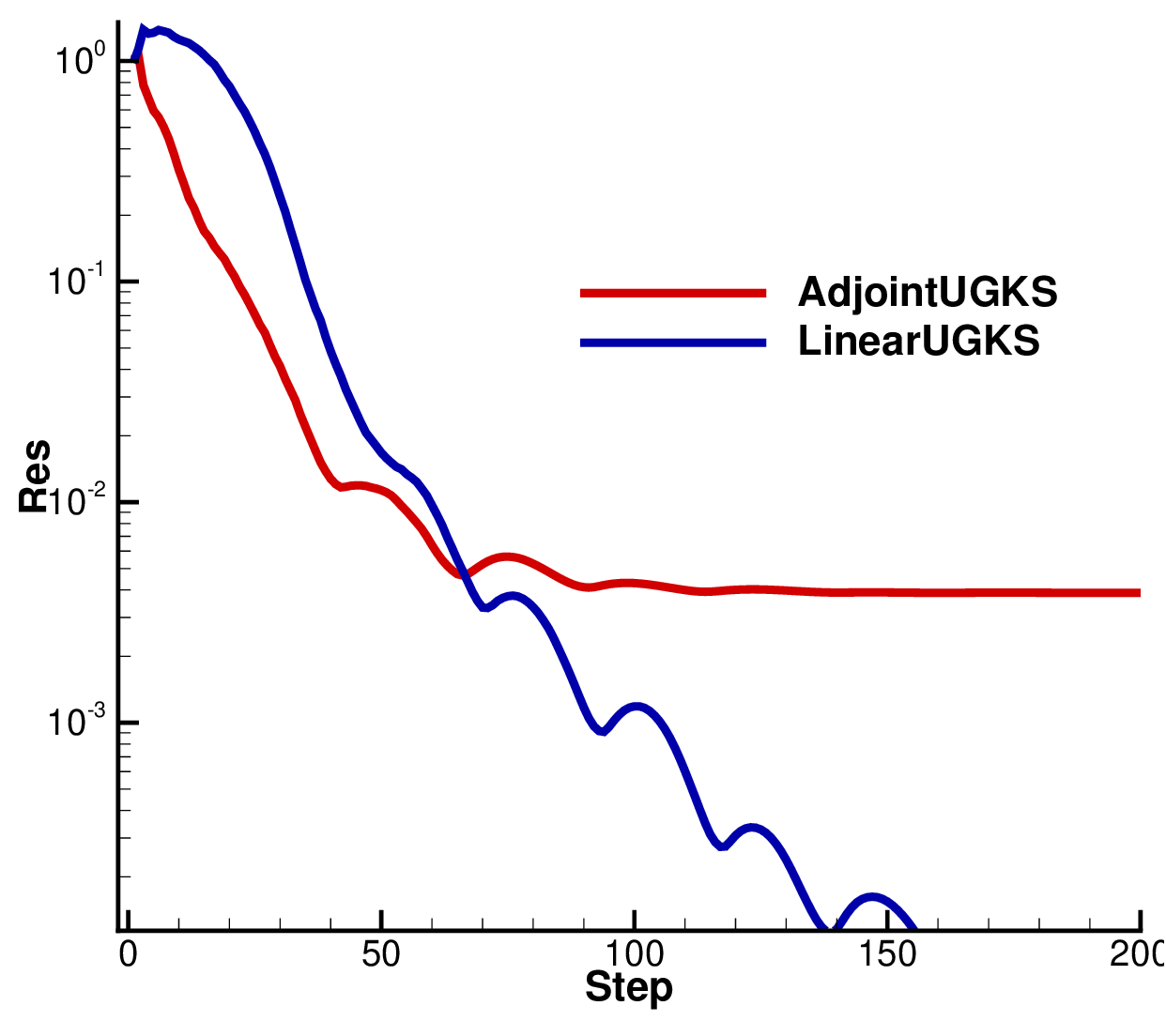}
    \caption{}
    \label{fig:cavityKn1_residual}
  \end{subfigure}
  \caption{Cavity flow at Knudsen number 1. (a) Sensitivity to the left wall temperature. (b) Residual of the adjoint system and the linear system.}
  \label{fig:cavityKn1}
\end{figure}

\begin{figure}[!htbp]
  \centering
  \begin{subfigure}[b]{0.32\textwidth}
    \centering
    \includegraphics[width=\linewidth]{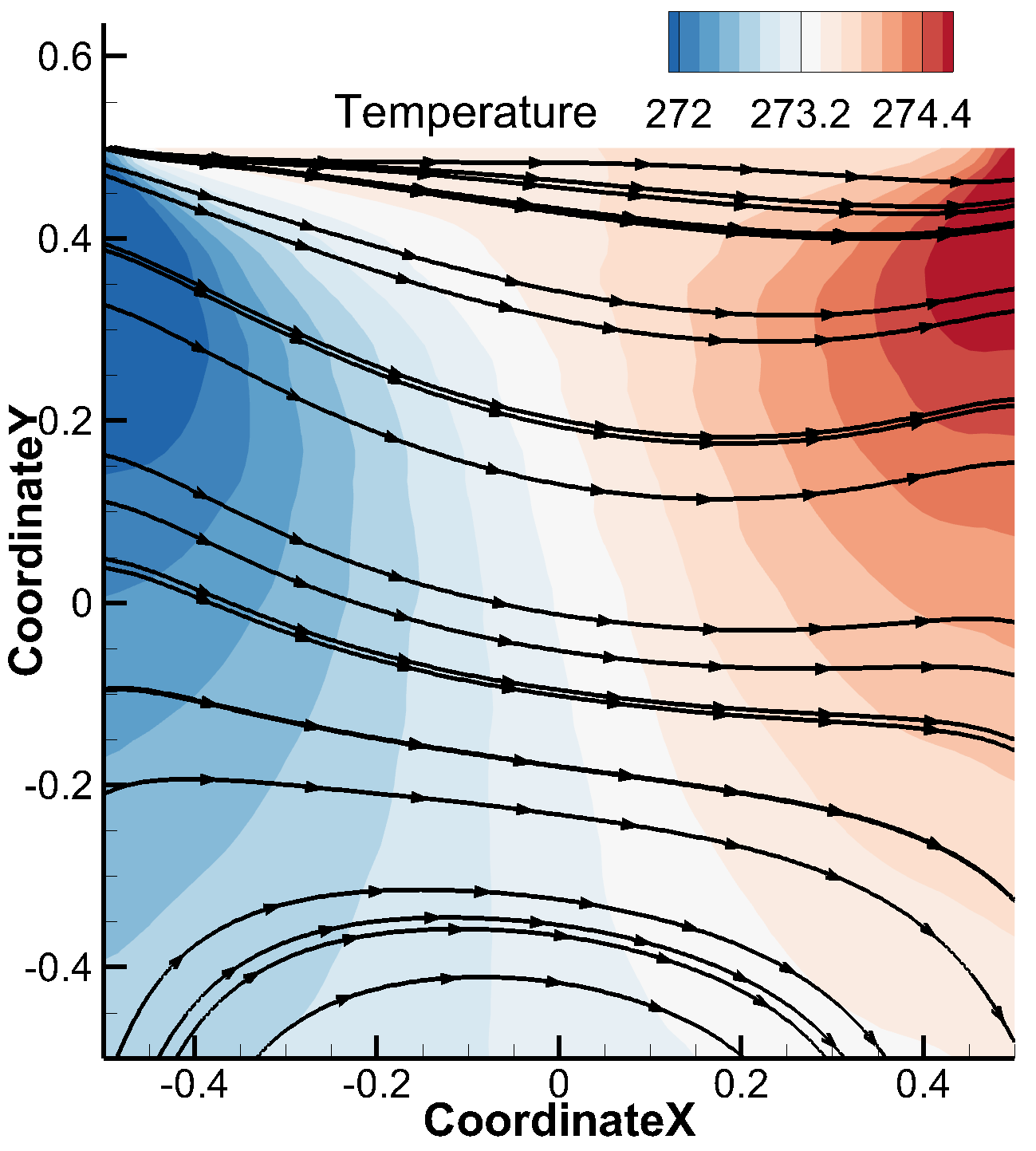}
    \caption{}
    \label{fig:cavityKn1_origin}
  \end{subfigure}
  \hfill
  \begin{subfigure}[b]{0.32\textwidth}
    \centering
    \includegraphics[width=\linewidth]{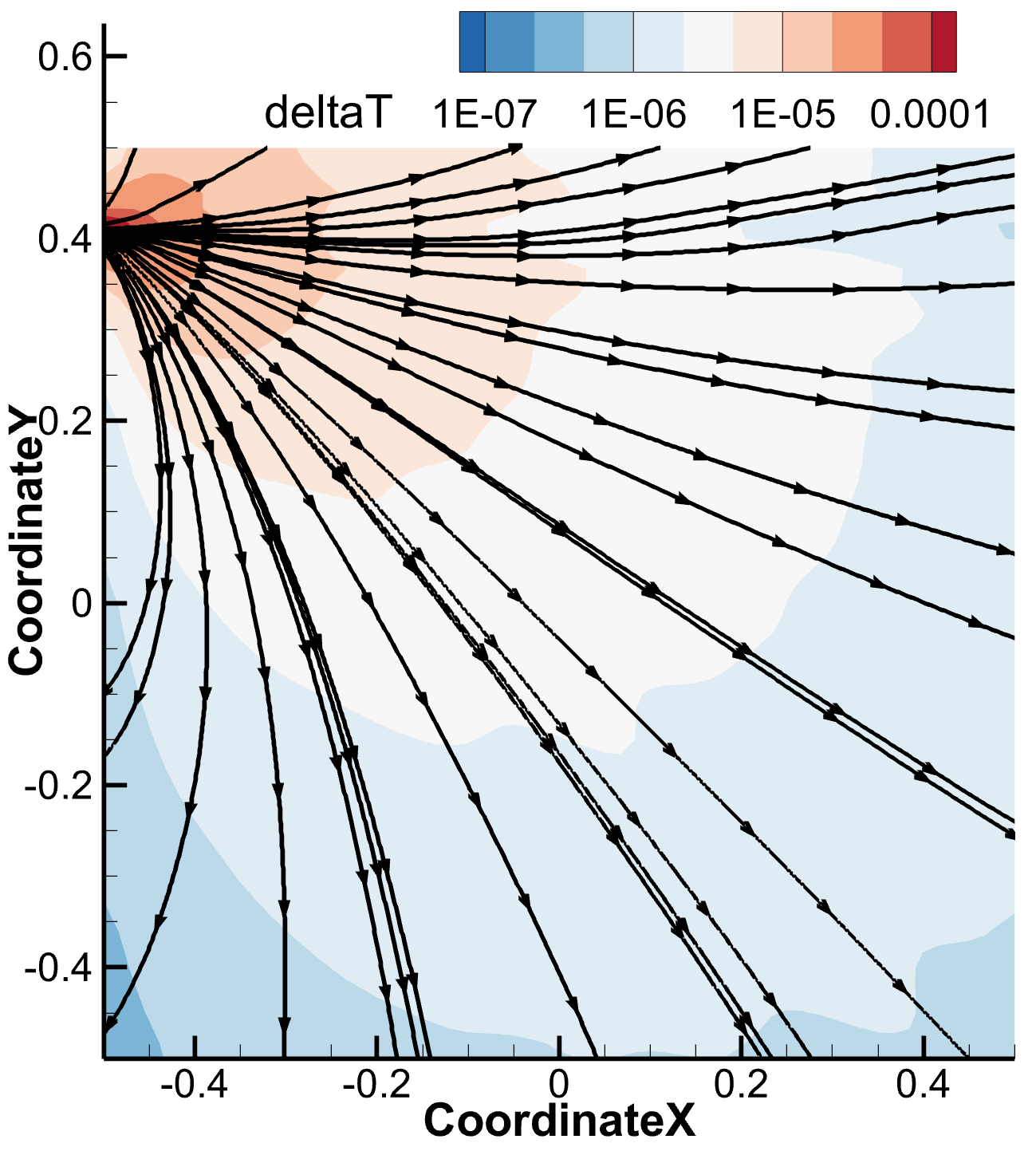}
    \caption{}
    \label{fig:cavityKn1_linear}
  \end{subfigure}
  \hfill
  \begin{subfigure}[b]{0.32\textwidth}
    \centering
    \includegraphics[width=\linewidth]{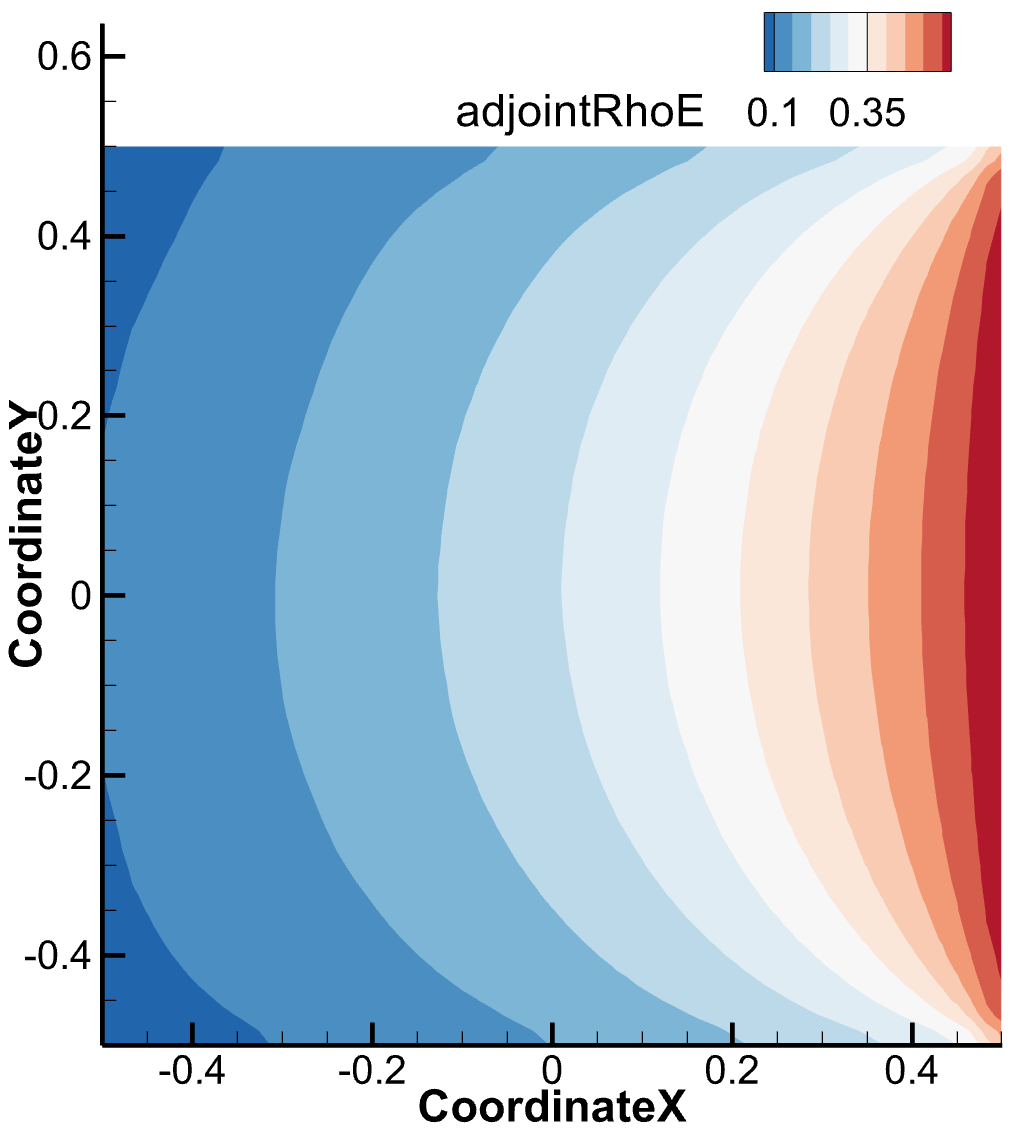}
    \caption{}
    \label{fig:cavityKn1_adjoint}
  \end{subfigure}
  \caption{Field distributions of cavity flow at Knudsen number 1. (a) Baseline temperature field and heat-flux streamlines. (b) Linearized field at perturbation location $y=0.40833$ (contours: perturbed temperature; streamlines: perturbed heat flux). (c) Adjoint energy field.}
  \label{fig:cavityKn1_contour}
\end{figure}

In contrast, the highly rarefied case $Kn=1$ exhibits a noticeable discrepancy between the adjoint and linearized sensitivities (Fig.~\ref{fig:cavityKn1}).
The corresponding field distributions are shown in Fig.~\ref{fig:cavityKn1_contour}.
A minor, mesh-independent discrepancy between the adjoint gradient and the linearized solution is observed at $Kn=1$. This behavior is attributed to the subtle algebraic dual inconsistency at the boundaries under the discrete-velocity framework. Specifically, the linearized boundary operator and its discrete adjoint counterpart do not form a strictly self-adjoint pair at the discrete level for intermediate Knudsen numbers. Furthermore, the convergence of the adjoint residual is limited by the clustering of eigenvalues in the coupled micro-macro system at $Kn=1$, leading to a premature stagnation of the linear solver compared to the forward L-UGKS. Nevertheless, the trend of the adjoint sensitivity curves remains in excellent agreement with the linearized results, verifying the practical utility of the proposed method for design optimization.

\subsection{Thermal creep flow with closed walls}

The second test case considers thermal creep flow in a sealed two-dimensional microchannel with a rectangular cross section suitable for MEMS applications, following the benchmark configuration reported in \cite{masters2007octant}. The two ends of the channel are maintained at different temperatures with $T_1=273\,\mathrm{K}$ and $T_2=573\,\mathrm{K}$, where $T_1<T_2$. The side-wall temperature varies linearly along the channel surface. The working fluid is argon, initially in thermal equilibrium with the walls according to
\begin{equation*}
T(x,y)=\frac{T_2-T_1}{L}x+T_1,
\end{equation*}
and at a uniform pressure of one atmosphere, i.e., $P(x,y)=P=1\,\mathrm{atm}$. Under these conditions, the mean free path of the gas is approximately $64\,\mathrm{nm}$. Two channel widths are considered: $1\,\mu\mathrm{m}$ and $20\,\mathrm{nm}$. The channel is discretized using $200$ cells along the length and $40$ cells along the width. The velocity space is discretized with a $60\times60$ grid over the domain $(-5\sqrt{2RT_2},5\sqrt{2RT_2})^2$.

\begin{figure}[!htbp]
  \centering
  \begin{subfigure}[b]{0.35\textwidth}
    \centering
    \includegraphics[width=\linewidth]{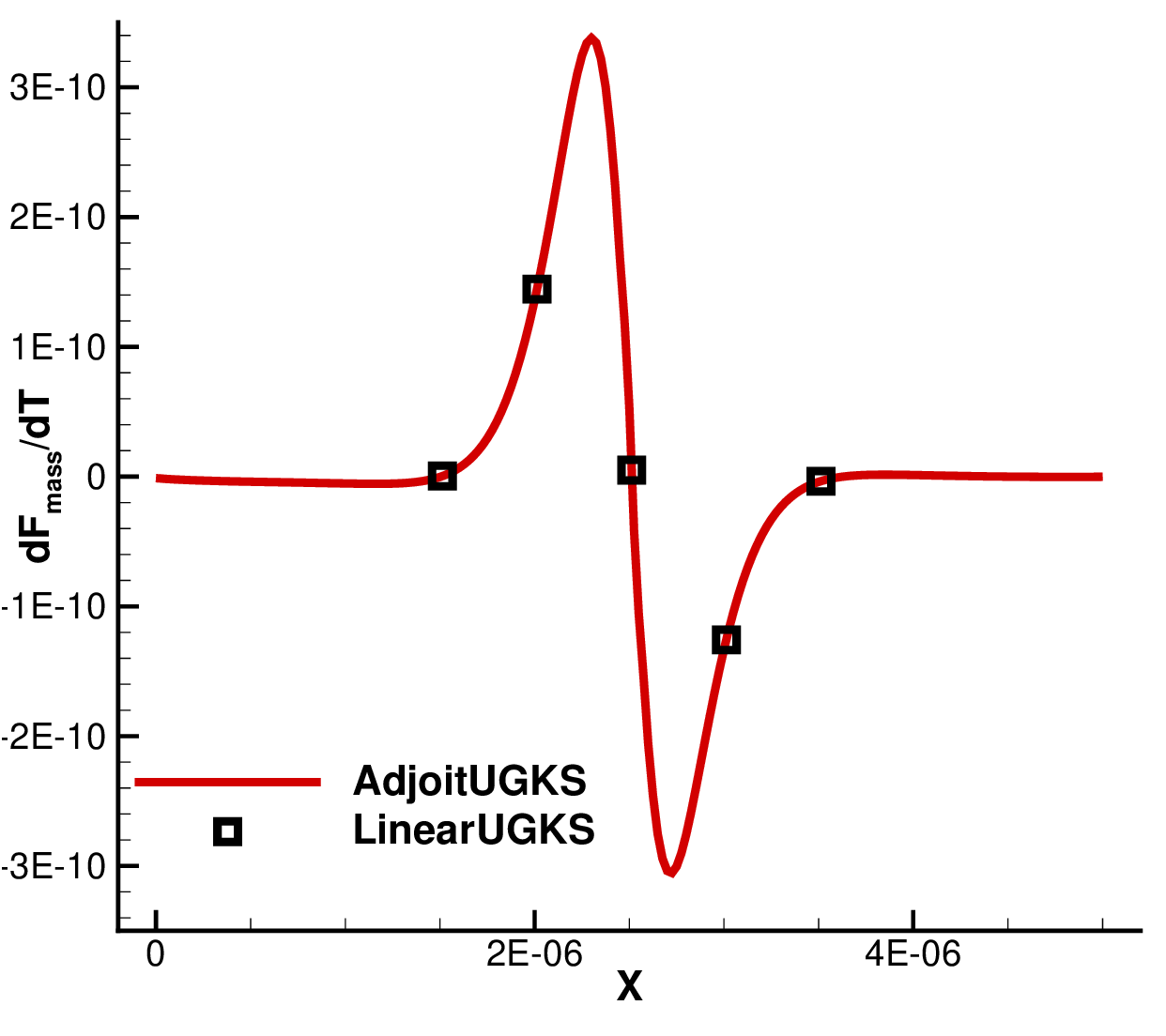}
    \caption{}
    \label{fig:microFlow1e6_sensity}
  \end{subfigure}
  \quad
  \begin{subfigure}[b]{0.35\textwidth}
    \centering
    \includegraphics[width=\linewidth]{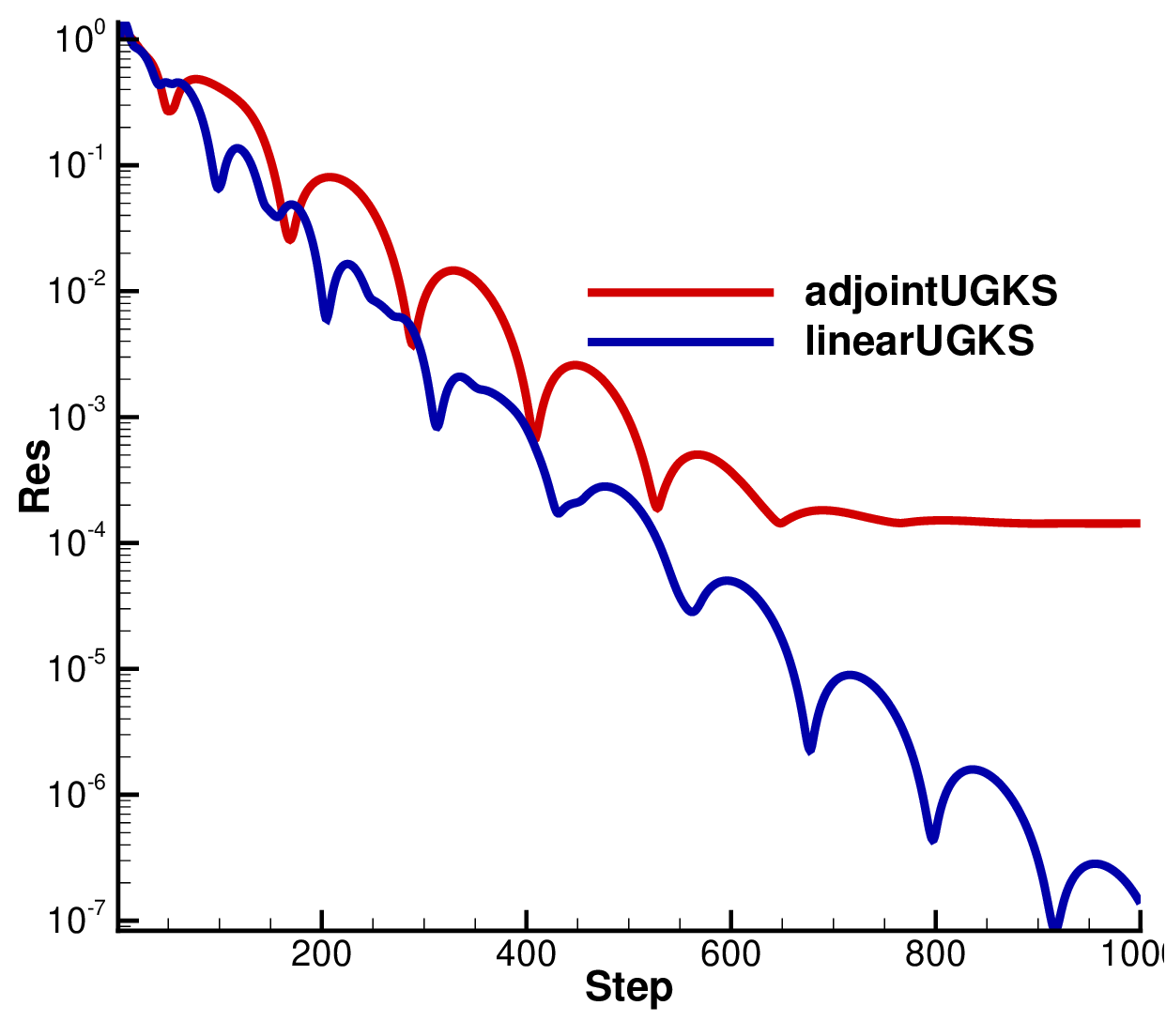}
    \caption{}
    \label{fig:microFlow1e6_residual}
  \end{subfigure}
  \caption{Thermal creep flow with closed walls with channel width 1$\mu$m. (a) Sensitivity with respect to the upper- and lower-wall temperatures. (b) Residual of the adjoint system and the linear system.}
  \label{fig:microFlow1e6}
\end{figure}

\begin{figure}[!htbp]
  \centering
  \begin{subfigure}[b]{0.9\textwidth}
    \centering
    \includegraphics[width=\linewidth]{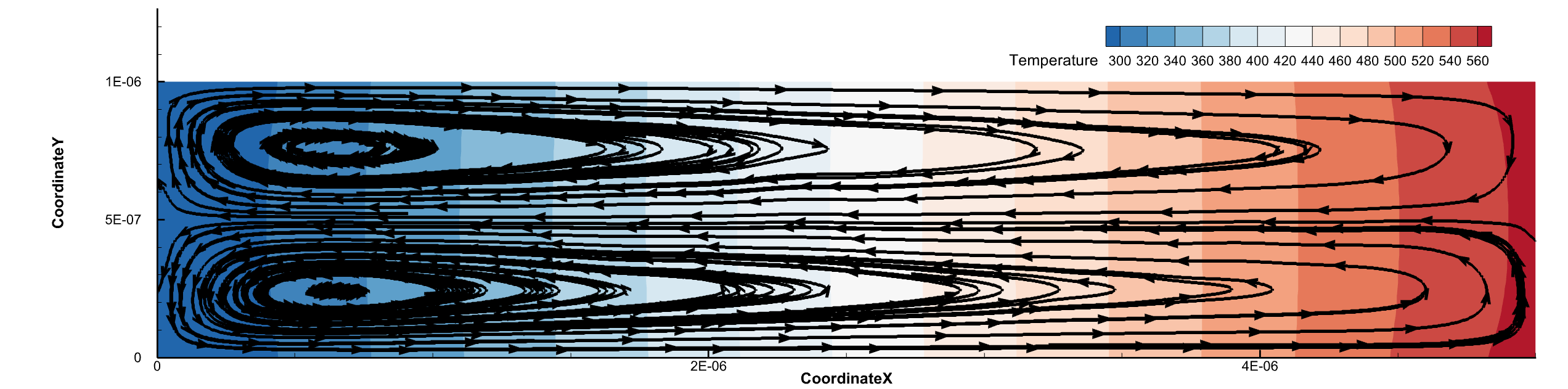}
    \caption{}
    \label{fig:microFlow1e6_origin}
  \end{subfigure}
  \hfill
  \begin{subfigure}[b]{0.9\textwidth}
    \centering
    \includegraphics[width=\linewidth]{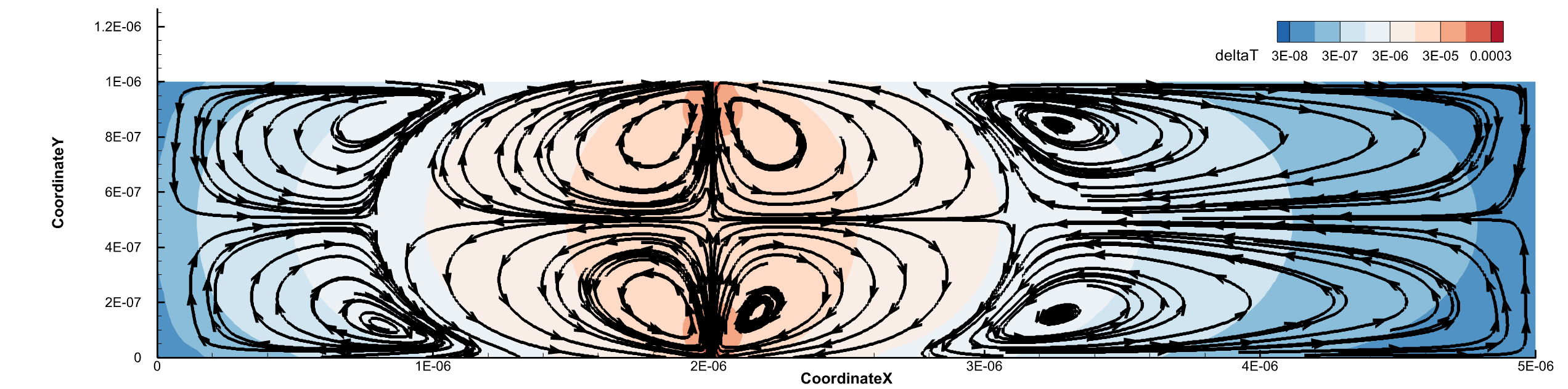}
    \caption{}
    \label{fig:microFlow1e6_linear}
  \end{subfigure}
  \hfill
  \begin{subfigure}[b]{0.9\textwidth}
    \centering
    \includegraphics[width=\linewidth]{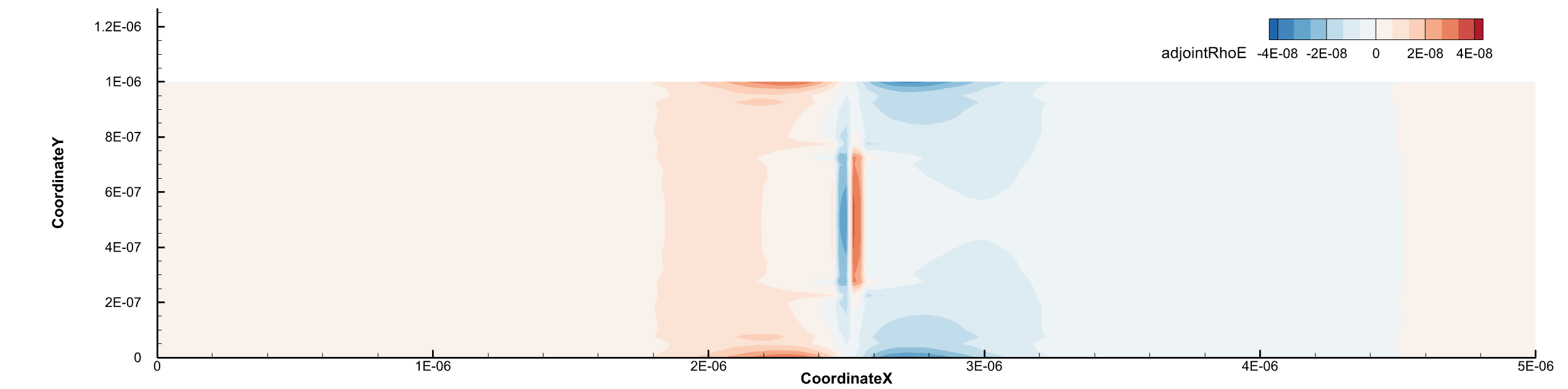}
    \caption{}
    \label{fig:microFlow1e6_adjoint}
  \end{subfigure}
  \caption{Field distributions of thermal creep flow with closed walls with channel width 1$\mu$m. (a) Baseline temperature field and streamlines. (b) Linearized field at perturbation location (contours: perturbed temperature; streamlines: perturbed streamlines). (c) Adjoint $x$-direction moment field.}
  \label{fig:microFlow1e6_contour}
\end{figure}

The objective function is defined as the mass flow rate through the vertical cross-section at the streamwise center $x=L/2$, integrated along the channel height from $y=H/4$ to $y=3H/4$, where $L$ is the channel length and $H$ is the channel width. In other words, $(\rho u)$ is evaluated on the mid-channel plane and integrated over the central half of the channel height:
\begin{equation*}
J = \int_{H/4}^{3H/4} (\rho u)\big|_{x=L/2} \,\mathrm{d}y.
\end{equation*}
The sensitivities of $J$ with respect to the upper- and lower-wall temperatures are investigated. A linearized solver is used for comparison.
For the linearized program, temperature perturbations are imposed on the upper and lower walls with $\delta T=1\times 10^{-3}\,\mathrm{K}$. The streamwise locations of the perturbation points are
\begin{equation*}
x_n=\frac{1}{40}\left(20n+\frac{1}{2}+60\right), \quad n=0,1,2,3,4.
\end{equation*}
Figure~\ref{fig:microFlow1e6} shows the sensitivities of the mass flow rate with respect to the upper- and lower-wall temperatures for the $1\,\mu\mathrm{m}$ channel width. The adjoint and linearized sensitivities agree closely, and the residual histories converge at essentially the same rate. This confirms that the adjoint formulation remains valid in the near-continuum regime.
The field distributions for the $1\,\mu\mathrm{m}$ channel are shown in Fig.~\ref{fig:microFlow1e6_contour}. The baseline temperature field and streamlines in Fig.~\ref{fig:microFlow1e6_origin} exhibit a counter-rotating vortex pair driven by the wall temperature gradient. Figure~\ref{fig:microFlow1e6_linear} presents the linearized field at the perturbation location, where the perturbed flow emanates from the wall and organizes into the vortex-pair structure of the baseline state, reflecting continuum-like transport behavior. The adjoint field in Fig.~\ref{fig:microFlow1e6_adjoint} displays the corresponding sensitivity response associated with the mass-flow objective.

\begin{figure}[!htbp]
  \centering
  \begin{subfigure}[b]{0.35\textwidth}
    \centering
    \includegraphics[width=\linewidth]{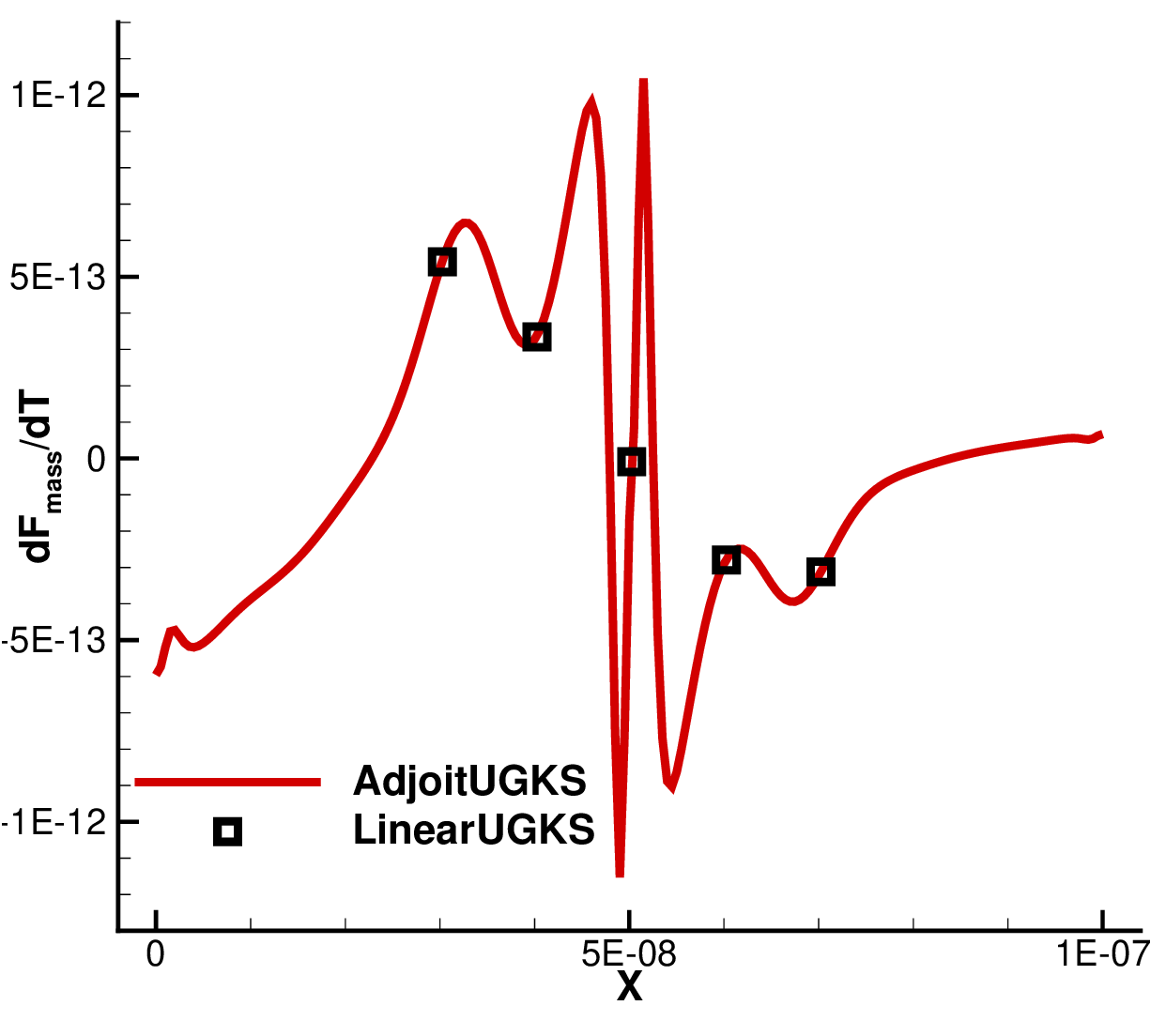}
    \caption{}
    \label{fig:microFlow2e8_sensity}
  \end{subfigure}
  \quad
  \begin{subfigure}[b]{0.35\textwidth}
    \centering
    \includegraphics[width=\linewidth]{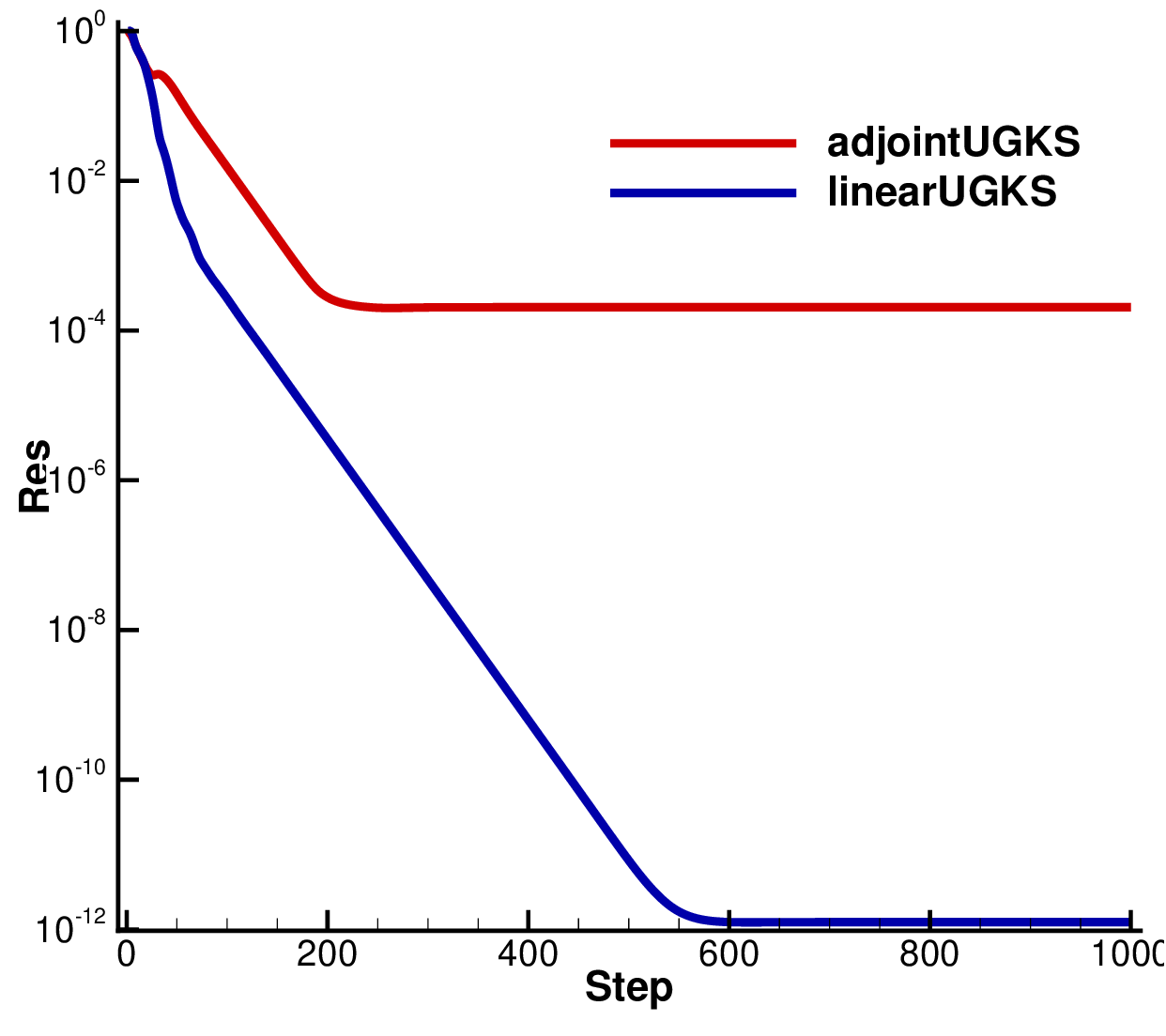}
    \caption{}
    \label{fig:microFlow2e8_residual}
  \end{subfigure}
  \caption{Thermal creep flow with closed walls with channel width 20\,nm. (a) Sensitivity with respect to the upper- and lower-wall temperatures. (b) Residual of adjoint system and linear system.}
  \label{fig:microFlow2e8}
\end{figure}

Figure~\ref{fig:microFlow2e8} presents the corresponding results for the $20\,\mathrm{nm}$ channel width. Again, the adjoint and linearized sensitivities agree closely, and the residual histories converge at essentially the same rate.

\begin{figure}[!htbp]
  \centering
  \begin{subfigure}[b]{0.9\textwidth}
    \centering
    \includegraphics[width=\linewidth]{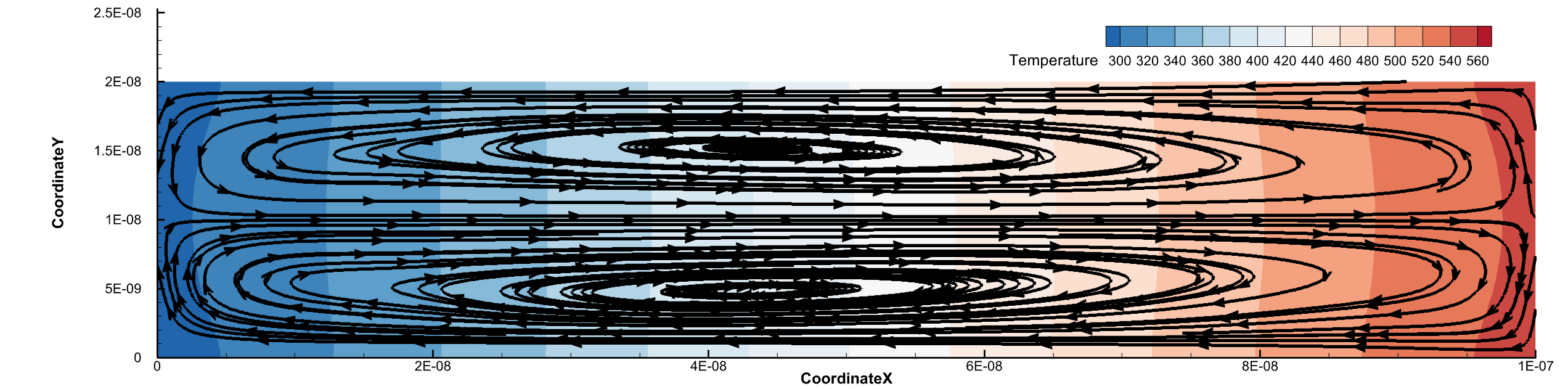}
    \caption{}
    \label{fig:microFlow2e8_origin}
  \end{subfigure}
  \hfill
  \begin{subfigure}[b]{0.9\textwidth}
    \centering
    \includegraphics[width=\linewidth]{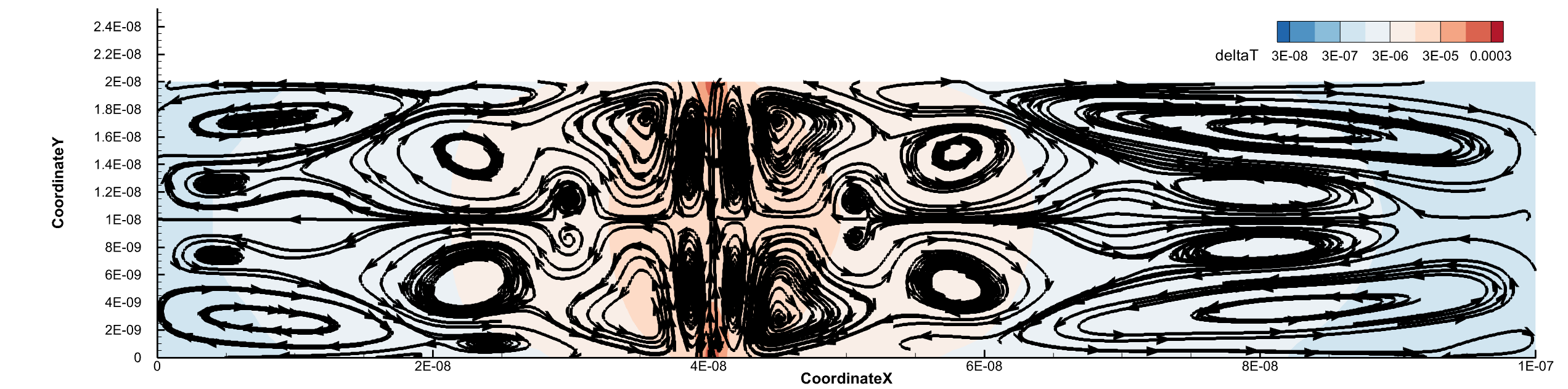}
    \caption{}
    \label{fig:microFlow2e8_linear}
  \end{subfigure}
  \hfill
  \begin{subfigure}[b]{0.9\textwidth}
    \centering
    \includegraphics[width=\linewidth]{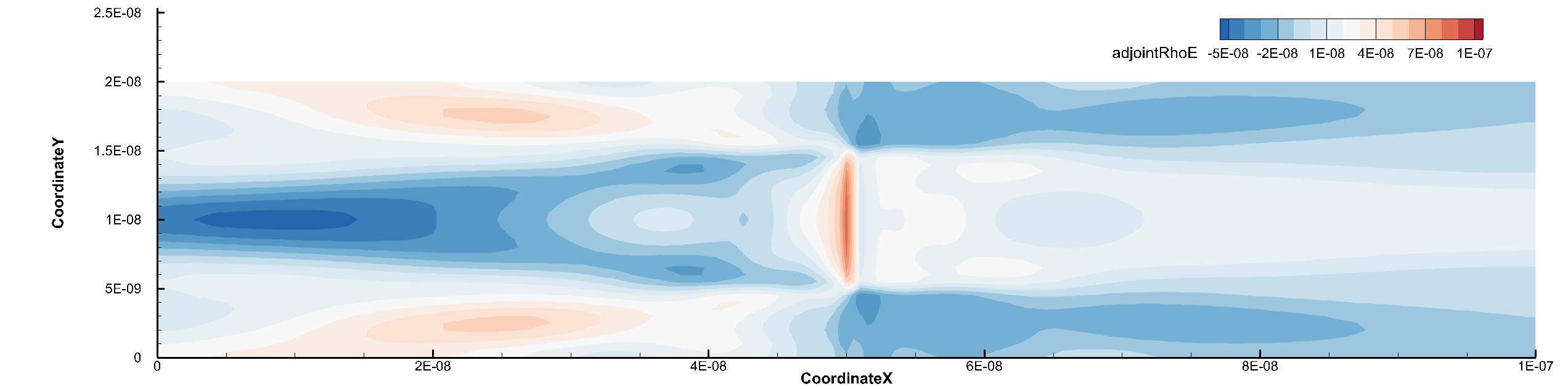}
    \caption{}
    \label{fig:microFlow2e8_adjoint}
  \end{subfigure}
  \caption{Field distributions of thermal creep flow with closed walls with channel width 20\,nm. (a) Baseline temperature field and streamlines. (b) Linearized field at perturbation location (contours: perturbed temperature; streamlines: perturbed streamlines). (c) Adjoint $x$-direction moment field.}
  \label{fig:microFlow2e8_contour}
\end{figure}

Unlike the $1\,\mu\mathrm{m}$ case, the baseline vortex pair in Fig.~\ref{fig:microFlow2e8_origin} rotates in the opposite direction. In the linearized field in Fig.~\ref{fig:microFlow2e8_linear}, the perturbed flow emanates from the wall and a large fraction of the gas streams directly toward the channel center; the upper and lower streams collide near the midline, producing a richer vortex structure characteristic of stronger rarefaction effects. The adjoint field in Fig.~\ref{fig:microFlow2e8_adjoint} exhibits a correspondingly more intricate sensitivity pattern.

\subsection{Hypersonic cylinder flow with surface roughness}

\begin{figure}
\centering
	 \begin{subfigure}[b]{0.4\textwidth}
    \centering
    \includegraphics[width=\linewidth]{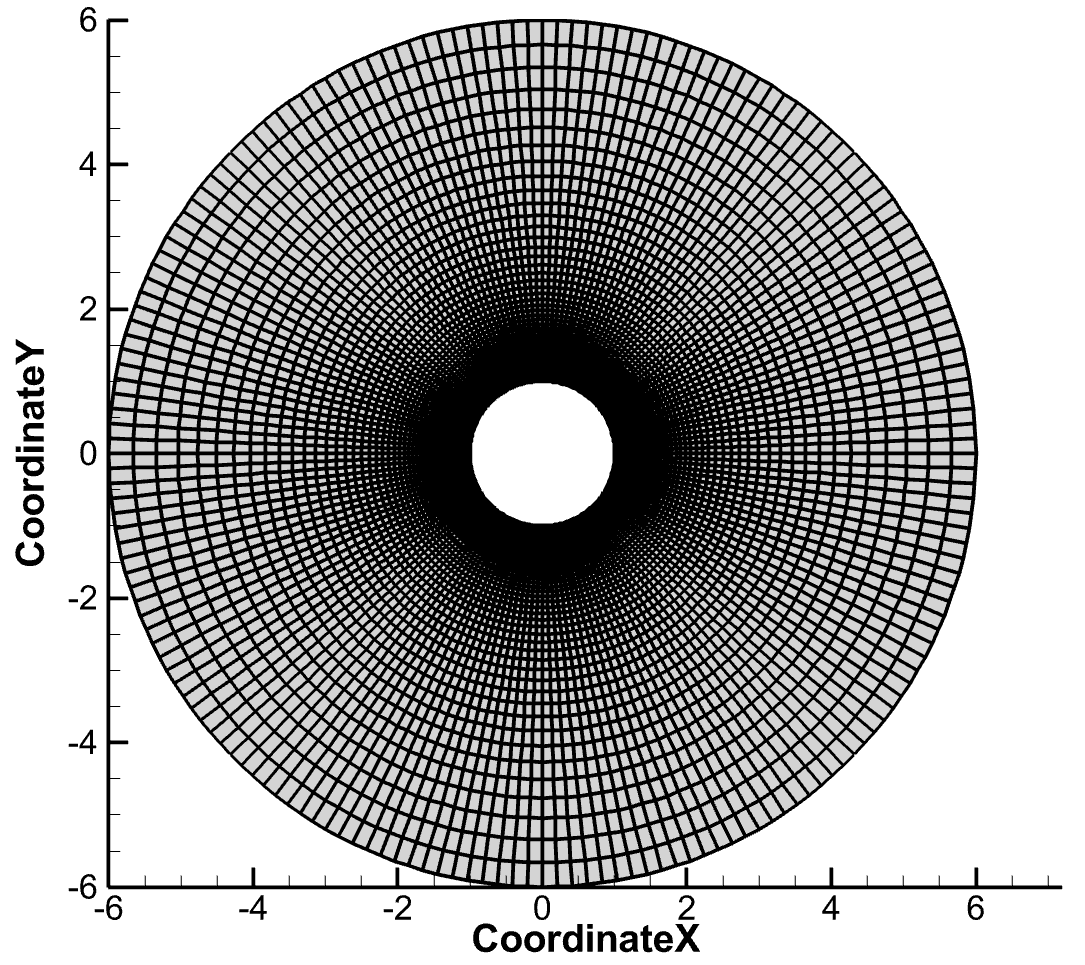}
    \caption{}
    \label{fig:cylindermesh}
  \end{subfigure}
  \quad
  \begin{subfigure}[b]{0.4\textwidth}
    \centering
    \includegraphics[width=\linewidth]{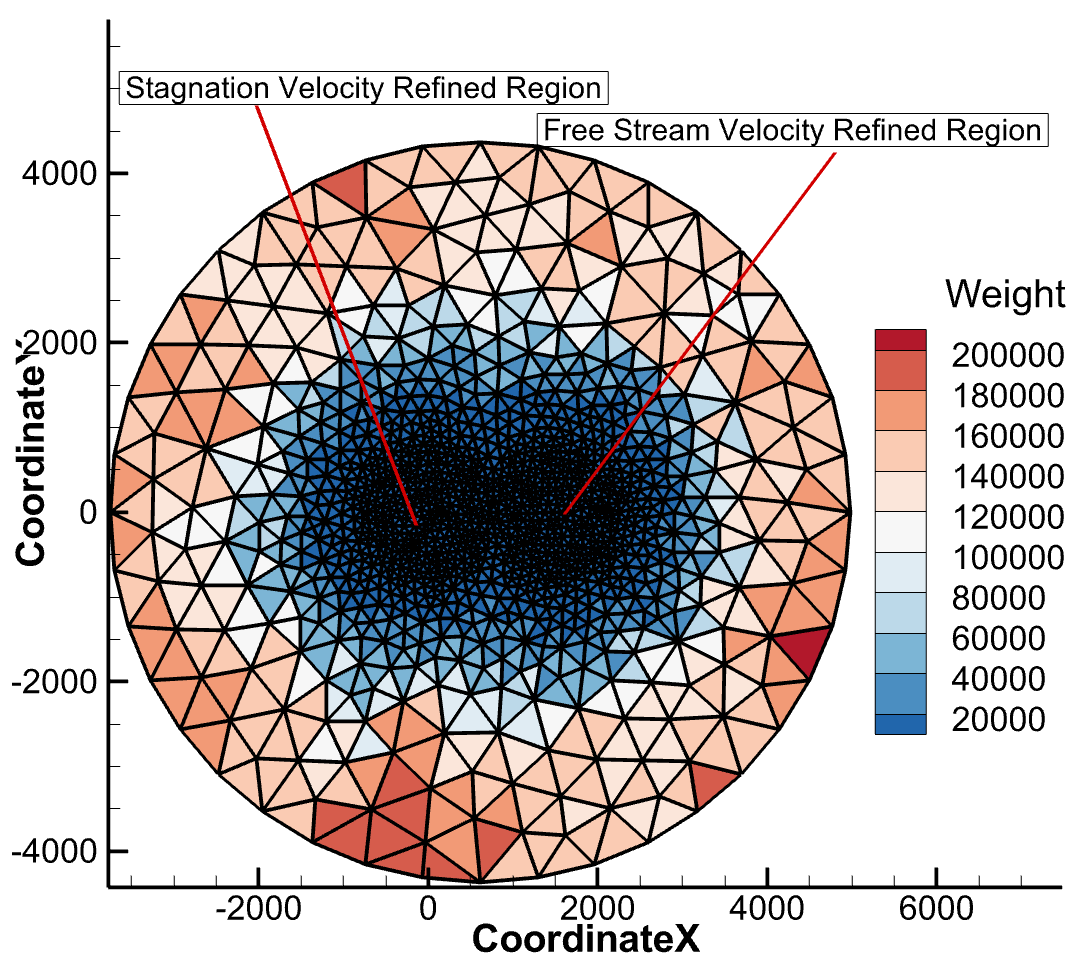}
    \caption{}
    \label{fig:cylinderdvs}
  \end{subfigure}
  \caption{Mesh and DVS mesh for the cylinder flow. (a) Physical space mesh discretized by a $180\times70$ mesh with first wall-cell thickness $0.005$. (b) DVS mesh consists of $2210$ cells.}
  \label{fig:cylinder_mesh}
\end{figure}

The third test case considers hypersonic rarefied flow past a circular cylinder and investigates the influence of surface roughness on the wall heat flux.
The freestream temperature is $T_\infty=273\,\mathrm{K}$, the freestream Mach number is $Ma=5.0$, and the Prandtl number is $Pr=1$.
The wall temperature is fixed at $T_w=273\,\mathrm{K}$, and the solid surface is treated by the Maxwell reflection model with a baseline accommodation coefficient $\alpha_0=0.8$.
The physical space is discretized by a $180\times70$ mesh, and the first wall-cell thickness is $0.005$ as shown in Fig.~\ref{fig:cylindermesh}, while the unstructured discrete velocity space (DVS) mesh consists of $2210$ cells as shown in Fig.~\ref{fig:cylinderdvs}.
The DVS is discretized in a circular region centered at $0.4\times(U_\infty,0)$ with a total radius of $6\sqrt{RT_s}$, where $T_s$ is the stagnation temperature of the freestream flow.
The unstructured DVS mesh is refined at the zero-velocity point within a radius of $3\sqrt{RT_w}$ and at the freestream point within a radius of $3\sqrt{RT_\infty}$.
Three freestream Knudsen numbers based on the cylinder diameter are examined, $Kn=0.4$, $0.1$, and $0.02$, spanning transitional to near-continuum regimes.
Surface roughness is modeled by a circumferential variation of the wall reflection coefficient,
\begin{equation*}
  \alpha(\theta)=\alpha_0+\alpha^\prime, \qquad \alpha^\prime=0.001\cos(2n\theta), \qquad \theta=\arctan(x/|y|),
\end{equation*}
where $\theta$ is the polar angle measured from the stagnation point and $n$ is the roughness mode index.
For each mode, the perturbation amplitude is fixed at $0.001$, and the corresponding heat-flux response is examined mode by mode.
The objective functional $J$ is the surface heat flux on the cylinder wall.
Sensitivities of $J$ with respect to the perturbation coefficient in $\alpha^\prime=a\cos(2n\theta)$ are computed by the adjoint UGKS and compared with an independent linearized UGKS solver.

\begin{figure}[!htbp]
  \centering
  \begin{subfigure}[b]{0.32\textwidth}
    \centering
    \includegraphics[width=\linewidth]{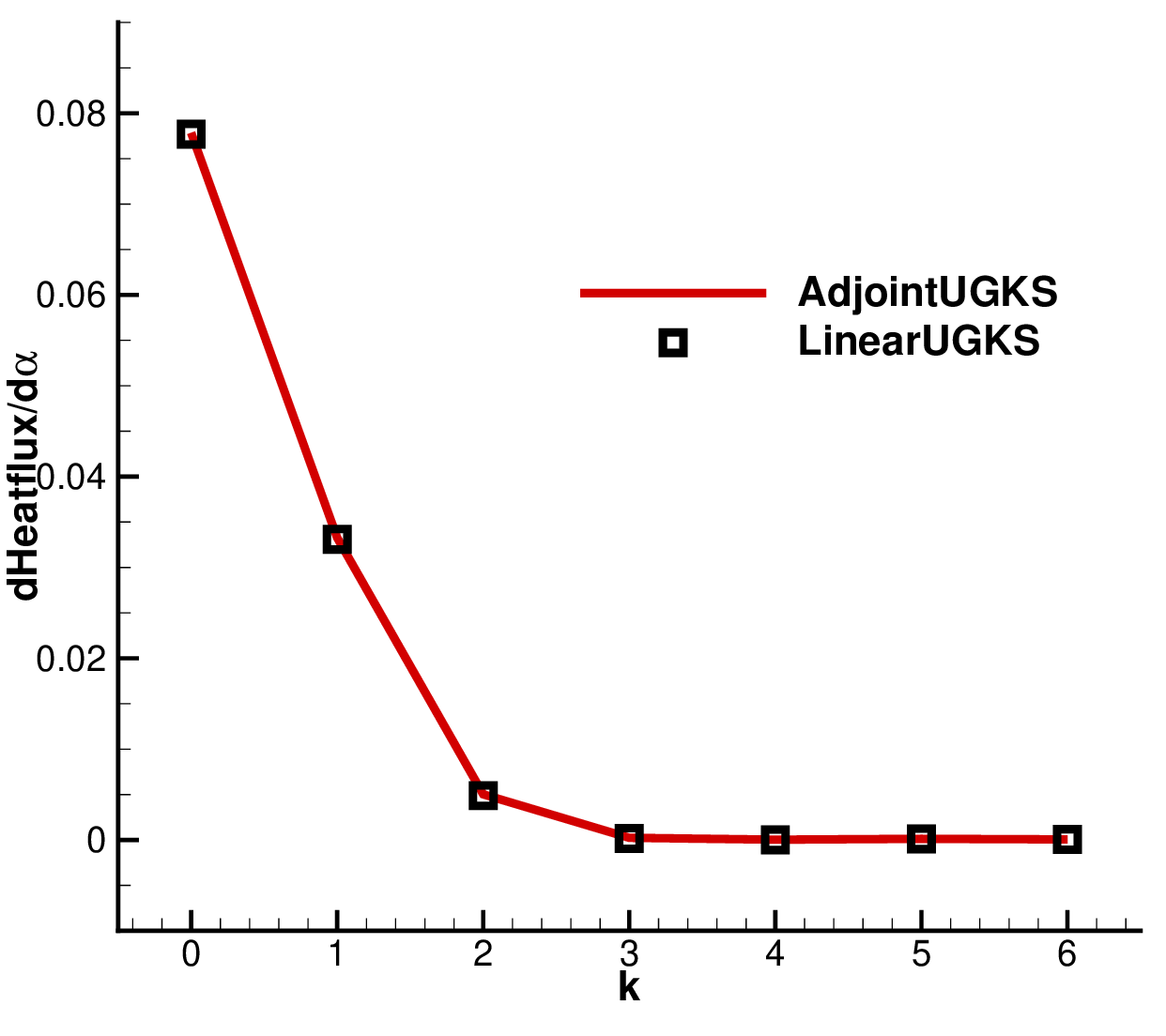}
    \caption{$Kn=0.4$}
    \label{fig:cylinder_sensitive_Kn04}
  \end{subfigure}
  \hfill
  \begin{subfigure}[b]{0.32\textwidth}
    \centering
    \includegraphics[width=\linewidth]{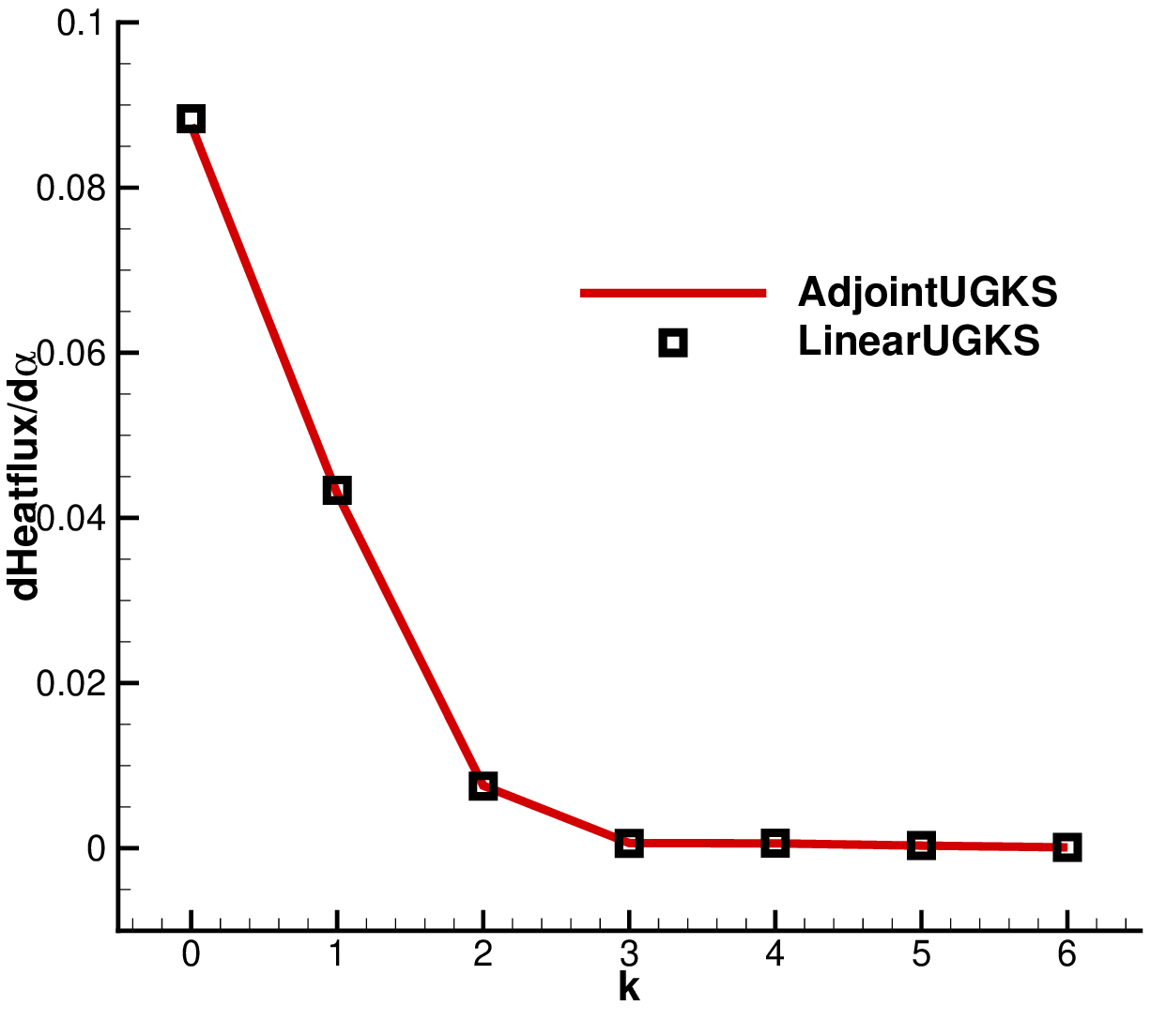}
    \caption{$Kn=0.1$}
    \label{fig:cylinder_sensitive_Kn01}
  \end{subfigure}
  \hfill
  \begin{subfigure}[b]{0.32\textwidth}
    \centering
    \includegraphics[width=\linewidth]{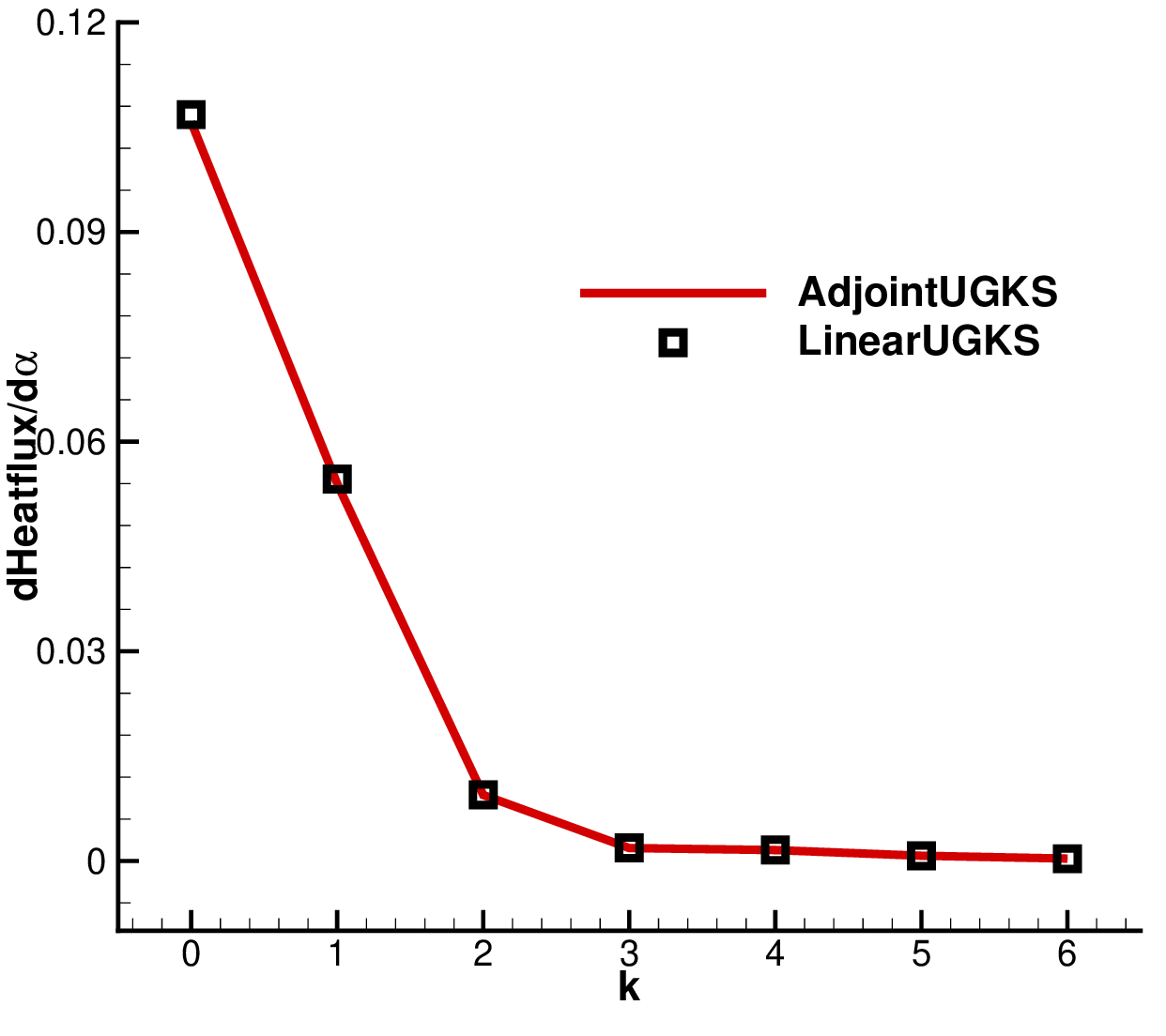}
    \caption{$Kn=0.02$}
    \label{fig:cylinder_sensitive_Kn002}
  \end{subfigure}
  \caption{Modal sensitivities of the cylinder surface heat flux with respect to the perturbation coefficient $a$ in $\alpha=\alpha_0+a\cos(2n\theta)$. Comparison between adjoint UGKS and linearized UGKS.}
  \label{fig:cylinder_sensitive}
\end{figure}

Figure~\ref{fig:cylinder_sensitive} compares the modal sensitivities $\partial J/\partial a$ obtained by the adjoint and linearized methods for $n=0,\ldots,6$.
At all three Knudsen numbers, the two curves coincide, and the sensitivity decays rapidly with increasing mode index: the lowest modes dominate the heat-flux response, whereas higher-order roughness components contribute only weakly.
This modal decay is strongest at $Kn=0.4$ and becomes milder as the flow approaches the continuum limit, consistent with the thicker kinetic boundary layer in the more rarefied regime.

\begin{figure}[!htbp]
  \centering
  \begin{subfigure}[b]{0.45\textwidth}
    \centering
    \includegraphics[width=\linewidth]{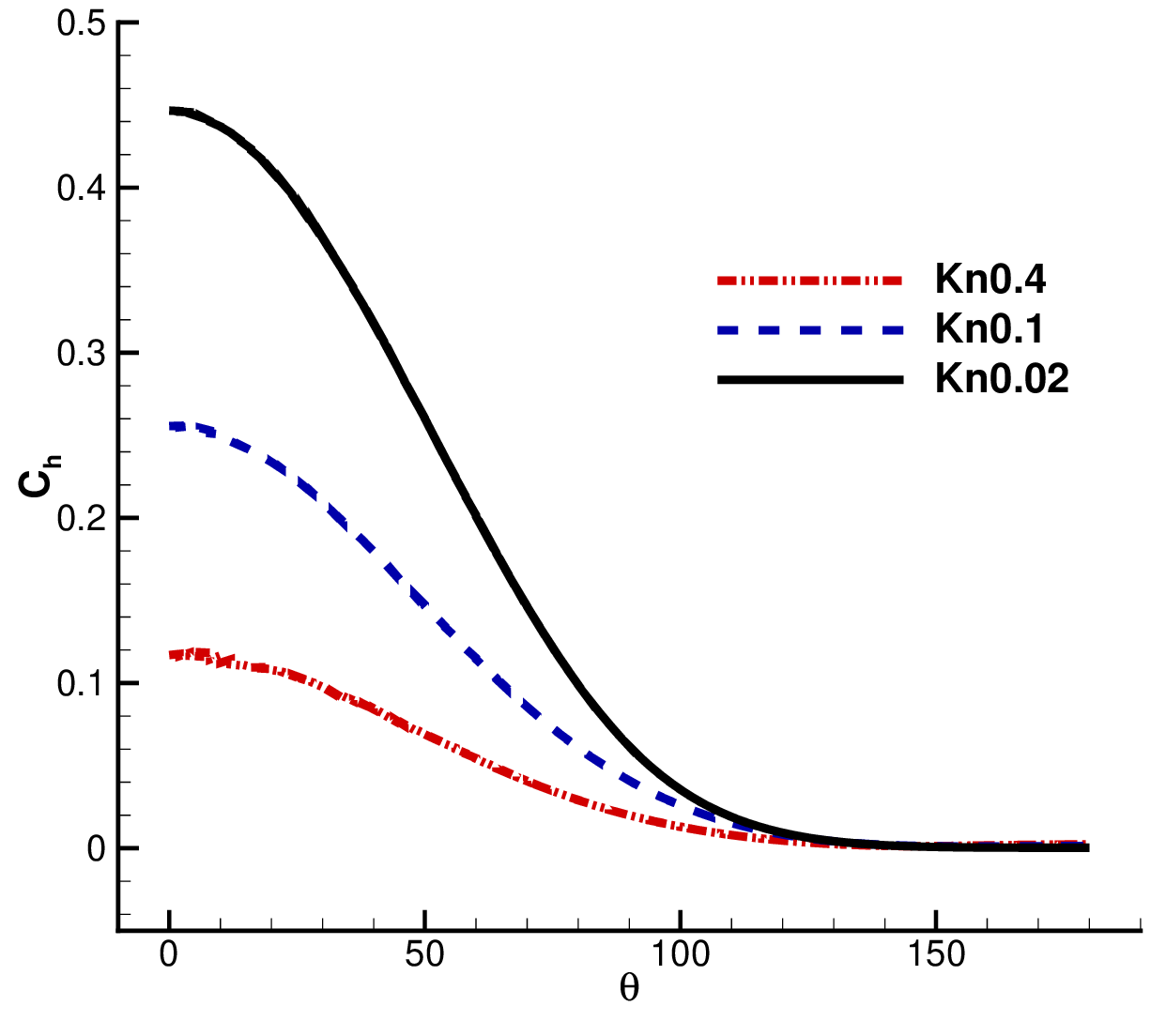}
    \caption{Baseline heat flux}
    \label{fig:cylinder_heatflux}
  \end{subfigure} 
  \quad
  \begin{subfigure}[b]{0.45\textwidth}
    \centering
    \includegraphics[width=\linewidth]{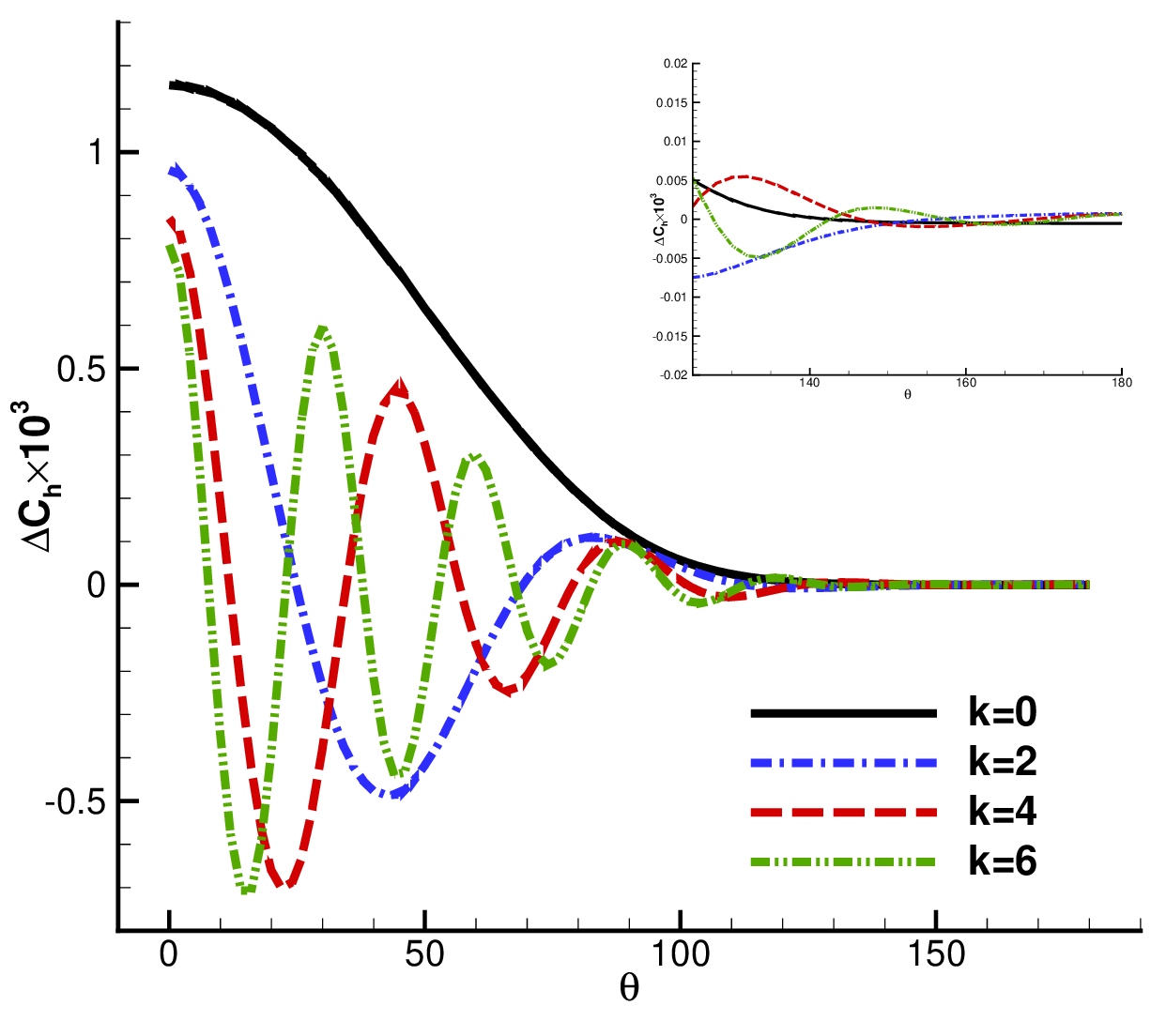}
    \caption{$Kn=0.4$}
    \label{fig:cylinder_linearHF_Kn04}
  \end{subfigure}
  \begin{subfigure}[b]{0.45\textwidth}
    \centering
    \includegraphics[width=\linewidth]{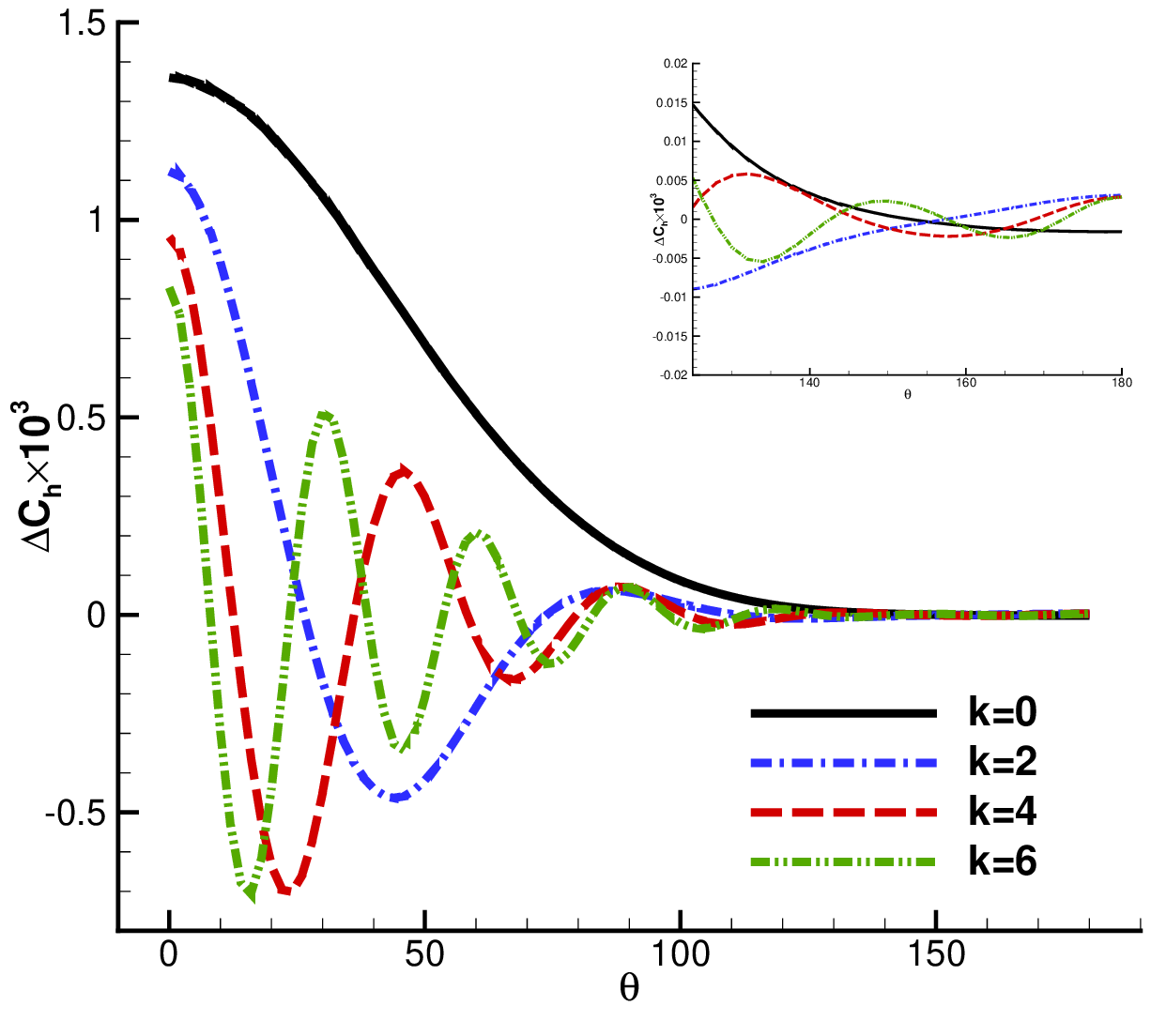}
    \caption{$Kn=0.1$}
    \label{fig:cylinder_linearHF_Kn01}
  \end{subfigure}
  \quad
    \begin{subfigure}[b]{0.45\textwidth}
    \centering
    \includegraphics[width=\linewidth]{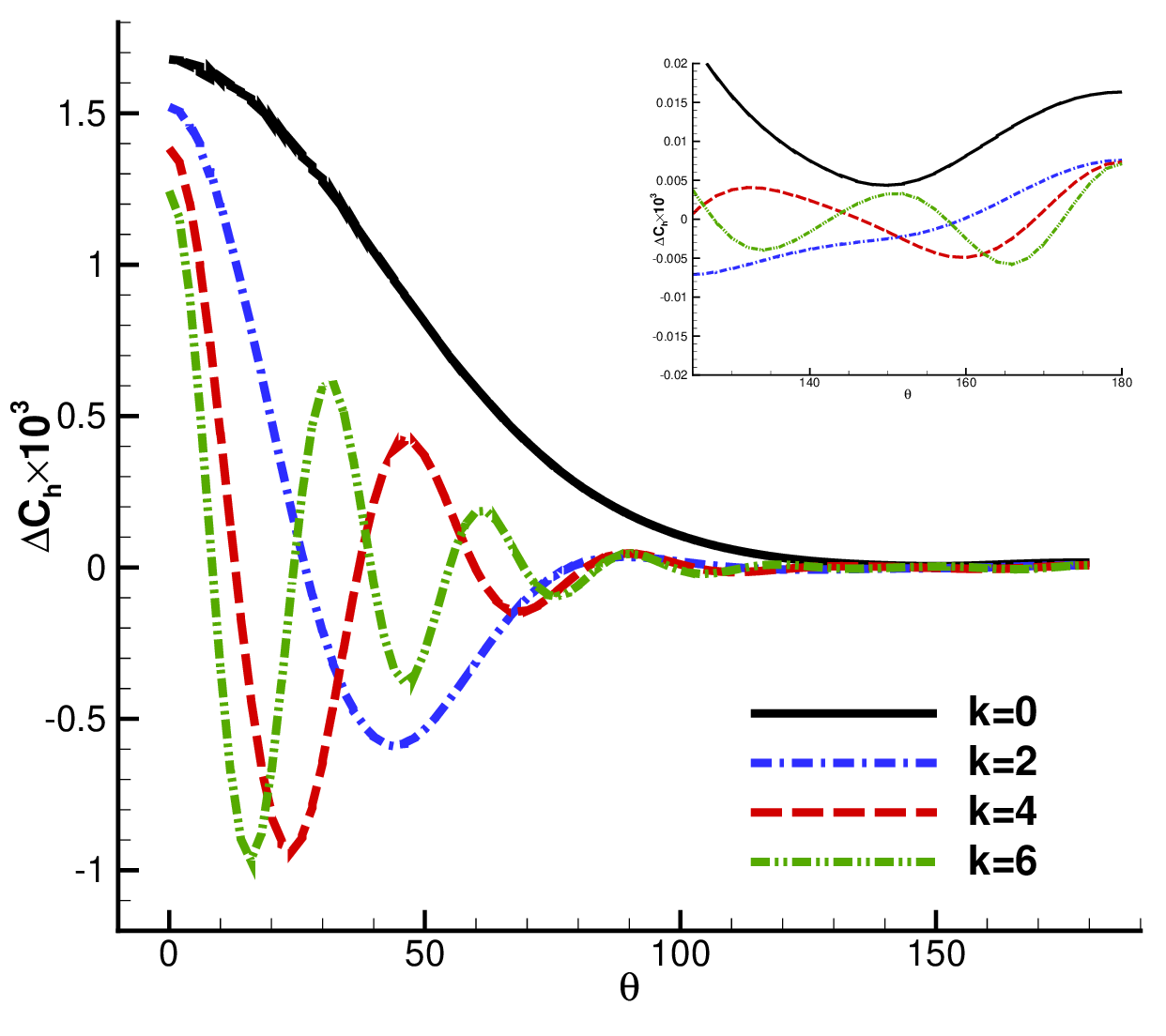}
    \caption{$Kn=0.02$}
    \label{fig:cylinder_linearHF_Kn002}
  \end{subfigure}
  \caption{(a)~Baseline wall heat flux. (b)--(d)~Linearized wall heat-flux perturbations for roughness modes $n=0,2,4,6$ with $\alpha'=0.001\cos(2n\theta)$ at $Kn=0.4$, $0.1$ and $0.02$.}
  \label{fig:cylinder_linear_heatflux}
\end{figure}

\begin{figure}[!htbp]
  \centering
  \includegraphics[width=0.45\textwidth]{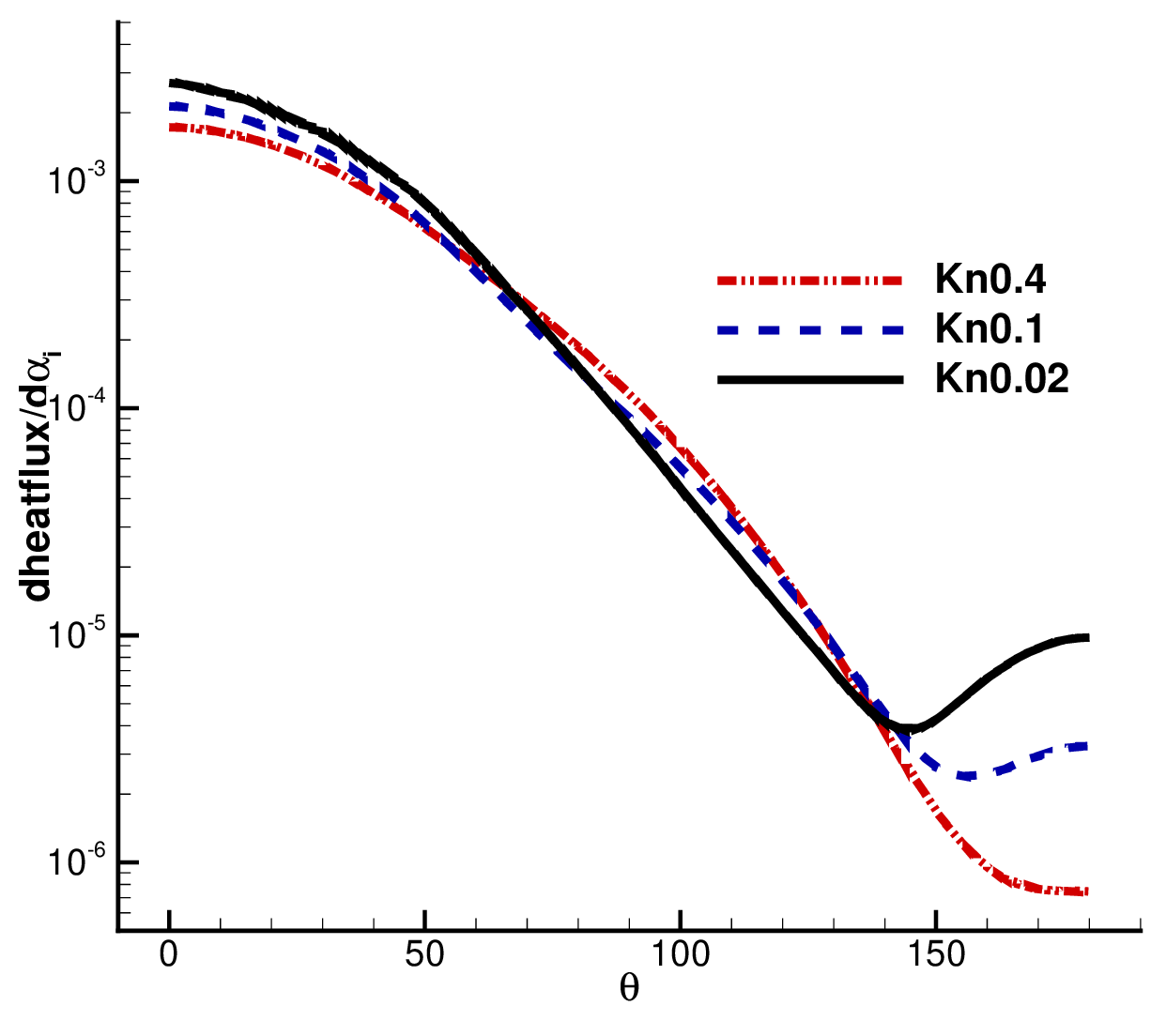}
  \caption{Adjoint sensitivities of the local wall heat flux with respect to a cellwise perturbation of the accommodation coefficient $\alpha$ along the cylinder surface for $Kn=0.4$, $0.1$ and $0.02$.}
  \label{fig:cylinder_adjoint_sensitive}
\end{figure}

The corresponding wall heat-flux distributions are shown in Fig.~\ref{fig:cylinder_linear_heatflux}.
Panel~(a) gives the baseline surface heat flux, while panels~(b)--(d) show the linearized heat-flux perturbations for modes $n=0,2,4,6$ with the fixed disturbance $\alpha'=0.001\cos(2n\theta)$ at $Kn=0.4$, $0.1$ and $0.02$.
Both the baseline heat flux and the linearized heat-flux response at the stagnation point increase as the Knudsen number decreases.
On the lateral surface, the high-$Kn$ cases exhibit larger perturbation amplitudes; as $Kn$ is reduced, these amplitudes diminish, and the disturbance becomes increasingly concentrated near the stagnation region.
On the leeward surface, the heat-flux perturbation grows with decreasing $Kn$.
This trend is associated with the nearly vacuum leeward region at high $Kn$, although the absolute leeward heat flux remains at a very low level in all cases.
Figure~\ref{fig:cylinder_adjoint_sensitive} shows the adjoint-based local sensitivities of the wall heat flux with respect to a cellwise perturbation of the accommodation coefficient $\alpha$ for the three Knudsen numbers.
In this way, the sensitivity on every wall cell is obtained from a single adjoint solve; the modal sensitivities in Fig.~\ref{fig:cylinder_sensitive} can then be recovered by multiplying the cellwise sensitivities by the corresponding mode shape $\cos(2n\theta)$ and summing over the wall surface.
The sensitivity is largest near the stagnation point, where the low-$Kn$ cases attain higher peak values, and decreases with increasing distance from the stagnation region.
On the lateral surface, the high-$Kn$ case yields a larger sensitivity, whereas on the leeward surface the sensitivity increases as $Kn$ decreases.
These trends are consistent with the linearized heat-flux response in Fig.~\ref{fig:cylinder_linear_heatflux}.

\begin{figure}[!htbp]
  \centering
  \begin{subfigure}[b]{0.32\textwidth}
    \centering
    \includegraphics[width=\linewidth]{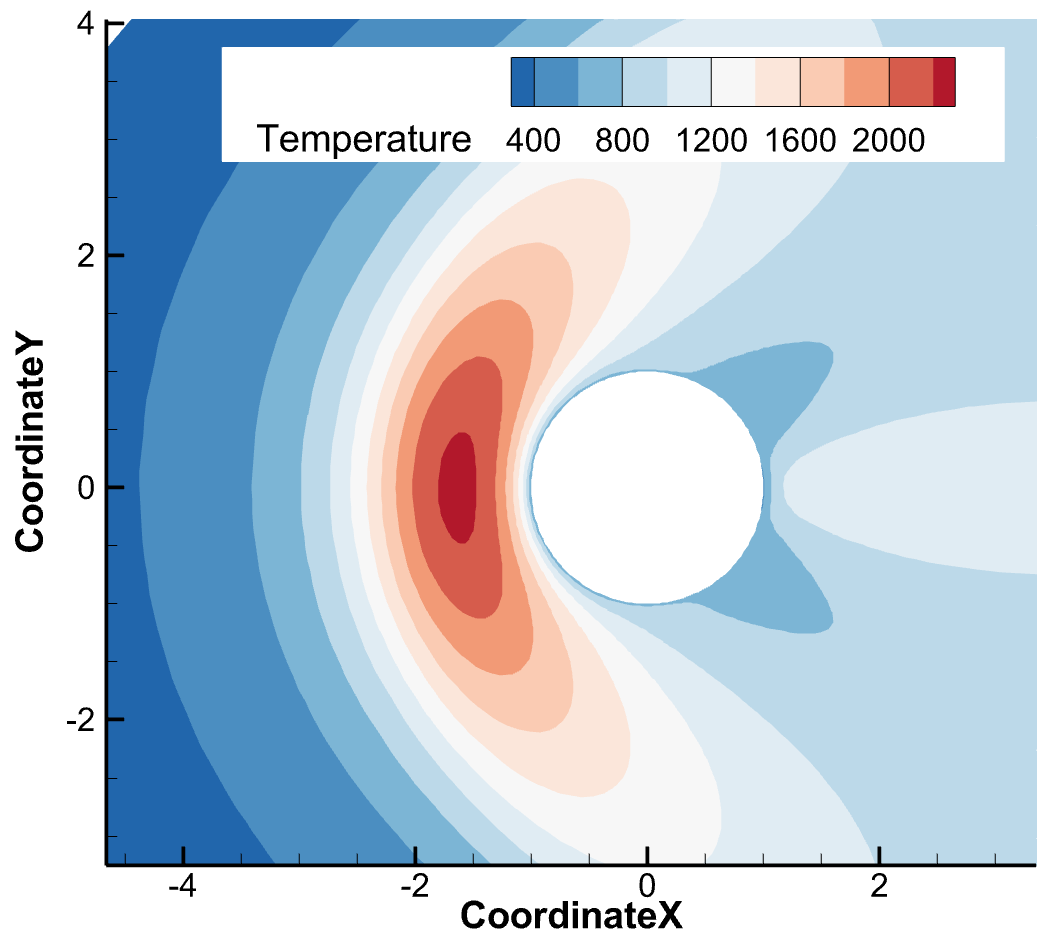}
    \caption{}
    \label{fig:cylinder_Kn04_origin}
  \end{subfigure}
  \hfill
  \begin{subfigure}[b]{0.32\textwidth}
    \centering
    \includegraphics[width=\linewidth]{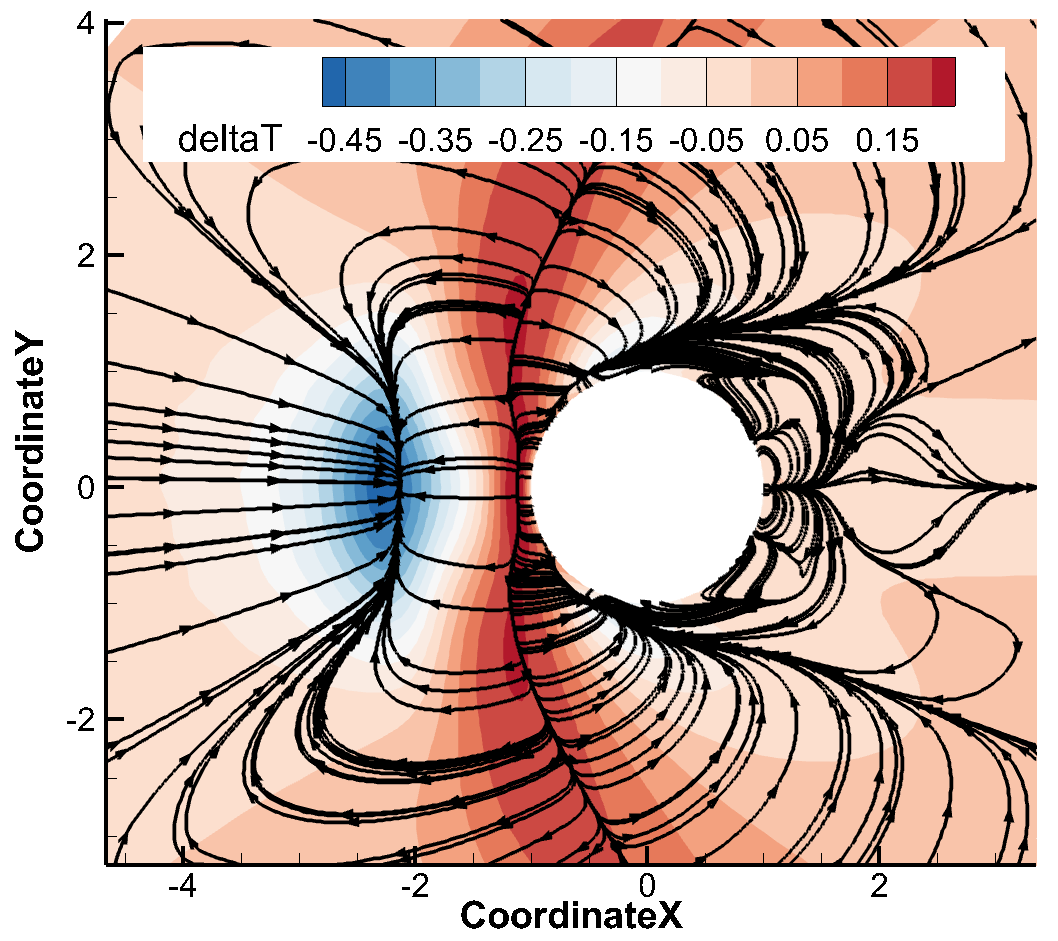}
    \caption{}
    \label{fig:cylinder_Kn04_linear}
  \end{subfigure}
  \hfill
  \begin{subfigure}[b]{0.32\textwidth}
    \centering
    \includegraphics[width=\linewidth]{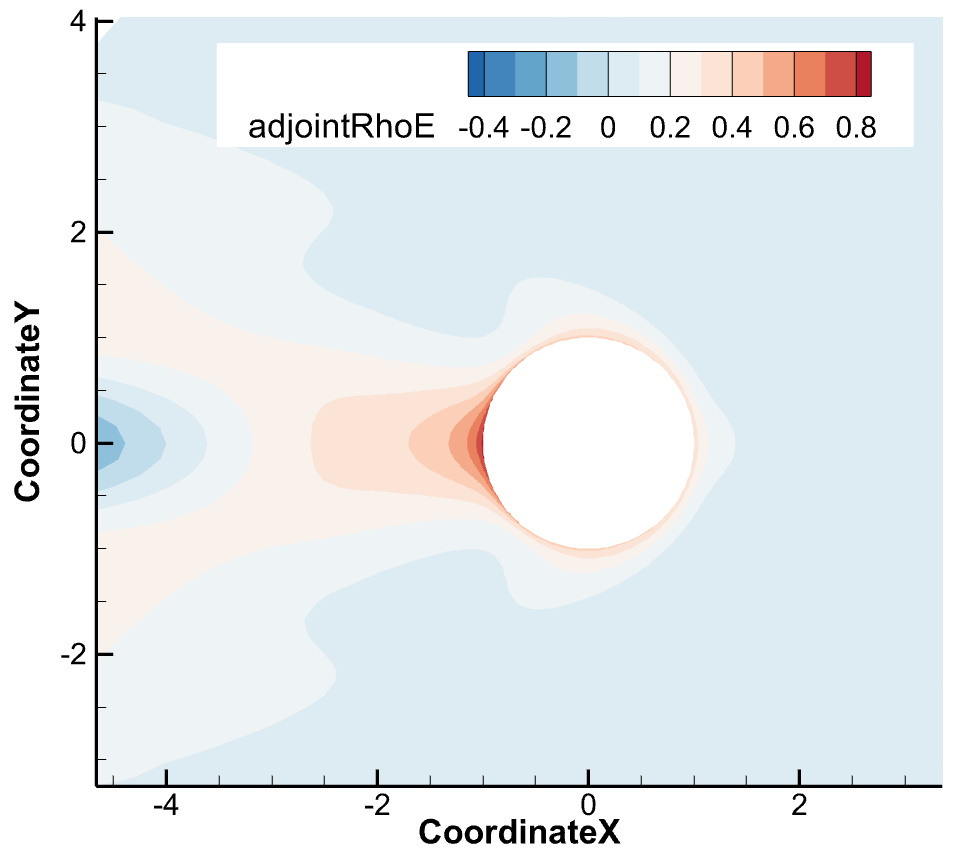}
    \caption{}
    \label{fig:cylinder_Kn04_adjoint}
  \end{subfigure}
  \caption{Field distributions for the cylinder flow at $Kn=0.4$. (a) Baseline temperature. (b) Linearized field (contours: perturbed temperature; streamlines: perturbed heat flux). (c) Adjoint energy field.}
  \label{fig:cylinder_Kn04_contour}
\end{figure}

The field distributions further illustrate the multiscale thermal response.
At $Kn=0.4$ (Fig.~\ref{fig:cylinder_Kn04_contour}), the baseline temperature field exhibits a diffuse bow shock and a broad high-temperature layer ahead of the cylinder.
The linearized temperature perturbation is concentrated on the windward side, and the overlaid streamlines represent the perturbed heat flux.
The adjoint energy field propagates the heat-flux sensitivity information upstream from the wall into the shock layer, with the strongest response at the stagnation surface.

\begin{figure}[!htbp]
  \centering
  \begin{subfigure}[b]{0.32\textwidth}
    \centering
    \includegraphics[width=\linewidth]{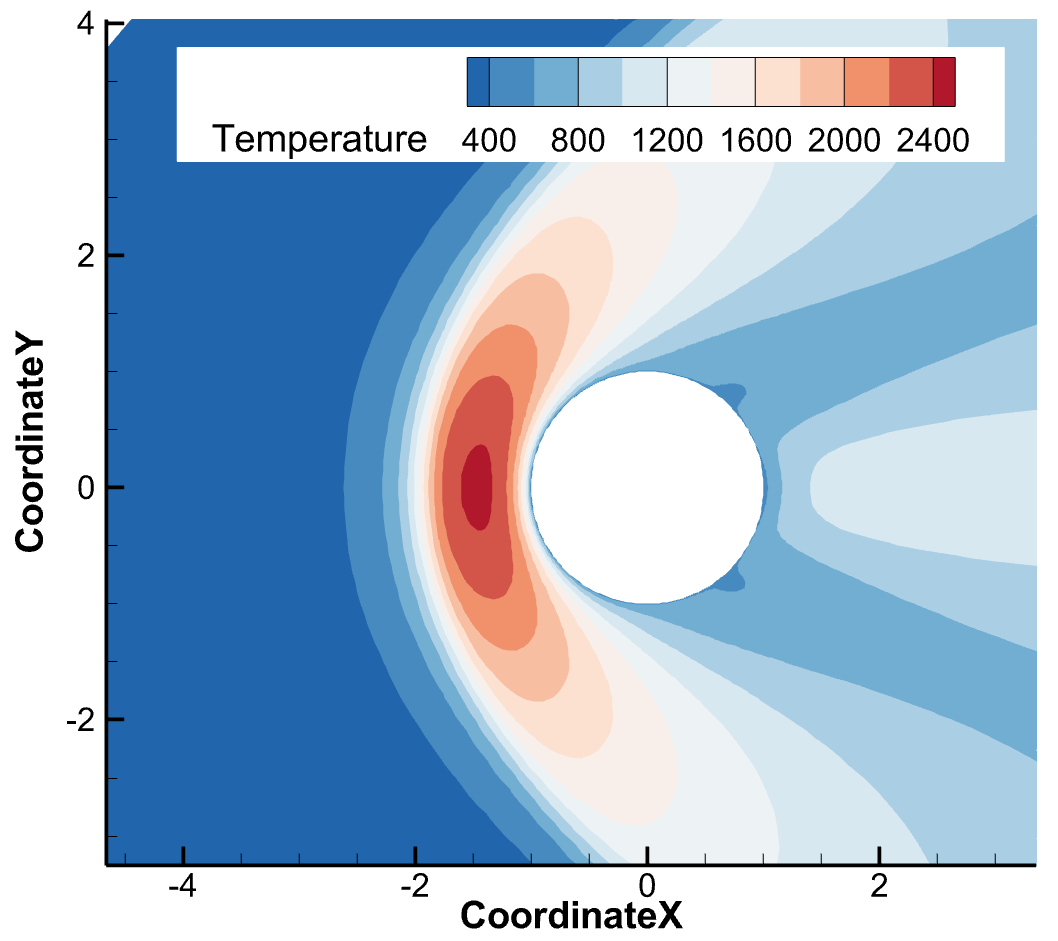}
    \caption{}
    \label{fig:cylinder_Kn01_origin}
  \end{subfigure}
  \hfill
  \begin{subfigure}[b]{0.32\textwidth}
    \centering
    \includegraphics[width=\linewidth]{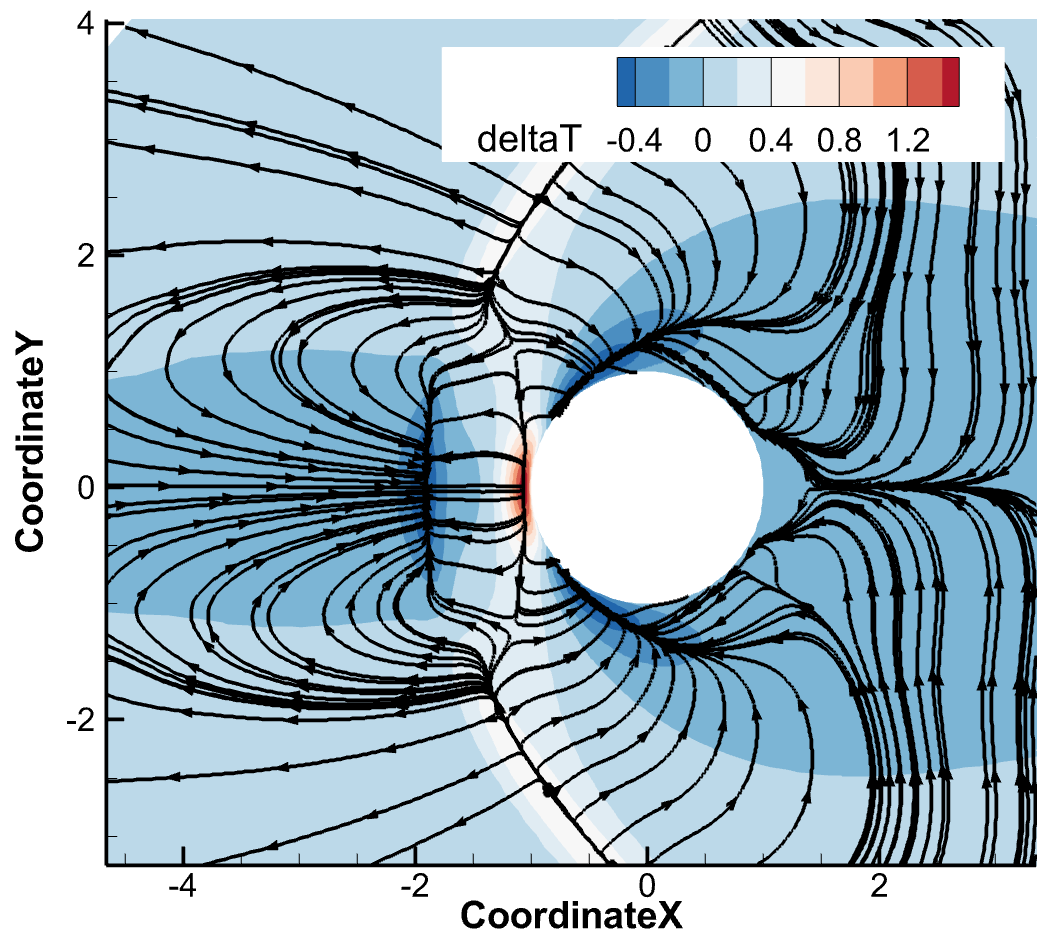}
    \caption{}
    \label{fig:cylinder_Kn01_linear}
  \end{subfigure}
  \hfill
  \begin{subfigure}[b]{0.32\textwidth}
    \centering
    \includegraphics[width=\linewidth]{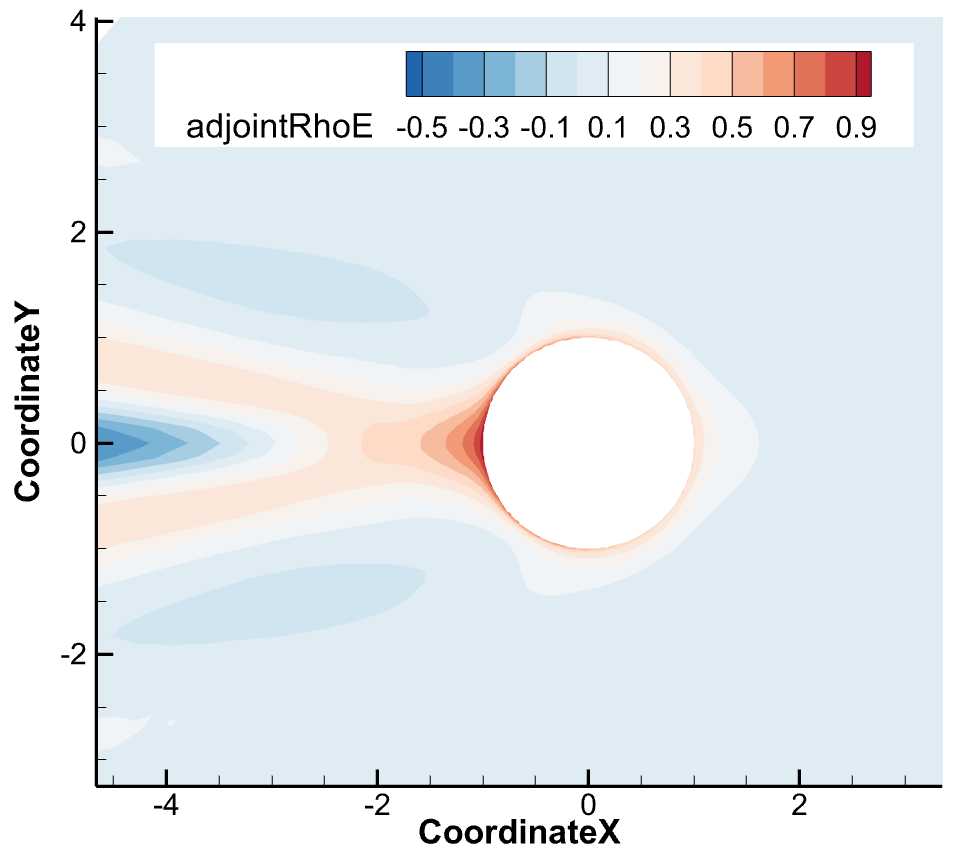}
    \caption{}
    \label{fig:cylinder_Kn01_adjoint}
  \end{subfigure}
  \caption{Field distributions for the cylinder flow at $Kn=0.1$. (a) Baseline temperature. (b) Linearized field (contours: perturbed temperature; streamlines: perturbed heat flux). (c) Adjoint energy field.}
  \label{fig:cylinder_Kn01_contour}
\end{figure}

As the Knudsen number is reduced to $Kn=0.1$ (Fig.~\ref{fig:cylinder_Kn01_contour}), the shock stand-off decreases, and the thermal layer tightens against the windward surface.
The linearized and adjoint fields remain focused on the stagnation region, while the wake contribution becomes relatively weaker than in the $Kn=0.4$ case.
The perturbed heat-flux streamlines in the linearized field likewise concentrate on the windward surface.

\begin{figure}[!htbp]
  \centering
  \begin{subfigure}[b]{0.32\textwidth}
    \centering
    \includegraphics[width=\linewidth]{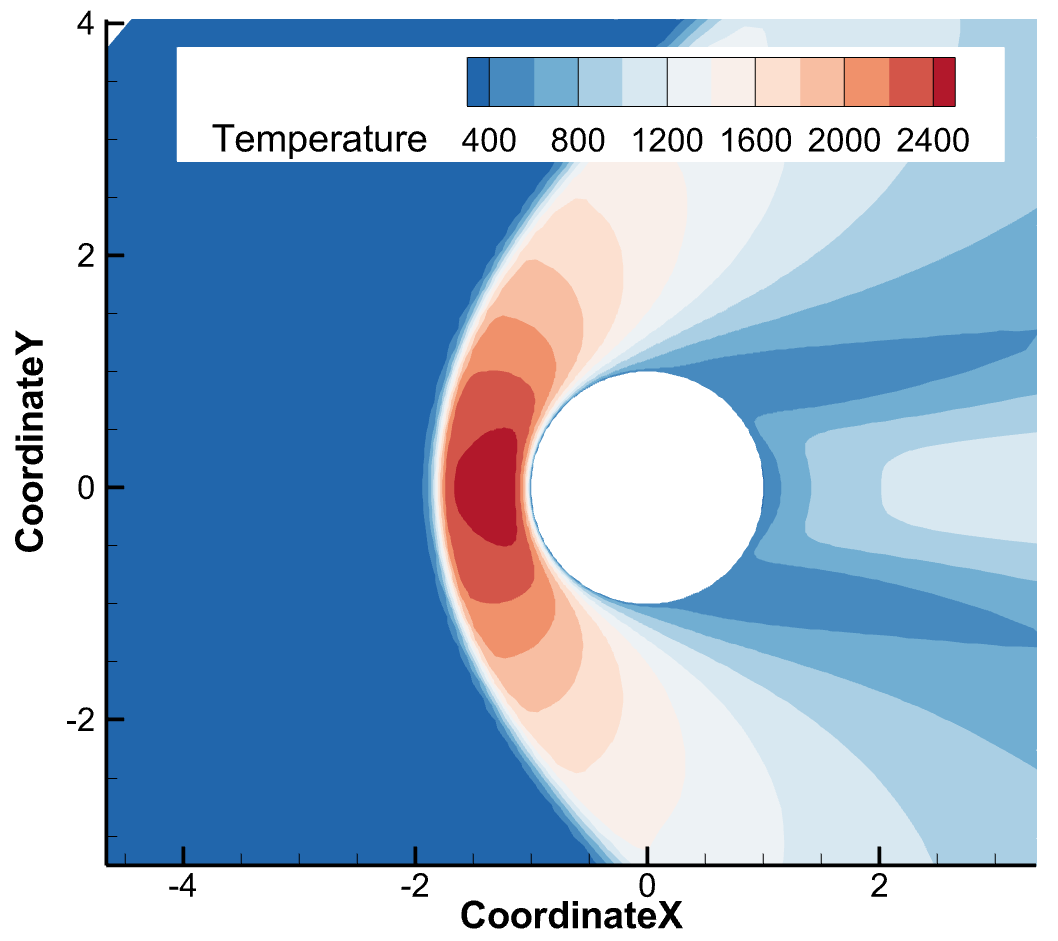}
    \caption{}
    \label{fig:cylinder_Kn002_origin}
  \end{subfigure}
  \hfill
  \begin{subfigure}[b]{0.32\textwidth}
    \centering
    \includegraphics[width=\linewidth]{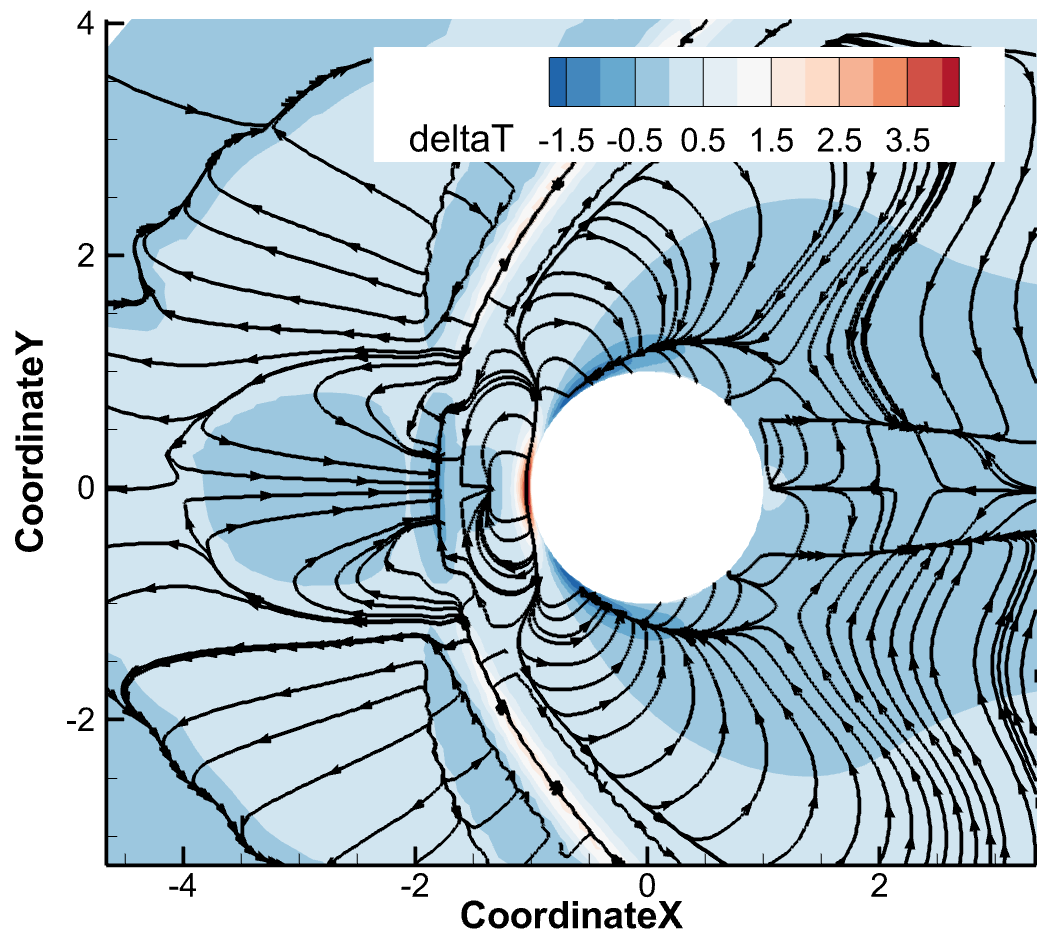}
    \caption{}
    \label{fig:cylinder_Kn002_linear}
  \end{subfigure}
  \hfill
  \begin{subfigure}[b]{0.32\textwidth}
    \centering
    \includegraphics[width=\linewidth]{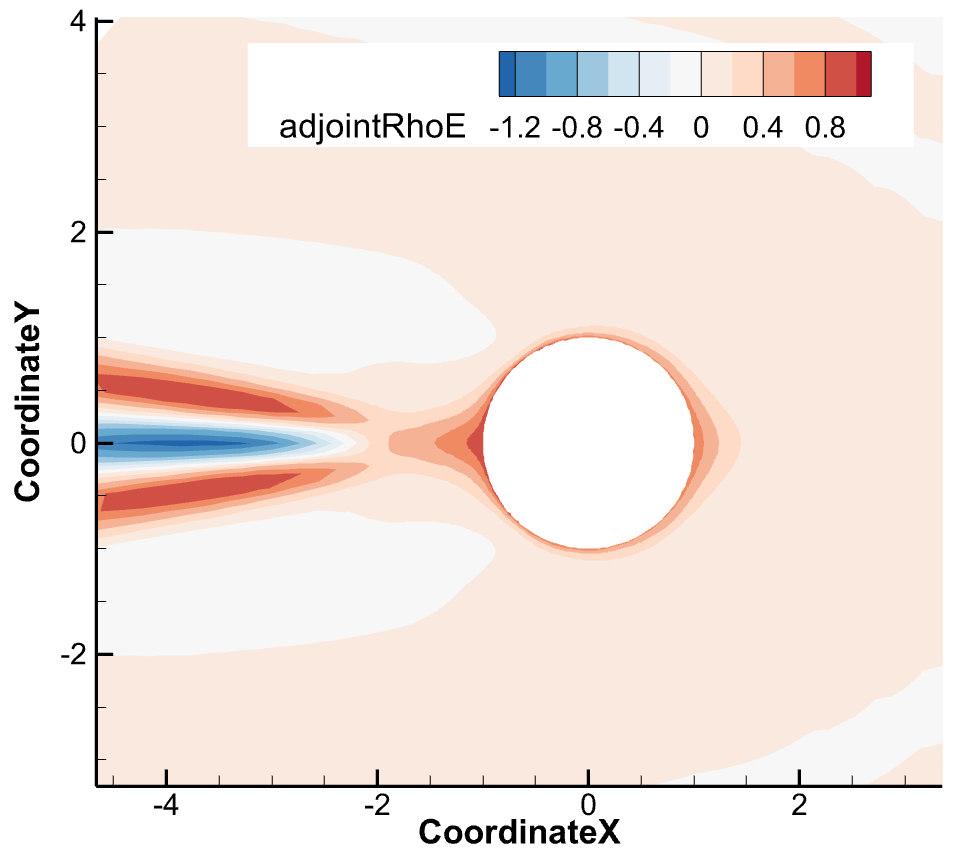}
    \caption{}
    \label{fig:cylinder_Kn002_adjoint}
  \end{subfigure}
  \caption{Field distributions for the cylinder flow at $Kn=0.02$. (a) Baseline temperature. (b) Linearized field (contours: perturbed temperature; streamlines: perturbed heat flux). (c) Adjoint energy field.}
  \label{fig:cylinder_Kn002_contour}
\end{figure}

In the near-continuum regime $Kn=0.02$ (Fig.~\ref{fig:cylinder_Kn002_contour}), a sharp bow shock and a thin high-temperature sheath form ahead of the cylinder.
The linearized perturbation is tightly attached to the windward surface, with the perturbed heat-flux streamlines forming a slender upstream structure, and the adjoint energy field develops a corresponding upstream lobe.
These trends confirm that roughness-induced heat-flux sensitivities are controlled by the lowest geometric modes and by the near-wall kinetic layer, and that the adjoint UGKS recovers the same sensitivities as the linearized solver across the transitional and continuum limits.

\section{Conclusion}
A discrete adjoint method for the UGKS has been developed for multiscale gas dynamics.
Because the macroscopic equation can be recovered by taking moments of the microscopic equation together with the collision compatibility condition, the adjoint system is derived from the discrete microscopic equation in conjunction with the compatibility relation.
To remove the stiff feedback in the coupled adjoint system, an intermediate macroscopic adjoint moment is introduced, and asymptotic-preserving (AP) adjoint equations are constructed through two algebraically equivalent routes: macroscopic-moment projection and microscopic lifting.
The former eliminates the stiff collision terms and yields the evolution equation for the macroscopic adjoint moment, while the latter produces a microscopic adjoint equation that is explicitly coupled to the macroscopic moment; the AP adjoint system is composed of these two equations.
In the continuum limit, the system reduces to the Euler adjoint equations, whereas in the free-molecular limit it recovers the discrete-velocity-method adjoint equations.
In addition, the adjoint form of the Maxwell wall reflection model with an accommodation coefficient is derived.

Numerically, a memory-efficient residual-evaluation architecture dual to the forward UGKS is adopted: the common $\mathbf{W}$ and $\mathbf{q}$ contributions are precomputed, and the pointwise residuals are then assembled in the velocity-space loop, with the adjoint residual assembled in the reverse order of the forward solver.
Under this framework, the microscopic adjoint solve requires only the storage of the microscopic distribution functions, the adjoint variables, and their residuals; all remaining quantities are independent of the velocity-space size, so that the memory cost remains manageable.
Time advancement is performed by a macroscopic--microscopic predictor--corrector implicit scheme, thereby avoiding a global matrix inversion over the full velocity space.

Using the linearized UGKS as an independent reference, the method has been systematically verified for lid-driven cavity heat conduction, closed-channel thermal creep flow, and hypersonic flow past a circular cylinder.
From continuum to transitional regimes, the adjoint sensitivities agree closely with the linearized results, and the residual convergence histories are likewise consistent.
The cylinder case further shows that a single adjoint solve yields the sensitivities of the wall heat flux with respect to the accommodation coefficient on every wall cell, from which the modal roughness sensitivities can be reconstructed.
In the strongly rarefied regime (e.g., $Kn=1$), discrepancies between the adjoint and linearized sensitivities may appear, and the adjoint residual may fail to converge to a satisfactorily low level; this behavior is associated with the incomplete duality between the adjoint and linearized formulations and with ray effects in the discrete velocity space, and warrants further investigation.

Overall, the present method provides a reliable foundation for sensitivity analysis, uncertainty quantification, and design optimization of multiscale gas flows.
Future work will address three-dimensional complex geometries, shape and topology optimization, and coupling with adaptive velocity-space discretization.

 \section*{Acknowledgements}
This work was supported by the National Natural Science Foundation of China (Grant No. 92371107), and the Hong Kong Research Grants Council (Grant No. 16208324).


\bibliographystyle{elsarticle-num}
\bibliography{sample}

\appendix

\section{Derivation of the adjoint system from the primal macroscopic and microscopic governing equations} \label{app:ill_conditioned}
This appendix complements the compatibility-based adjoint derivation in the main text by working directly with the explicit macroscopic update and the implicit microscopic update used in the forward UGKS.
The primal macroscopic operator is the conservative flux-residual form in Eq.~\eqref{macroscopicUpdata}:
\begin{equation*}
  \mathcal{L}_{\text{macro}}( \mathbf{W}_i, f_{i,k}) = \frac{\mathbf{W}_i^{n+1}-\mathbf{W}_i^{n}}{\Delta t}+\frac{1}{|\Omega_i|}\sum_k w_k \mathcal{R}_{i,k}^n \mathbf{\Psi}_k = 0,
\end{equation*}
and the corresponding microscopic operator is the implicit collision-residual form in Eq.~\eqref{UGKSmicro}:
\begin{equation*}
  \mathcal{L}_{\text{micro}}(\mathbf{W}_i, f_{i,k}) = \frac{f_{i,k}^{n+1}-f_{i,k}^{n}}{\Delta t}+\frac{1}{|\Omega_i|}\mathcal{R}_{i,k}^n - \frac{f_{i,k}^{+,n+1} - f_{i,k}^{n+1}}{\tau_i} = 0.
\end{equation*}
To evaluate the sensitivity of the objective function $J$, the Lagrangian is formed with the macroscopic and microscopic inner products defined in the main text.
Taking first-order variations with respect to $\mathbf{W}_i$ and $f_{i,k}$, and normalizing by $|\Omega_i|$ and $w_k$, yields the inhomogeneous adjoint system
\begin{equation}
  \label{eq:app_inhomogeneous}
\begin{cases}
\mathcal{L}_{\text{macro}}^*(\mathbf{\Lambda}_i^{*,n},\lambda_{i,k}^{*,n}) + \dfrac{\partial J}{\partial \mathbf{W}_i^n}/|\Omega_i| = 0, \\
\mathcal{L}_{\text{micro}}^*(\mathbf{\Lambda}_i^{*,n},\lambda_{i,k}^{*,n}) + \dfrac{\partial J}{\partial f_{i,k}^n}/(w_k|\Omega_i|) = 0,
\end{cases}
\end{equation}
where the microscopic adjoint operator is
\begin{equation}
  \label{eq:app_scaled_micro}
  \begin{aligned}
  \mathcal{L}_{\text{micro}}^*=\frac{\lambda_{i,k}^{*,n-1}-\lambda_{i,k}^{*,n}}{\Delta t}&+\frac{1}{|\Omega_i|}\sum_j\frac{\partial \mathcal{R}_{j,k}}{\partial f_{i,k}}(\mathbf{\Psi}_k^T \mathbf{\Lambda}_j^{*,n}+\lambda_{j,k}^{*,n})\\
  &+ \frac{1}{|\Omega_i|} \frac{1}{w_k}\frac{\partial \mathbf{q}_i}{\partial f_{i,k}}\sum_j \sum_{k^\prime}w_{k^\prime}\frac{\partial \mathcal{R}_{j,k^\prime}}{\partial \mathbf{q}_i}(\mathbf{\Psi}_k^T \mathbf{\Lambda}_j^{*,n}+\lambda_{j,k^\prime}^{*,n}) \\
  &+\frac{\lambda_{i,k}^{*,n-1}}{\tau_i}-\frac{1}{\tau_i}\frac{1}{w_k}\frac{\partial \mathbf{q}_i}{\partial f_{i,k}}\sum_{k^\prime} w_{k^\prime}  \frac{\partial f_{i,k^\prime}^+}{\partial \mathbf{q}_i}  \lambda_{i,k^\prime}^{*,n}= 0,
  \end{aligned}
\end{equation}
and the macroscopic adjoint operator is
\begin{equation}
  \label{eq:app_macro_adjoint}
  \begin{aligned}
  \mathcal{L}_{\text{macro}}^*=\frac{\mathbf{\Lambda}_i^{*,n-1}-\mathbf{\Lambda}_i^{*,n}}{\Delta t}&+\frac{1}{|\Omega_i|}\sum_j\sum_{k} w_k\frac{\partial \mathcal{R}_{j,k}}{\partial \mathbf{W}_i}(\mathbf{\Psi}_k^T \mathbf{\Lambda}_j^{*,n}+\lambda_{j,k}^{*,n})\\
  &+ \frac{1}{|\Omega_i|}\frac{\partial \mathbf{q}_i}{\partial \mathbf{W}_i}\sum_j\sum_{k} w_k\frac{\partial \mathcal{R}_{j,k}}{\partial \mathbf{q}_{i}}(\mathbf{\Psi}_k^T \mathbf{\Lambda}_j^{*,n}+\lambda_{j,k}^{*,n})\\
  &-\frac{1}{\tau_i}\sum_{k} w_k  \frac{\partial f_{i,k}^+}{\partial \mathbf{W}_i}\lambda_{i,k}^{*,n-1}-\frac{1}{\tau_i}\frac{\partial \mathbf{q}_i}{\partial \mathbf{W}_i}\sum_{k} w_k  \frac{\partial f_{i,k}^+}{\partial \mathbf{q}_i} \lambda_{i,k}^{*,n}\\
  &+\frac{1}{\tau_i}\frac{\partial \tau_i}{\partial \mathbf{W}_i}\sum_{k} w_k\frac{f_{i,k}^+-f_{i,k}}{\tau_i}\lambda_{i,k}^{*,n}= 0.
  \end{aligned}
\end{equation}
In Eq.~\eqref{eq:app_scaled_micro} and Eq.~\eqref{eq:app_macro_adjoint}, the flux and heat-flux Jacobians are evaluated through the combined microscopic-macroscopic adjoint state $\mathbf{\Psi}_k^T \mathbf{\Lambda}_j^{*,n}+\lambda_{j,k}^{*,n}$.
In Eqs.~\eqref{eq:app_inhomogeneous}--\eqref{eq:app_macro_adjoint}, the macroscopic and microscopic adjoint equations remain formally decoupled: the microscopic block must be solved first, and the macroscopic block is then updated sequentially.
This one-way structure differs from the coupled multiscale formulation adopted in the main text, where explicit cross-scale coupling appears already in Eqs.~\eqref{scaled_micro} and~\eqref{adjointmacro_final}.
Nevertheless, the same algebraic AP reconstruction can be applied directly to Eqs.~\eqref{eq:app_macro_adjoint} and~\eqref{eq:app_scaled_micro}, and yields an AP adjoint system equivalent to that derived in the main text.

Two equivalent reconstruction routes are constructed from Eqs.~\eqref{eq:app_macro_adjoint} and~\eqref{eq:app_scaled_micro}.
Throughout both routes, the combined microscopic adjoint variable
\begin{equation*}
  \widetilde{\lambda}_{i,k}^{*,n}=\lambda_{i,k}^{*,n}+ \mathbf{\Psi}_k^T \mathbf{\Lambda}_i^{*,n},
\end{equation*}
is used so that $\mathbf{\Psi}_k^T \mathbf{\Lambda}_j^{*,n}+\lambda_{j,k}^{*,n}=\widetilde{\lambda}_{j,k}^{*,n}$.
Because the present route retains the explicit macroscopic time derivative in Eq.~\eqref{eq:app_macro_adjoint}, the intermediate macroscopic adjoint moment is defined as
\begin{equation*}
  \widetilde{\mathbf{\Lambda}}_i^{*,n}=\sum_{k} w_k \frac{\partial f_{i,k}^+}{\partial \mathbf{W}_i}\lambda_{i,k}^{*,n}+\mathbf{\Lambda}_i^{*,n},
\end{equation*}
which reduces to the definition used in the main text once the steady-state backward-marching residual is absorbed.

\textbf{Macroscopic-moment projection route}

Multiplying Eq.~\eqref{eq:app_scaled_micro} by $w_k\partial f_{i,k}^+/\partial \mathbf{W}_i$ and summing over $k$ gives
\begin{equation}
  \label{eq:app_AP_projection}
  \begin{aligned}
  \sum_{k} w_k \frac{\partial f_{i,k}^+}{\partial \mathbf{W}_i}\frac{\lambda_{i,k}^{*,n-1}-\lambda_{i,k}^{*,n}}{\Delta t} &+ \frac{1}{|\Omega_i|}\sum_{k} w_k\frac{\partial f_{i,k}^+}{\partial \mathbf{W}_i}\sum_j \frac{\partial \mathcal{R}_{j,k}}{\partial f_{i,k}}(\mathbf{\Psi}_k^T \mathbf{\Lambda}_j^{*,n}+\lambda_{j,k}^{*,n}) \\
  &+ \frac{1}{|\Omega_i|}\sum_{k} \frac{\partial f_{i,k}^+}{\partial \mathbf{W}_i}\frac{\partial \mathbf{q}_i}{\partial f_{i,k}}\sum_j \sum_{k^\prime}w_{k^\prime}\frac{\partial \mathcal{R}_{j,k^\prime}}{\partial \mathbf{q}_i}(\mathbf{\Psi}_k^T \mathbf{\Lambda}_j^{*,n}+\lambda_{j,k^\prime}^{*,n}) \\
  &+ \frac{1}{\tau_i}\sum_{k} w_k \lambda_{i,k}^{*,n-1} \frac{\partial f_{i,k}^+}{\partial \mathbf{W}_i} - \frac{1}{\tau_i}\sum_{k} \frac{\partial f_{i,k}^+}{\partial \mathbf{W}_i}\frac{\partial \mathbf{q}_i}{\partial f_{i,k}}\sum_{k^\prime} w_{k^\prime}  \frac{\partial f_{i,k^\prime}^+}{\partial \mathbf{q}_i} \lambda_{i,k^\prime}^{*,n} =0.
  \end{aligned}
\end{equation}
Adding Eq.~\eqref{eq:app_macro_adjoint} to Eq.~\eqref{eq:app_AP_projection}, the terms proportional to $\partial f_{i,k}^+/\partial \mathbf{W}_i$ and $\lambda_{i,k}^{*,n-1}$ cancel exactly.
As in the main text, the remaining $\mathbf{q}$-dependent blocks are multiplied by the total heat-flux derivative $\mathrm{d}\mathbf{q}_i/\mathrm{d}\mathbf{W}_i|_{f=f^+}$, which vanishes identically, and therefore drop out, yielding
\begin{equation}
  \label{eq:app_AP_system}
  \begin{aligned}
  \sum_{k} w_k  \frac{\partial f_{i,k}^+}{\partial \mathbf{W}_i}\frac{\lambda_{i,k}^{*,n-1}-\lambda_{i,k}^{*,n}}{\Delta t} &+ \frac{\mathbf{\Lambda}_i^{*,n-1}-\mathbf{\Lambda}_i^{*,n}}{\Delta t}\\
  &+ \frac{1}{|\Omega_i|}\sum_{k} w_k\frac{\partial f_{i,k}^+}{\partial \mathbf{W}_i}\sum_j \frac{\partial \mathcal{R}_{j,k}}{\partial f_{i,k}}(\mathbf{\Psi}_k^T \mathbf{\Lambda}_j^{*,n}+\lambda_{j,k}^{*,n})\\
  &+\frac{1}{|\Omega_i|}\sum_j\sum_{k} w_k \frac{\partial \mathcal{R}_{j,k}}{\partial \mathbf{W}_i}(\mathbf{\Psi}_k^T \mathbf{\Lambda}_j^{*,n}+\lambda_{j,k}^{*,n})\\
  &+\frac{1}{\tau_i}\frac{\partial \tau_i}{\partial \mathbf{W}_i}\sum_{k} w_k\frac{f_{i,k}^+-f_{i,k}}{\tau_i}\lambda_{i,k}^{*,n} =0.
  \end{aligned}
\end{equation}
Substituting $\widetilde{\lambda}_{j,k}^{*,n}$ into the flux blocks and replacing $\sum_{k} w_k (\partial f_{i,k}^+/\partial \mathbf{W}_i)\lambda_{i,k}^{*,n}$ by $\widetilde{\mathbf{\Lambda}}_i^{*,n}$, Eq.~\eqref{eq:app_AP_system} becomes
\begin{equation}
  \label{eq:app_AP_tilde}
  \begin{aligned}
  \frac{\widetilde{\mathbf{\Lambda}}_i^{*,n-1}-\widetilde{\mathbf{\Lambda}}_i^{*,n}}{\Delta t} 
  &+ \frac{1}{|\Omega_i|}\sum_{k} w_k\frac{\partial f_{i,k}^+}{\partial \mathbf{W}_i}\sum_j \frac{\partial \mathcal{R}_{j,k}}{\partial f_{i,k}}\widetilde{\lambda}_{j,k}^{*,n}+\frac{1}{|\Omega_i|}\sum_j\sum_{k} w_k \frac{\partial \mathcal{R}_{j,k}}{\partial \mathbf{W}_i}\widetilde{\lambda}_{j,k}^{*,n}\\
  &+\frac{1}{\tau_i}\frac{\partial \tau_i}{\partial \mathbf{W}_i}\sum_{k} w_k\frac{f_{i,k}^+-f_{i,k}}{\tau_i}\widetilde{\lambda}_{i,k}^{*,n} =0,
  \end{aligned}
\end{equation}
because the macroscopic moment projection of the collision residual vanishes identically,
\begin{equation*}
  \sum_{k^\prime} w_{k^\prime}\frac{f_{i,k^\prime}^+-f_{i,k^\prime}}{\tau_i}\mathbf{\Psi}_{k^\prime}^T\mathbf{\Lambda}_i^{*,n}=\left(\mathbf{\Lambda}_i^{*,n}\right)^T\sum_{k^\prime} w_{k^\prime}\frac{f_{i,k^\prime}^+-f_{i,k^\prime}}{\tau_i}\mathbf{\Psi}_{k^\prime}=0.
\end{equation*}
Under the steady-state backward marching adopted in the main text, the macroscopic marching residual $(\mathbf{\Lambda}_i^{*,n-1}-\mathbf{\Lambda}_i^{*,n})/\Delta t$ is absorbed into $\widetilde{\mathbf{\Lambda}}_i^{*,n}$, and the remaining $\mathbf{\Lambda}_i^{*,n}$-difference terms cancel when $\mathbf{\Psi}_k^T \mathbf{\Lambda}_j^{*,n}$ is separated from $\widetilde{\lambda}_{j,k}^{*,n}$ in the flux blocks.
Eq.~\eqref{eq:app_AP_tilde} then reduces to Eq.~\eqref{asymptotic-preservingadjointsystem2}.

\textbf{Microscopic lifting route}

Multiplying Eq.~\eqref{eq:app_macro_adjoint} by $\mathbf{\Psi}_k^T$ and adding the result to Eq.~\eqref{eq:app_scaled_micro} for each discrete velocity point $k$ gives
\begin{equation}
  \label{eq:app_AP_micro}
  \begin{aligned}
  \frac{\lambda_{i,k}^{*,n-1}-\lambda_{i,k}^{*,n}}{\Delta t}&+\mathbf{\Psi}_k^T\frac{\mathbf{\Lambda}_i^{*,n-1}-\mathbf{\Lambda}_i^{*,n}}{\Delta t}\\
  &+\frac{1}{|\Omega_i|}\sum_j\frac{\partial \mathcal{R}_{j,k}}{\partial f_{i,k}}(\mathbf{\Psi}_k^T \mathbf{\Lambda}_j^{*,n}+\lambda_{j,k}^{*,n})+\mathbf{\Psi}_k^T\frac{1}{|\Omega_i|}\sum_j\sum_{k^\prime} w_{k^\prime}\frac{\partial \mathcal{R}_{j,k^\prime}}{\partial \mathbf{W}_i}(\mathbf{\Psi}_{k^\prime}^T \mathbf{\Lambda}_j^{*,n}+\lambda_{j,k^\prime}^{*,n})\\
  &+\frac{1}{|\Omega_i|}\left(\frac{1}{w_k}\frac{\partial \mathbf{q}_i}{\partial f_{i,k}}+\mathbf{\Psi}_k^T\frac{\partial \mathbf{q}_i}{\partial \mathbf{W}_i}\right)\sum_j\sum_{k^\prime}w_{k^\prime}\frac{\partial \mathcal{R}_{j,k^\prime}}{\partial \mathbf{q}_i}(\mathbf{\Psi}_k^T \mathbf{\Lambda}_j^{*,n}+\lambda_{j,k^\prime}^{*,n})\\
  &-\frac{\mathbf{\Psi}_k^T{\mathbf{\Lambda}}_i^{*,n-1}-{\lambda}_{i,k}^{*,n-1}}{\tau_i}\\
  &-\frac{1}{\tau_i}\left(\frac{1}{w_k}\frac{\partial \mathbf{q}_i}{\partial f_{i,k}}+\mathbf{\Psi}_k^T\frac{\partial \mathbf{q}_i}{\partial \mathbf{W}_i}\right)\sum_{k^\prime} w_{k^\prime}  \frac{\partial f_{i,k^\prime}^+}{\partial \mathbf{q}_i} \lambda_{i,k^\prime}^{*,n}\\
  &+\frac{1}{\tau_i}\mathbf{\Psi}_k^T\frac{\partial \tau_i}{\partial \mathbf{W}_i}\sum_{k^\prime} w_{k^\prime}\frac{f_{i,k^\prime}^+-f_{i,k^\prime}}{\tau_i}\lambda_{i,k^\prime}^{*,n}=0.
  \end{aligned}
\end{equation}
The collision feedback terms in Eq.~\eqref{eq:app_AP_micro} follow directly from adding $\mathbf{\Psi}_k^T$ times the macroscopic collision block in Eq.~\eqref{eq:app_macro_adjoint} to the microscopic collision term $\lambda_{i,k}^{*,n-1}/\tau_i$ in Eq.~\eqref{eq:app_scaled_micro}.
Introducing $\widetilde{\lambda}_{i,k}^{*,n}$, Eq.~\eqref{eq:app_AP_micro} can be rewritten as
\begin{equation}
  \label{eq:app_AP_micro_tilde}
  \begin{aligned}
  \frac{\widetilde{\lambda}_{i,k}^{*,n-1}-\widetilde{\lambda}_{i,k}^{*,n}}{\Delta t}&+\frac{1}{|\Omega_i|}\sum_j\frac{\partial \mathcal{R}_{j,k}}{\partial f_{i,k}}\widetilde{\lambda}_{j,k}^{*,n}+\mathbf{\Psi}_k^T\frac{1}{|\Omega_i|}\sum_j\sum_{k^\prime} w_{k^\prime}\frac{\partial \mathcal{R}_{j,k^\prime}}{\partial \mathbf{W}_i}\widetilde{\lambda}_{j,k^\prime}^{*,n}\\
  &+\frac{1}{|\Omega_i|}\left(\frac{1}{w_k}\frac{\partial \mathbf{q}_i}{\partial f_{i,k}}+\mathbf{\Psi}_k^T\frac{\partial \mathbf{q}_i}{\partial \mathbf{W}_i}\right)\sum_j\sum_{k^\prime}w_{k^\prime}\frac{\partial \mathcal{R}_{j,k^\prime}}{\partial \mathbf{q}_i}\widetilde{\lambda}_{j,k^\prime}^{*,n}\\
  &-\frac{\mathbf{\Psi}_k^T\widetilde{\mathbf{\Lambda}}_i^{*,n-1}-\widetilde{\lambda}_{i,k}^{*,n-1}}{\tau_i}\\
  &-\frac{1}{\tau_i}\left(\frac{1}{w_k}\frac{\partial \mathbf{q}_i}{\partial f_{i,k}}+\mathbf{\Psi}_k^T\frac{\partial \mathbf{q}_i}{\partial \mathbf{W}_i}\right)\sum_{k^\prime} w_{k^\prime}  \frac{\partial f_{i,k^\prime}^+}{\partial \mathbf{q}_i} \widetilde{\lambda}_{i,k^\prime}^{*,n}\\
  &+\frac{1}{\tau_i}\mathbf{\Psi}_k^T\frac{\partial \tau_i}{\partial \mathbf{W}_i}\sum_{k^\prime} w_{k^\prime}\frac{f_{i,k^\prime}^+-f_{i,k^\prime}}{\tau_i}\widetilde{\lambda}_{i,k^\prime}^{*,n}=0.
  \end{aligned}
\end{equation}

Under the steady-state backward marching adopted in the main text, the macroscopic time-derivative term is absorbed into $\widetilde{\mathbf{\Lambda}}_i^{*,n}$, and Eq.~\eqref{eq:app_AP_micro_tilde} reduces to Eq.~\eqref{eq:micro_coupled_W}.
The two appendix routes are algebraically equivalent in the same way as in the main text: multiplying Eq.~\eqref{eq:app_AP_micro} by $w_k\partial f_{i,k}^+/\partial \mathbf{W}_i$, summing over $k$, and adding Eq.~\eqref{eq:app_macro_adjoint} reproduces Eq.~\eqref{eq:app_AP_system}.
Therefore, although the direct primal derivation yields a sequentially solvable adjoint system, the same AP adjoint algorithm derived in the main text is recovered after AP reconstruction.

\section{Jacobi Matrices for the linear UGKS} \label{app:jacobi}
In the development of the Linearized Unified Gas-Kinetic Scheme (L-UGKS), the linearization of the Shakhov equilibrium state $f^+ = g + g^+$ involves complex dependencies on both the macroscopic conservative variables $\mathbf{W} = [\rho, \rho \mathbf{V}, \rho E]^T$ and the heat flux vector $\mathbf{q}$. To facilitate a fully implicit or linearized solver, the Jacobi matrices representing the partial derivatives of the discrete distribution functions with respect to the macroscopic variables must be derived.

To reduce memory usage, two reduced distribution functions $h$ and $b$ are introduced,
\begin{equation*}
  \begin{aligned}
 h(\mathbf{x},t,\mathbf{v})&=\int f\td \mathbf{\xi}, \\
 b(\mathbf{x},t,\mathbf{v})&=\int \mathbf{\xi}^2f\td\mathbf{\xi}.
  \end{aligned}
\end{equation*}
Multiplying Eq.~(\ref{BGK}) by $1$ and $\mathbf{\xi}^2$ and integrating over the inner degrees of freedom yields
\begin{equation*}
  \begin{aligned}
  \frac{\partial h}{\partial t} + \mathbf{v}\cdot\nabla h =\frac{h^+-h}{\tau},\\
  \frac{\partial b}{\partial t} + \mathbf{v}\cdot\nabla b =\frac{b^+-b}{\tau},
  \end{aligned}
\end{equation*}
where the reduced equilibrium distributions are
\begin{equation*}
  \begin{aligned}
 h^+&=H+H^+,\\
 b^+&=B+B^+.
  \end{aligned}
\end{equation*}
The corresponding reduced Maxwellian distribution $g$ becomes
\begin{equation*}
  \begin{aligned}
H&=\int g\td\mathbf{\xi}=\rho\left(\frac{\lambda}{\pi}\right)^{D/2}e^{-\lambda(\mathbf{v}-\mathbf{v})^2},\\
B&=\int\mathbf{\xi}^2 g\td\mathbf{\xi}=\frac{3-D+K}{2\lambda}H,
  \end{aligned}
\end{equation*}
and the corresponding terms related to $g^+$ become
\begin{equation*}
  \begin{aligned}
 H^+&=\int g^+\td\mathbf{\xi} = \frac{4(1-Pr)\lambda^2}{5\rho}(\mathbf{v}-\mathbf{v})\cdot \mathbf{q}(2\lambda(\mathbf{v}-\mathbf{v})^2-2-D)H,\\
 B^+&=\int\mathbf{\xi}^2 g^+\td\mathbf{\xi} = \frac{4(1-Pr)\lambda^2}{5\rho}(\mathbf{v}-\mathbf{v})\cdot \mathbf{q}\{[ 2\lambda(\mathbf{v}-\mathbf{v})^2-D](3-D+K)-2K\}\frac{H}{2\lambda}.
  \end{aligned}
\end{equation*}

\subsection{Derivatives with Respect to Primitive Variables}
The primitive variable vector is defined as $\mathbf{Q} = [\rho, \mathbf{V}, \lambda]^T$, where $\lambda = \frac{m}{2k_B T}$. Taking the logarithm of the Maxwellian distribution functions and differentiating yields
\begin{equation*}
    \frac{\partial H}{\partial \mathbf{Q}} = H \left[ \frac{1}{\rho}, \quad 2\lambda(\mathbf{v}-\mathbf{V}), \quad \frac{D}{2\lambda} - (\mathbf{v}-\mathbf{V})^2 \right],
\end{equation*}
\begin{equation*}
    \frac{\partial B}{\partial \mathbf{Q}} = \frac{3-D+K}{2\lambda}\frac{\partial H}{\partial \mathbf{Q}} + H \left[ 0, \quad \mathbf{0}, \quad -\frac{3-D+K}{2\lambda^2} \right].
\end{equation*}

\subsection{Jacobian Matrix between Primitive and Conservative Variables}
These derivatives are transformed with respect to the conservative variable vector $\mathbf{W}$ by the chain rule through the transformation matrix $\frac{\partial \mathbf{Q}}{\partial \mathbf{W}}$. Let $\mathcal{C} = \frac{4\lambda^2}{(K+3)\rho}$. In a general three-dimensional velocity space context, this matrix is given by
\begin{equation*}
    \frac{\partial[\rho,\mathbf{V},\lambda]^T}{\partial[\rho,\rho \mathbf{V},\rho E]^T} = \begin{bmatrix}
        1 & \mathbf{0} & 0 \\
        -\mathbf{V}/\rho & \mathbf{I}/\rho & \mathbf{0} \\
        \frac{\lambda}{\rho}-\frac{1}{2}\mathcal{C}\mathbf{V}^2 & \mathcal{C}\mathbf{V} & -\mathcal{C}
    \end{bmatrix}.
\end{equation*}

\subsection{Shakhov Equilibrium Correction Term Derivatives}
Under the Shakhov model, the corrections $H^+$ and $B^+$ account for the correct Prandtl number ($\text{Pr}$). Let $\mathbf{c} = \mathbf{v}-\mathbf{V}$ denote the peculiar velocity vector. The derivatives of the Shakhov correction terms with respect to heat flux $\mathbf{q}$ and primitive variables are written as follows.

\subsubsection{For the $H^+$ Correction:}
\begin{equation*}
    \frac{\partial H^+}{\partial \mathbf{q}} = H^+ \frac{\mathbf{c}}{\mathbf{c}\cdot \mathbf{q}},
\end{equation*}
\begin{equation*}
    \frac{\partial H^+}{\partial \rho} = \frac{H^+}{H}\frac{\partial H}{\partial \rho} - \frac{H^+}{\rho},
\end{equation*}
\begin{equation*}
    \frac{\partial H^+}{\partial \mathbf{V}} = \frac{H^+}{H}\frac{\partial H}{\partial \mathbf{V}} - \frac{4\lambda \mathbf{c}}{2\lambda c^2 - 2 - D}H^+,
\end{equation*}
\begin{equation*}
    \frac{\partial H^+}{\partial \lambda} = \frac{H^+}{H}\frac{\partial H}{\partial \lambda} + \left( \frac{2c^2}{2\lambda c^2 - 2 - D} + \frac{2}{\lambda} \right) H^+.
\end{equation*}

\subsubsection{For the $B^+$ Correction:}
\begin{equation*}
    \frac{\partial B^+}{\partial \mathbf{q}} = B^+ \frac{\mathbf{c}}{\mathbf{c}\cdot \mathbf{q}},
\end{equation*}
\begin{equation*}
    \frac{\partial B^+}{\partial \rho} = \frac{B^+}{H}\frac{\partial H}{\partial \rho} - \frac{B^+}{\rho},
\end{equation*}
\begin{equation*}
    \frac{\partial B^+}{\partial \mathbf{V}} = \frac{B^+}{H}\frac{\partial H}{\partial \mathbf{V}}  - \frac{4\lambda \mathbf{c}(3-D+K)}{(2\lambda c^2 - D)(3-D+K) - 2K}B^+,
\end{equation*}
\begin{equation*}
    \frac{\partial B^+}{\partial \lambda} = \frac{B^+}{H}\frac{\partial H}{\partial \lambda} + \left( \frac{2c^2(3-D+K)}{(2\lambda c^2 - D)(3-D+K) - 2K} + \frac{1}{\lambda} \right) B^+.
\end{equation*}

\subsection{Expansion Matrix $M$ for Second-Order Space-Time Derivatives}
In the second-order UGKS flux evaluation, spatial and temporal derivatives of the Maxwellian state are expanded using a relation vector $\mathbf{a} = [a_1, a_2, a_3, a_4, a_5]^T$. The coefficient vector $\mathbf{a}$ relates directly to the spatial gradients of the conservative variables $\partial \mathbf{W}$ via the matrix $M$, structured as $\rho \mathbf{a} = M \partial \mathbf{W}$.
The matrix $M$ and its derivatives below are symmetric; for brevity, only the independent entries are written out, and the symbols $\ast$ denote the entries recovered by symmetry, i.e.\ $M_{ij}=M_{ji}$.
\begin{equation*}
    M = \begin{bmatrix}
        \dfrac{K+5}{2}+\dfrac{2\lambda^2\mathbf{V}^4}{K+3} & \ast & \ast & \ast & \ast \\
        -\dfrac{4\lambda^2U \mathbf{V}^2}{3 + K} & 2\lambda +\dfrac{8\lambda^2 U^2}{3 + K} & \dfrac{8\lambda^2 UV}{3 + K} & \dfrac{8 \lambda^2 UW}{3 + K} & -\dfrac{8 \lambda^2U}{3 + K} \\
        -\dfrac{4\lambda^2V \mathbf{V}^2}{3 + K} & \dfrac{8\lambda^2 UV}{3 + K} & 2\lambda + \dfrac{8\lambda^2 V^2}{3 + K} & \dfrac{8 \lambda^2 VW}{3 + K} & -\dfrac{8 \lambda^2V}{3 + K} \\
        -\dfrac{4\lambda^2W \mathbf{V}^2}{3 + K} & \dfrac{8\lambda^2 UW}{3 + K} & \dfrac{8\lambda^2 VW}{3 + K} & 2\lambda +\dfrac{8 \lambda^2 W^2}{3 + K} & -\dfrac{8 \lambda^2W}{3 + K} \\
        -2\lambda+\dfrac{4\lambda^2 \mathbf{V}^2}{3 + K} & -\dfrac{8\lambda^2 U}{3 + K} & -\dfrac{8\lambda^2 V}{3 + K} & -\dfrac{8 \lambda^2 W}{3 + K} & \dfrac{8 \lambda^2}{3 + K}
    \end{bmatrix}
\end{equation*}
The analytical sensitivities of $M$ with respect to the primitive parameters $\lambda, \mathbf{V}=(U, V, W)$ are
\begin{equation*}
    \frac{\partial M}{\partial \lambda} = \begin{bmatrix}
        \dfrac{4\lambda\mathbf{V}^4}{K+3} & \ast & \ast & \ast & \ast \\
        -\dfrac{8\lambda U \mathbf{V}^2}{3 + K} & 2 +\dfrac{16\lambda U^2}{3 + K} & \dfrac{16\lambda UV}{3 + K} & \dfrac{16 \lambda UW}{3 + K} & -\dfrac{16 \lambda U}{3 + K} \\
        -\dfrac{8\lambda V \mathbf{V}^2}{3 + K} & \dfrac{16\lambda UV}{3 + K} & 2+ \dfrac{16\lambda V^2}{3 + K} & \dfrac{16 \lambda VW}{3 + K} & -\dfrac{16 \lambda V}{3 + K} \\
        -\dfrac{8\lambda W \mathbf{V}^2}{3 + K} & \dfrac{16\lambda UW}{3 + K} & \dfrac{16\lambda VW}{3 + K} & 2 +\dfrac{16 \lambda W^2}{3 + K} & -\dfrac{16 \lambda W}{3 + K} \\
        -2+\dfrac{8\lambda \mathbf{V}^2}{3 + K} & -\dfrac{16\lambda U}{3 + K} & -\dfrac{16\lambda V}{3 + K} & -\dfrac{16 \lambda W}{3 + K} & \dfrac{16 \lambda}{3 +K}
    \end{bmatrix}
\end{equation*}

\begin{equation*}
    \frac{\partial M}{\partial U} = \begin{bmatrix}
        \dfrac{8\lambda^2U\mathbf{V}^2}{K+3} & \ast & \ast & \ast & \ast \\
        -\dfrac{4\lambda^2(3U^2 + V^2 + W^2)}{3 + K} & \dfrac{16\lambda^2 U}{3 + K} & \dfrac{8\lambda^2 V}{3 + K} & \dfrac{8 \lambda^2 W}{3 + K} & -\dfrac{8 \lambda^2}{3 + K} \\
        -\dfrac{8\lambda^2 UV}{3 + K} & \dfrac{8\lambda^2 V}{3 + K} & 0 & 0 & 0 \\
        -\dfrac{8\lambda^2UW }{3 + K} & \dfrac{8\lambda^2 W}{3 + K} & 0 & 0 & 0 \\
        \dfrac{8\lambda^2 U}{3 + K} & -\dfrac{8\lambda^2 }{3 + K} & 0 & 0 & 0
    \end{bmatrix}
\end{equation*}

\begin{equation*}
    \frac{\partial M}{\partial V} = \begin{bmatrix}
        \dfrac{8\lambda^2V\mathbf{V}^2}{K+3} & \ast & \ast & \ast & \ast \\
        -\dfrac{8\lambda^2 UV}{3 + K} & 0 & \dfrac{8\lambda^2 U}{3 + K} & 0 & 0 \\
        -\dfrac{4\lambda^2(U^2 + 3V^2 + W^2)}{3 + K} & \dfrac{8\lambda^2 U}{3 + K} & \dfrac{16\lambda^2 V}{3 + K} & \dfrac{8 \lambda^2 W}{3 + K} & -\dfrac{8 \lambda^2}{3 + K} \\
        -\dfrac{8\lambda^2VW }{3 + K} & 0 & \dfrac{8\lambda^2 W}{3 + K} & 0 & 0 \\
        \dfrac{8\lambda^2 V}{3 + K} & 0 & -\dfrac{8\lambda^2 }{3 + K} & 0 & 0
    \end{bmatrix}
\end{equation*}

\begin{equation*}
    \frac{\partial M}{\partial W} = \begin{bmatrix}
        \dfrac{8\lambda^2W\mathbf{V}^2}{K+3} & \ast & \ast & \ast & \ast \\
        -\dfrac{8\lambda^2 UW}{3 + K} & 0 & 0 & \dfrac{8\lambda^2 U}{3 + K} & 0 \\
        -\dfrac{8\lambda^2VW }{3 + K} & 0 & 0 & \dfrac{8\lambda^2 V}{3 + K} & 0 \\
        -\dfrac{4\lambda^2(U^2 + V^2 + 3W^2)}{3 + K} & \dfrac{8\lambda^2 U}{3 + K} & \dfrac{8\lambda^2 V}{3 + K} & \dfrac{16 \lambda^2 W}{3 + K} & -\dfrac{8 \lambda^2}{3 + K} \\
        \dfrac{8\lambda^2 W}{3 + K} & 0 & 0 & -\dfrac{8\lambda^2 }{3 + K} & 0
    \end{bmatrix}
\end{equation*}

\section{Linear unified gas-kinetic scheme} \label{app:linear_UGKS}

The linearized unified gas-kinetic scheme (L-UGKS) is constructed by taking first-order perturbations of the discrete residual operators about a converged base state $(\mathbf{W}_i,f_{i,k})$.
It provides an independent reference for verifying the adjoint sensitivities reported in Sec.~4.
The perturbation of the heat flux is
\begin{equation*}
\delta\mathbf{q}_i = \sum_{k^\prime} w_{k^\prime}\,\frac{\partial \mathbf{q}_i}{\partial f_{i,k^\prime}}\,\delta f_{i,k^\prime}+\frac{\partial \mathbf{q}_i}{\partial \mathbf{W}_i}\delta \mathbf{W}_i.
\end{equation*}
Consistent with the compatibility-based formulation in the main text, the linearized macroscopic constraint is
\begin{equation}
  \label{eq:macro_lin_compact}
  \delta\mathbf{W}_i- \sum_k w_k \delta f_{i,k}\mathbf{\Psi}_k = 0.
\end{equation}
The linearization of the discrete microscopic governing equation~\eqref{UGKSmicro} reads
\begin{equation}
  \label{eq:L_micro_lin}
  \begin{aligned} \frac{\delta f_{i,k}^{n+1}-\delta f_{i,k}^{n}}{\Delta t}
  &+ \frac{1}{|\Omega_i|}\sum_j
  \left(
  \frac{\partial \mathcal{R}_{i,k}}{\partial \mathbf{W}_j}\,\delta\mathbf{W}_j^n
  + \frac{\partial \mathcal{R}_{i,k}}{\partial f_{j,k}}\,\delta f_{j,k}^n
  + \frac{\partial \mathcal{R}_{i,k}}{\partial \mathbf{q}_j^n}\,\delta\mathbf{q}_j
  \right) \\
  &- \frac{1}{\tau_i}\left(
  \frac{\partial f_{i,k}^+}{\partial \mathbf{W}_i}\,\delta\mathbf{W}_i^{n+1}
  + \frac{\partial f_{i,k}^+}{\partial \mathbf{q}_i}\,\delta\mathbf{q}_i^n
  - \delta f_{i,k}^{n+1}
  \right) \\
  &+ \frac{1}{\tau_i}\,\frac{\partial \tau_i}{\partial \mathbf{W}_i} \delta\mathbf{W}_i^{n+1}
  \frac{f_{i,k}^+-f_{i,k}}{\tau_i}\,
  = 0.
  \end{aligned}
\end{equation}
For the steady-state problems considered in this work, the time indices associated with the base-flow variables are omitted below.

Multiplying Eq.~\eqref{eq:L_micro_lin} by $w_k\mathbf{\Psi}_k$ and summing over all discrete velocity points, and using the identity
\begin{equation*}
  \sum_{k} w_k \mathbf{\Psi}_k \frac{\partial f_{i,k}^+}{\partial \mathbf{W}_i}
  = \frac{\partial}{\partial \mathbf{W}_i}\sum_{k} w_k \mathbf{\Psi}_k f_{i,k}^+
  = \mathbf{I},
\end{equation*}
one recovers the linearized macroscopic evolution
\begin{equation}
  \label{eq:ap_linear_macro}
  \begin{aligned}
  \frac{\delta \mathbf{W}_i^{n+1}-\delta\mathbf{W}_i^{n}}{\Delta t}
  &+ \frac{1}{|\Omega_i|}\sum_k w_k\mathbf{\Psi}_k\sum_j
  \left(
  \frac{\partial \mathcal{R}_{i,k}}{\partial \mathbf{W}_j}\,\delta\mathbf{W}_j^n
  + \frac{\partial \mathcal{R}_{i,k}}{\partial f_{j,k}}\,\delta f_{j,k}^n
  + \frac{\partial \mathcal{R}_{i,k}}{\partial \mathbf{q}_j^n}\,\delta\mathbf{q}_j
\right)  = 0.
  \end{aligned}
\end{equation}
Thus the L-UGKS system is composed of Eqs.~\eqref{eq:ap_linear_macro} and~\eqref{eq:L_micro_lin}, subject to the linearized compatibility constraint~\eqref{eq:macro_lin_compact}.

\subsection{Residual evaluation}
\begin{figure}[ht]
  \centering
  \includegraphics[width=0.5\textwidth]{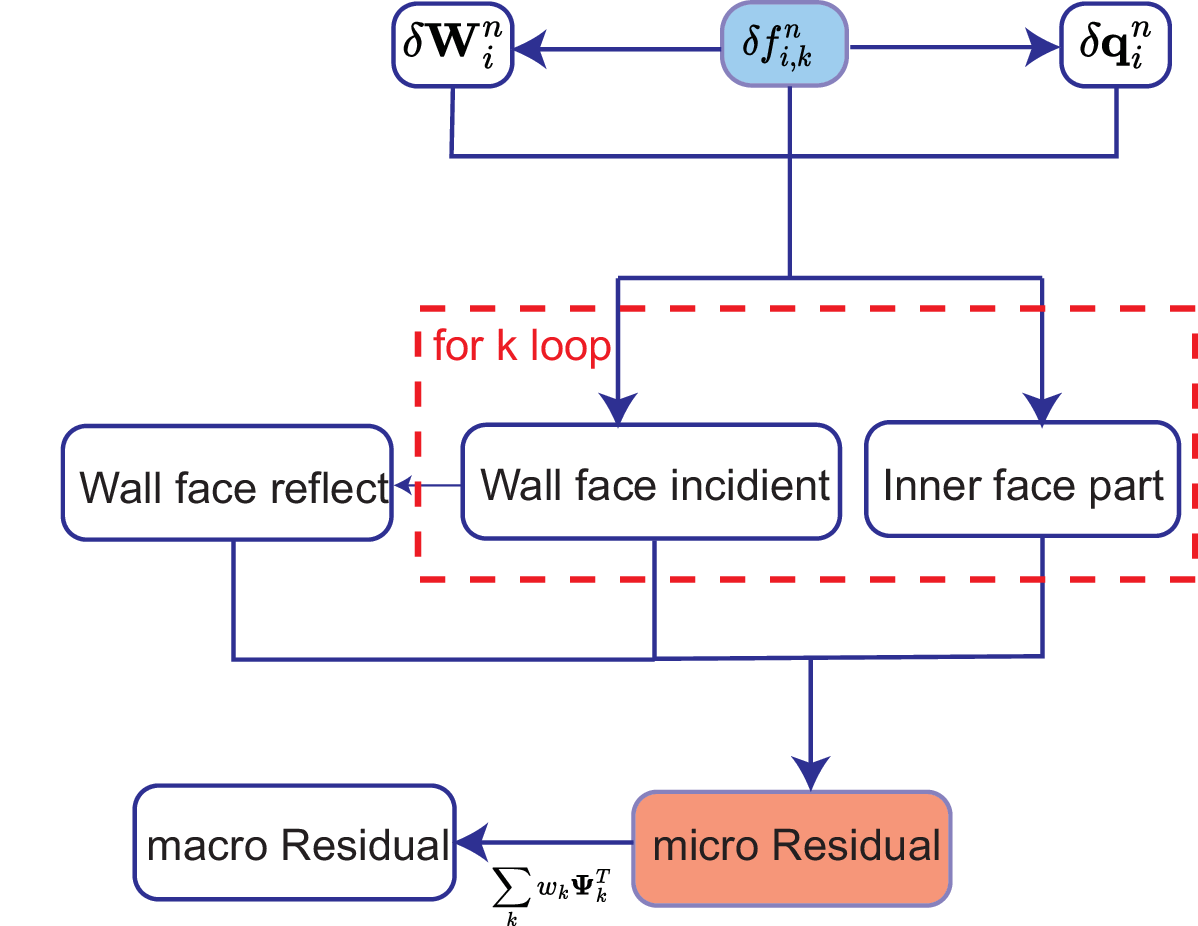}
  \caption{Organization of linear residual evaluation following the UGKS programming paradigm~\cite{zhang2025efficiency}. The procedure is reversed to the adjoint residual evaluation in Fig.~\ref{fig:resAlg}.}
  \label{fig:linearResAlg}
\end{figure}

The linear residual is evaluated with the same UGKS programming paradigm~\cite{zhang2025efficiency} as in the forward solver.
The common contributions associated with $\mathbf{W}_i$ and $\mathbf{q}_i$ are precomputed, after which the velocity-space loop assembles the pointwise microscopic residual.
This ordering is the reverse of the adjoint residual evaluation in Sec.~3 (Fig.~\ref{fig:resAlg}): the adjoint residual accumulates the common part first within the velocity loop and then adds the local $f_{i,k}$ contribution, whereas the linear residual follows the forward sequence.
The overall organization of the linear residual evaluation is summarized in Fig.~\ref{fig:linearResAlg}.

\subsection{Implicit marching}
The L-UGKS time advancement employs the same implicit predictor--corrector strategy as the adjoint solver~\cite{zhu2016implicit}.
At each step, the macroscopic equation~\eqref{eq:ap_linear_macro} is solved first by LU-SGS to obtain a predicted macroscopic perturbation $\delta\mathbf{W}_i^{n+1/2}$; the microscopic equation~\eqref{eq:L_micro_lin} is then solved for $\delta f_{i,k}^{n+1}$; a final correction updates $\delta\mathbf{W}_i^{n+1}$ through the linearized compatibility condition~\eqref{eq:macro_lin_compact}.
Relative to the adjoint formulation, the neighbor coupling uses the primal (non-transposed) flux Jacobian and the forward cell-sweep order, while the migration--relaxation splitting and the treatment of the stiff collision terms remain analogous.



\end{document}